\documentclass[10pt, reqno]{amsart}
\usepackage[a4paper, margin=3.1cm]{geometry}
\usepackage[english]{babel}
\usepackage{microtype}
\usepackage{graphicx}
\usepackage{amsfonts}
\usepackage{amsmath}
\usepackage{amsthm}
\usepackage{setspace}
\usepackage{amssymb}
\usepackage{pst-node}
\usepackage{mathrsfs}
\usepackage{tikz}
\usepackage{pgfplots}
\pgfplotsset{compat=1.18}
\usetikzlibrary{arrows.meta, positioning, decorations.pathmorphing}
\usetikzlibrary{fit,calc,decorations.pathreplacing}
\usepackage{hyperref}
\colorlet{darkblue}{blue!50!black}
\hypersetup{
	colorlinks,%
	citecolor=darkblue,%
	filecolor=red,%
	linkcolor=darkblue,%
	urlcolor=darkblue,%
	pdfnewwindow=true,%
	pdfstartview={FitH}
}
\usepackage{comment}
\usepackage{tabularx}
\usepackage{pifont}
\usepackage{bbm}
\usepackage{appendix}
\usetikzlibrary{automata,arrows,positioning,calc}
\usepackage{accents}
\usepackage{enumitem}

\usepackage[framemethod=TikZ]{mdframed}
\usepackage{thmtools}
\usepackage{xcolor}
\usepackage{cleveref}
\usepackage{booktabs}
\usepackage{makecell}
\usepackage{stmaryrd}
\usepackage{caption}

\usetikzlibrary{patterns}

\definecolor{convexdiffusive}{HTML}{B9D8EA}
\definecolor{nonconvexdiffusive}{HTML}{F5D28A}
\definecolor{nonconvexsuperdiffusive}{HTML}{B8DDCF}
\definecolor{convexsuperdiffusive}{HTML}{E5C1D7}

\makeatletter
\renewcommand{\section}{\@startsection{section}{1}{\z@}%
  {-35\p@ \@plus -4\p@ \@minus -4\p@}%
  {26\p@ \@plus 4\p@ \@minus 4\p@}%
  {\normalfont\LARGE\bfseries\centering}}

\renewcommand{\subsection}{\@startsection{subsection}{2}{\z@}%
  {-22\p@ \@plus -4\p@ \@minus -4\p@}%
  {15\p@ \@plus 4\p@ \@minus 4\p@}%
  {\normalfont\Large\bfseries}}
  
\renewcommand{\subsubsection}{\@startsection{subsubsection}{3}{\z@}%
  {-14\p@ \@plus -3\p@ \@minus -2\p@}%
  {7\p@ \@plus 2\p@ \@minus 2\p@}%
  {\normalfont\large\bfseries}}

\renewcommand{\@seccntformat}[1]{\textbf{\csname the#1\endcsname}\quad}
\renewcommand{\sectionmark}[1]{\markboth{#1}{\@title}}
\makeatother

\colorlet{myblue}{blue!50!white}
\colorlet{myblue2}{blue!30!white}

\theoremstyle{plain}
\newtheorem{theorem}{Theorem}[section]
\newtheorem{lemma}[theorem]{Lemma}
\newtheorem{proposition}[theorem]{Proposition}

\newtheorem{corollary}[theorem]{Corollary}

\theoremstyle{definition}
\newtheorem{definition}[theorem]{Definition}
\newtheorem{example}[theorem]{Example}
\newtheorem{remark}[theorem]{Remark}

\newtheorem*{lemma*}{Lemma}

\declaretheorem[name=Condition]{hyp}

\numberwithin{equation}{section}

\newcommand{\NN}{\mathbb{N}}
\newcommand{\ZZ}{\mathbb{Z}}

\newcommand{\RR}{\mathbb{R}}

\newcommand{\PP}{\mathbb{P}}
\newcommand{\EE}{\mathbb{E}}

\newcommand{\dd}{\text{d}}
\newcommand{\eqdef}{ \overset{\textup{\tiny def.}}{=}}

\newcommand{\dbl}{\llbracket}
\newcommand{\dbr}{\rrbracket}

\newcommand{\notation}[4]{#1 & #2 & #3 & #4 \\}

\title{Large Deviations and Convexity Transition for the Multidimensional Elephant Random Walk}

\author{Jan-Luka Fatras}

\makeatletter
\let\oldsetauthors\@setauthors
\renewcommand{\@setauthors}{%
    \oldsetauthors
    \vspace{0.8em}
    \begin{center}
        \small
        Institut de Mathématiques de Toulouse, Université de Toulouse, \\
        CNRS, UPS, F-31062 Toulouse Cedex 9, France. \\
        \vspace{0.4 em}
        \texttt{jan-luka.fatras@math.univ-toulouse.fr}
    \end{center}
    \vspace{-0.1cm}
 
}
\makeatother



\begin{document}

\begin{abstract}
\noindent In this article we establish a weak large deviation principle for the empirical measures of the elephant random walk on $\ZZ^d$, for every reinforcement parameter $p \in (0,1)$. We identify a qualitative change in the behaviour of the associated rate function between dimension $d=1$ and dimensions $d \geq 2$. In dimension $d=1$, the rate function coincides with the Donsker--Varadhan rate function of the simple symmetric random walk and is independent of the reinforcement parameter $p$. By contrast, in every dimension $d\geq 2$, the rate function differs from that of the simple random walk as soon as $p\neq 1/(2d)$. Finally, we prove that for all $d \geq 2$ the rate function is convex if and only if $p \leq 1/(2d)$. This allows us to identify a critical point in the large-deviation behaviour of the model which occurs strictly before the standard diffusive-to-superdiffusive crossover.

\vspace{-0.5cm}
\end{abstract}

\maketitle

\markboth{}{}

\tableofcontents

\section{Introduction}\label{section_introduction}

The study of large deviation principles (LDPs) for empirical measures originates with Sanov's theorem \cite{Sanov} for independent and identically distributed random variables, which was subsequently extended to more general spaces in \cite{Lanford} and \cite{Bahadur}. The transition to dependent processes began with the seminal work of Donsker and Varadhan \cite{Donsker_varadhan,Donsker_varadhan_2,Donsker_varadhan_3,Donsker_varadhan_4}, who established the variational characterisation of the rate function for Markov chains. A classic example for which this theory yields a weak LDP is that of the symmetric random walk on $\ZZ^d$. A natural question is whether such results extend beyond the Markovian setting. 

Reinforced random walks provide a natural class of such models, where memory effects lead to qualitatively different behaviour. A canonical example is the elephant random walk (ERW), introduced in \cite{Schutz_Trimper}, in which each step depends on the entire past through a reinforcement parameter $p \in (0,1)$. This long-range memory induces a crossover as the memory parameter varies, separating diffusive and superdiffusive regimes. The corresponding limit theorems have been established in several directions, including laws of large numbers, central limit theorems and functional limit theorems \cite{Baur_Bertoin,Coletti_Gava_Schutz,Bercu_2018,Kubota_Naoki,Guerin_Laulin_Raschel}, moderate deviations and Berry--Esseen bounds~\cite{Fan_2021}, and recurrence and transience criteria~\cite{Coletti_2021}.  Multidimensional versions of the model have also been introduced and studied at the level of scaling limits~\cite{Lucile_multi_dim_ERW}, functional limit theorems~\cite{Bertenghi} and recurrence criteria~\cite{Curien_Laulin_recurrence,Shuo_Qin_MERW}.

Regarding large deviations, memory effects have been shown to alter the asymptotic behaviour of rare events, see \cite{Touchette_Harris_speed,Robert_Harris_giant_leaps}. Heuristically, this phenomenon stems from the fact that reinforcement can make certain atypical trajectories more or less costly, thereby modifying the likelihood of large fluctuations. The impact of memory on the large deviations of the empirical measures of the ERW has, to the best of our knowledge, not been explored. Since the empirical measure encodes the time spent by the walk at different locations, it naturally captures the cumulative effect of reinforcement and thus provides a natural observable to probe the effect of memory.

Several large-deviation results are already available for observables related to the increments of the ERW. In dimension one, the Pólya-type urn representation recalled in \cite{Baur_Bertoin}, combined with the sample-path LDP of \cite{Franchini_urns}, yields an LDP for the rescaled urn trajectory and hence, by contraction, for the endpoint $S_n/n$. A different approach is provided by \cite{Budhiraja_Waterbury,Zoubouloglou}. Applied to the finite-state increment process, their results yield an LDP for the empirical measure of the increments, which again contracts to an LDP for $S_n/n$. More recently, González Rangel~\cite{Juan_Esaul_Gonzalez_Rangel} recovered the one-dimensional endpoint LDP from precise asymptotics of the Laplace transform via the Gärtner--Ellis theorem. Two recent works, \cite{Fan_Liang_Ma_2026,Liang_2026}, further investigate one-dimensional variants with random step sizes (and, in the latter, stops), focusing on an LDP for the endpoint rather than empirical-measure large deviations.

None of these approaches directly addresses the empirical measure of the positions considered here. In the urn representation, the relevant path records the cumulative numbers of positive and negative increments, while the empirical measure of the walk depends on their microscopic ordering through the resulting local occupation times. Likewise, the empirical measure of the increments records only the frequencies of the different steps and does not determine the spatial occupation of the walk. The empirical-measure problem considered here is therefore genuinely spatial and requires a different approach.

In this article, we establish the existence of a weak LDP for the empirical measures of the multidimensional elephant random walk (MERW). In dimension $d=1$, we identify the rate function by showing that it coincides with the Donsker--Varadhan rate function of the symmetric simple random walk on $\ZZ$, for every value of the reinforcement parameter  $p$. In other words, the memory introduced by the reinforcement mechanism leaves no trace at the level of large deviations of the empirical measure in dimension one. The picture in dimension $d \geq 2$ is strikingly different: the rate function depends on $p$ and does not coincide with the Donsker--Varadhan rate function of the symmetric simple random walk on $\ZZ^d$, except when $p=1/(2d)$. Furthermore, the rate function is convex if and only if $p \leq 1/(2d)$. The non-convexity of the rate function in the self-reinforcing regime $p>1/(2d)$ indicates a fundamental shift in the large-deviation behaviour of the empirical measure. In addition, the transition from convexity to non-convexity occurs at the critical value $p=1/(2d)$, strictly before the standard diffusive-to-superdiffusive crossover, and is therefore distinct from it.

\medskip

\noindent Throughout the article, we write $\NN=\{1,2,\dots\}$ and $\NN_0=\NN \cup \{0\}$. Moreover, for $m,n \in \ZZ$, we set $\dbl m,n \dbr=\{m,m+1, \dots , n \}$ if $m \leq n$, and $\dbl m,n \dbr=\emptyset$ if $n<m$.

\subsection{The multidimensional Elephant Random Walk}\label{subsection_introducing_the_ERW}

We recall the definition of the multidimensional elephant random walk (MERW), following \cite[p.~1147]{Lucile_multi_dim_ERW}. Let $p \in (0,1)$ and $d \in \NN$. Denoting by $(e_1,\dots,e_d)$ the canonical basis of $\ZZ^d$, let
\begin{equation*}\label{eq_def_mathscr_S_d}
\mathscr{S}_d \eqdef \{\pm e_i\}_{i=1}^d
\end{equation*}
be the set of nearest-neighbour increments on $\ZZ^d$. Let $(\Omega,\mathcal{F},\PP_{d,p})$ be a probability space supporting the following construction. 

The MERW is defined as follows. The first increment is sampled according to $\nu_{\text{unif}}$, the uniform measure on $\mathscr{S}_d$. At each subsequent step $n \geq 2$, one of the previous $n-1$ increments is selected uniformly at random. With probability $p$, the walker repeats the selected increment. With probability $1-p$, it instead chooses uniformly at random one of the remaining $2d-1$ elements of $\mathscr{S}_d$.

More formally, let $(U_n)_{n \geq 1}$ be a sequence of independent random variables on $(\Omega,\mathcal{F})$, such that under $\PP_{d,p}$,
\begin{equation*}\label{eq_def_U_n}
U_n \sim \mathrm{Unif}(\dbl 1,n \dbr).
\end{equation*}
Let $(A_n)_{n \geq 2}$ be a sequence of i.i.d.~random matrices, independent of $(U_n)_{n \in \NN}$, taking values in $\{\pm \mathbb{I}_d, \pm J_d, 
\ldots, \pm J_d^{\,d-1}\}$ with distribution
\begin{align*}
\PP_{d,p}(A_n = +\mathbb{I}_d) = p, 
\qquad 
\PP_{d,p}(A_n = M) = \frac{1-p}{2d-1} \ \text{ for all } M \neq + \mathbb{I}_d.
\end{align*}
Here, $\mathbb{I}_d=J^0_d$ denotes the identity matrix and $J_d$ is the cyclic permutation matrix
\begin{align*}
J_d =
{\setlength{\arraycolsep}{3pt}
\renewcommand{\arraystretch}{0.8}
\begin{pmatrix}
0 & 1 & 0 & \cdots & 0 \\
0 & 0 & 1 & \ddots & \vdots \\
\vdots & \vdots & \ddots & \ddots & 0 \\
0 & 0 & \cdots & 0 & 1 \\
1 & 0 & \cdots & 0 & 0
\end{pmatrix}},
\end{align*}
which acts on $\mathscr{S}_d$ by $J_d^k e_i = e_{i-k}$, with indices taken modulo $d$. The MERW $(S_n)_{n \in \NN}$ on $\ZZ^d$ is defined recursively by $S_1 =0$ and, for every $n \geq 1$,
\begin{equation*}\label{eq_def_ERW}
S_{n+1}=S_n+X_{n},
\end{equation*}
where the sequence of increments $(X_n)_{n \in \NN}$ is given by
\begin{equation}\label{eq_def_ERW_increments}
X_1 \sim \nu_{\text{unif}} \qquad \text{and} \qquad X_n = A_n X_{U_{n-1}}, \quad n \geq 2,
\end{equation}
and the random variable $X_1$ and the families $(A_n)_{n \geq 2}$ and $(U_n)_{n \in \NN}$ are mutually independent. This model is the natural $d$-dimensional analogue of the one-dimensional ERW of \cite{Schutz_Trimper}. Note that, with this indexing convention, the process has only $n-1$ increments by time $n$. Instead of considering the initial step to be distributed according to the uniform measure $\nu_{\text{unif}}$, one can more generally define the MERW with an arbitrary initial-step distribution $\nu$ on $\mathscr{S}_d$. In Section~\ref{subsection_invariance_of_rate_function}, we show that this choice has no impact at the level of the large deviations of the empirical measures. For the sake of convenience, we therefore choose $\nu_{\text{unif}}$ for the initial step distribution.

\medskip

The MERW is a model in which the walker retains information about its past. Nevertheless, in dimension $d=1$, the ERW is a Markov chain, albeit a time-inhomogeneous one \cite[p.4]{Coletti_2021}. This stems from the fact that the position at time $n$ fully determines the proportions of upward and downward steps. This feature is completely lost in higher dimension. When $d \geq 2$ and $p \neq 1/(2d)$, the MERW is non-Markovian \cite[p.1059]{Shuo_Qin_MERW}.

It was shown in \cite{Schutz_Trimper} and later in \cite{Baur_Bertoin,Coletti_Gava_Schutz,Bercu_2018} that the asymptotic behaviour of the one-dimensional ERW depends on the value of the reinforcement parameter $p$. This work has been extended to higher dimensions by \cite{Lucile_multi_dim_ERW}. Consider 
\begin{equation}\label{eq_def_diffusive_to_superdiff_crit_param}
p_d \eqdef \frac{2d+1}{4d},
\end{equation}
which we will refer to as the \textit{diffusive-to-superdiffusive critical parameter}. In \cite{Lucile_multi_dim_ERW}, it is shown that
\begin{itemize}
\item[-] When $0<p<p_d$, the MERW is \textit{diffusive}:
\begin{align*}
\frac{S_n}{\sqrt{n}} \xrightarrow[n \to \infty]{\mathcal{L}} \mathcal{N}\left( 0 , \frac{2d-1}{d[2d(1-2p)+1]} \mathbb{I}_d \right).
\end{align*}

\item[-] When $p=p_d$, the MERW is \textit{critical}: 
\begin{align*}
\frac{S_n}{\sqrt{n \log n}} \xrightarrow[n \to \infty]{\mathcal{L}} \mathcal{N}\left( 0 , \frac{1}{d} \mathbb{I}_d \right).
\end{align*}

\item[-] For $p>p_d$, the MERW is \textit{superdiffusive}. There exists a non-degenerate random vector $L$, which does not collapse to a deterministic constant, such that
\begin{align*}
\frac{S_n}{n^{\frac{2dp-1}{2d-1}}} \xrightarrow[n \to \infty]{\text{a.s.}} L.
\end{align*}
\end{itemize}
This critical value also governs the recurrence--transience transition in dimensions one and two. More precisely, Coletti and Papageorgiou proved in \cite[Theorems~3.1 and 3.3]{Coletti_2021} that, for $d=1$, the ERW is recurrent when $p \leq 3/4 = p_1$ and transient when $p>p_1$. In dimension $d=2$, the MERW is recurrent when $p<5/8 = p_2$ and transient when $p \geq p_2$, as established in \cite[Theorem~1.1]{Shuo_Qin_MERW}. An alternative proof of recurrence in the diffusive regime is given in \cite[Theorem~1]{Curien_Laulin_recurrence}. By contrast, in dimensions $d \geq 3$, the MERW is transient for every $p \in (0,1)$, see \cite[Theorem~1.2]{Shuo_Qin_MERW}.

\subsection{Main results}\label{subsection_main_results}

We denote by $\mathcal{P}(\ZZ^d)$ the space of probability measures on $\ZZ^d$, endowed with the topology of weak convergence and the associated Borel $\sigma$-algebra $\mathcal{B}(\mathcal{P}(\ZZ^d))$.

We are interested in the sequence of empirical measures $(\ell_n)_{n \in \NN}$ associated with the MERW, defined by
\begin{equation}\label{eq_def_empirical_measure}
\ell_n \eqdef \frac{1}{n} \sum_{j=1}^n \delta_{S_j} , \quad n \in \NN.
\end{equation}
The sequence $(\ell_n)_{n \in \NN}$ is a sequence of $\mathcal{P}(\ZZ^d)$-valued random variables. Our aim is to establish a weak large deviation principle (LDP) for this sequence.

We recall that a sequence of random variables $(X_n)_{n \in \NN}$ on a topological space $\mathcal{X}$, endowed with its Borel $\sigma$-algebra, satisfies a (full) LDP with rate function $I:\mathcal{X} \to [0,+\infty]$, where $I$ is lower semicontinuous, if 
\begin{align}
-\inf \limits_{x \in U} I(x) \leq \liminf \limits_{n \to \infty} \frac{1}{n} \log \PP \left( X_n \in U \right) & \qquad \text{for all open sets } U \subseteq \mathcal{X} \, , \label{eq_def_LDP_lower_bound} \\
-\inf \limits_{x \in F} I(x) \geq \limsup \limits_{n \to \infty} \frac{1}{n} \log \PP \left( X_n \in F \right)   & \qquad  \text{for all closed sets } F \subseteq \mathcal{X}. \label{eq_def_LDP_upper_bound}
\end{align}
If the lower bound~(\ref{eq_def_LDP_lower_bound}) holds for all open sets and the upper bound~(\ref{eq_def_LDP_upper_bound}) holds only for compact sets, rather than for all closed sets, we say that the LDP is \textit{weak} rather than full.

For each dimension $d \in \NN$ and memory parameter $p \in (0,1)$, define the candidate rate function $I_{d,p}:\mathcal{P}(\ZZ^d) \longrightarrow [0,+\infty]$ by 	  
\begin{equation}\label{expression_rate_function}
I_{d,p}(\mu)= -\inf_{\varepsilon>0} \limsup_{n \to \infty} \frac{1}{n} \log \PP_{d,p}\big( \ell_n \in B(\mu,\varepsilon) \big), \qquad \mu \in \mathcal{P}(\ZZ^d),
\end{equation}
where $B(\mu,\varepsilon)$ denotes the open ball of radius $\varepsilon>0$ centred at $\mu$ with respect to the metric on $\mathcal{P}(\ZZ^d)$ introduced in Section~\ref{subsection_topology_space_of_proba_measures}, which metrises the topology of weak convergence. The main results of this article are stated below.

\begin{theorem}[Existence of a weak LDP]\label{theorem_existence_LDP}
Let $d \in \NN$ and $p \in (0,1)$. Then, under $\PP_{d,p}$, the sequence $(\ell_n)_{n \in \NN}$ satisfies a weak large deviation principle on $\mathcal{P}(\ZZ^d)$, endowed with the topology of weak convergence, with rate function $I_{d,p}$.
\end{theorem}

\begin{proposition}[Failure of the full LDP]\label{prop_no_full_LDP}
For every $d \in \NN$ and $p \in (0,1)$, the sequence of empirical measures $(\ell_n)_{n \in \NN}$ does not satisfy a full LDP on $\mathcal{P}(\ZZ^d)$ under $\PP_{d,p}$, with any rate function.
\end{proposition}

\noindent The proofs of Theorem~\ref{theorem_existence_LDP} and Proposition~\ref{prop_no_full_LDP} are given in Section~\ref{subsection_existence_of_a_weak_LDP_measures}. In Section~\ref{subsection_profile_and_path_observables}, we show that the method developed here can be adapted to establish an LDP for the cumulative step distribution, from which we derive LDPs for the rescaled endpoint $S_n/n$, as in \cite{Franchini_urns}, and for the maximal displacement of the MERW.

\medskip

The following proposition compares the rate function $I_{d,p}$ to that of the symmetric simple random walk. Denote by $I_{d,\textup{DV}}$ the usual Donsker--Varadhan rate function \cite{Donsker_varadhan} for the symmetric simple random walk on $\ZZ^d$, given by 
\begin{align}\label{Donsker_Varadhan_our_case}
\begin{array}{cccc}
I_{d,\textup{DV}}: & \mathcal{P}\big(\ZZ^d \big) & \longrightarrow & [0,+\infty] \\
& \mu & \longmapsto & \sup \limits_{\substack{(u_k)_{k \in \ZZ^d} \, : \\ 0< \inf u_k \leq \sup u_k < + \infty }} \displaystyle \sum \limits_{k \in \ZZ^d} \mu(k) \log \left( \frac{ u_k}{\frac{1}{2d}\sum_{s \in \mathscr{S}_d}u_{k+s}} \right) 
\end{array}.
\end{align}

\begin{proposition}[Dimensional dichotomy]\label{prop_dimensional_dichotomy}
    The rate function $I_{d,p}$ satisfies the following:
    \begin{itemize}
        \item[-] In dimension $d=1$, for every $p \in (0,1)$,
        \begin{equation}\label{eq_prop_dimensional_dichotomy_dim_1}
        I_{1,p}=I_{1,\textup{DV}}.
        \end{equation}
        \item[-] In dimensions $d \geq 2$, for every $p \in (0,1)$ such that $p \neq 1/(2d)$,
        \begin{equation}\label{eq_prop_dimensional_dichotomy_dim_geq_1}
        I_{d,p} \neq I_{d,\textup{DV}}.
        \end{equation}
    \end{itemize}
\end{proposition}

\noindent Proposition~\ref{prop_dimensional_dichotomy} highlights a dichotomy between dimension one and higher dimensions. On $\ZZ$, at the level of the large deviations for the empirical measure, the ERW and the simple symmetric random walk are equivalent. The effect of memory disappears. In particular, since \eqref{eq_prop_dimensional_dichotomy_dim_1} holds for all $p \in (0,1)$, the diffusive-to-superdiffusive crossover observed in \cite{Schutz_Trimper} is not visible at this level. In contrast, in dimensions $d \geq 2$, memory affects the rate function.

\begin{remark}\label{remark_sharpness_diff_I_I_DV}
When $p=1/(2d)$ and the initial step distribution is the uniform distribution over $\mathscr{S}_d$, a MERW reduces to a symmetric random walk, yielding $I_{d,1/(2d)}=I_{d,\text{DV}}$, see Lemma~\ref{lemma_MERW_p_half_is_SRW}. Thus, the condition $p \neq 1/(2d)$ in the second assertion of Proposition~\ref{prop_dimensional_dichotomy} is sharp. In particular, when $d \geq 2$, the rate function $I_{d,p}$ is not constant in the memory parameter $p$. In other words, memory plays a role in higher dimension. 
\end{remark}

\noindent The two assertions of Proposition~\ref{prop_dimensional_dichotomy} are proved in Sections~\ref{subsubsection_one_dimensional_rate_function} and \ref{subsubsection_quali_study_in_higher_dimension}, respectively.

\medskip

Finally, the following proposition shows that there is a qualitative change in the behaviour of the rate function when the memory parameter is large. 

\begin{proposition}[Change of regime in the convexity]\label{prop_change_of_regime}
For $d \geq 2$ and $p \in (0,1)$, we have
\begin{equation}
I_{d,p} \text{ is convex } \Longleftrightarrow p \in \left(0,\frac{1}{2d}\right].
\end{equation}
\end{proposition}

\noindent Proposition~\ref{prop_change_of_regime} shows that convexity of the rate function is monotone with respect to the memory parameter $p$. This motivates, for $d \geq 2$, the definition of the following convexity threshold 
\begin{equation}\label{eq_def_new_critical_point}
p_d^* \eqdef \inf \big\{ q \in (0,1) \, : \, I_{d,p} \text{ is non-convex for every } p \in [q,1) \big\}.
\end{equation}
The same proposition yields 
\begin{align*}
p_d^*=\frac{1}{2d} < \frac{2d+1}{4d} = p_d.
\end{align*}
Consequently, the loss of convexity of the large-deviation rate function occurs strictly before the diffusive-to-superdiffusive transition. In particular, non-convexity already arises in the diffusive regime. This strict inequality shows that $p_d^*$ differs from the classical critical parameter $p_d$, indicating that the change of regime identified here is distinct from the phase transition first observed in \cite{Schutz_Trimper} and subsequently extended to higher dimensions in \cite{Lucile_multi_dim_ERW}. This is illustrated in the phase diagram of Figure~\ref{fig_phase_diagram} below. The proof of Proposition~\ref{prop_change_of_regime} is given in Section~\ref{subsection_phase_transition_for_MERW}.

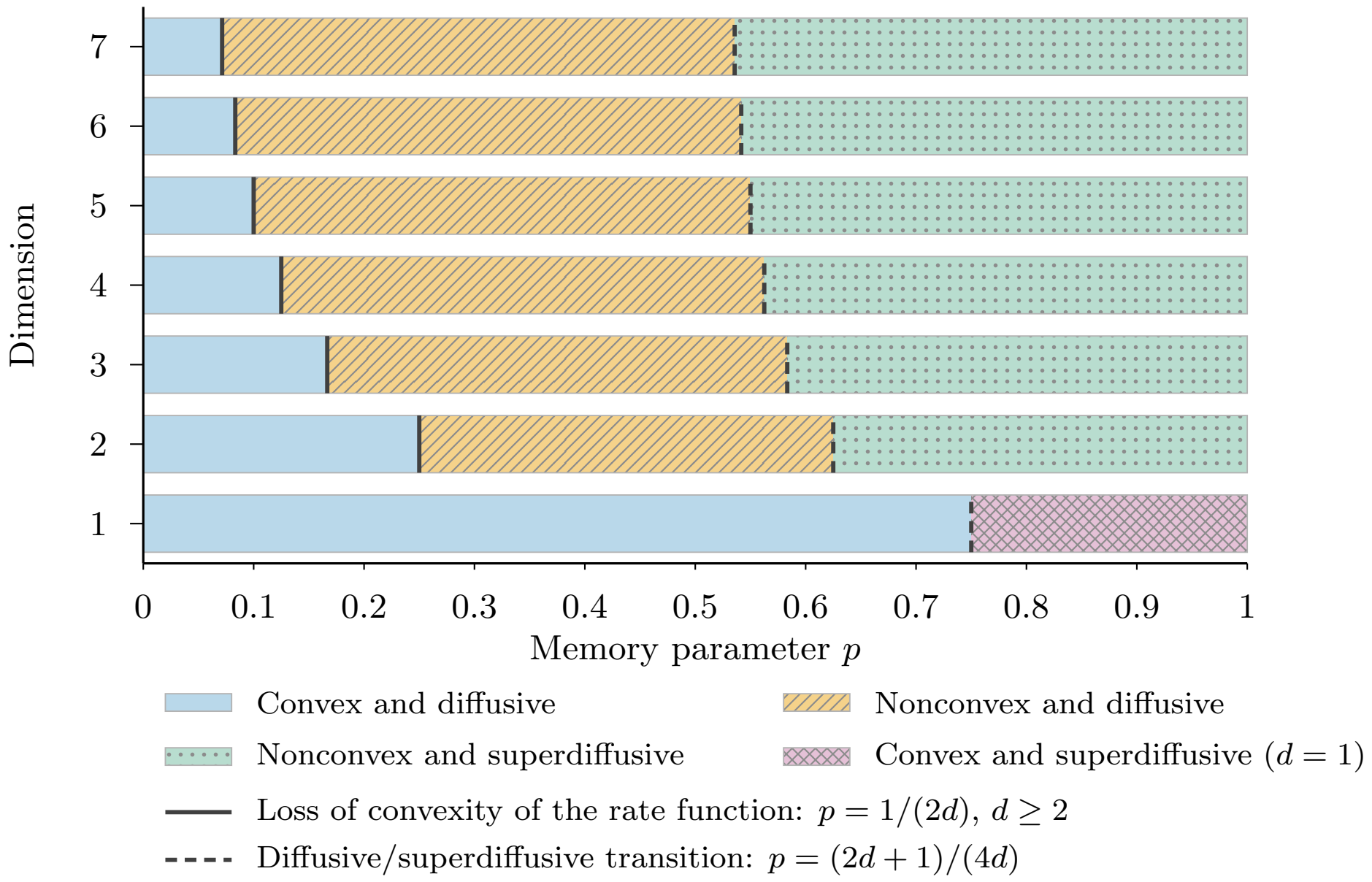
\begin{figure}[htbp]
\centering
\begin{tikzpicture}[
    x=10cm,
    y=0.72cm,
    font=\small
]
\foreach \dimn in {1,...,7} {
    \pgfmathsetmacro{\ymin}{\dimn-0.36}
    \pgfmathsetmacro{\ymax}{\dimn+0.36}
    \pgfmathsetmacro{\pdiff}{(2*\dimn+1)/(4*\dimn)}
    \ifnum\dimn=1
        \path[fill=convexdiffusive]
            (0,\ymin) rectangle (\pdiff,\ymax);
        \path[fill=convexsuperdiffusive]
            (\pdiff,\ymin) rectangle (1,\ymax);
        \path[
            pattern=crosshatch,
            pattern color=black!45
        ]
            (\pdiff,\ymin) rectangle (1,\ymax);
    \else
        \pgfmathsetmacro{\pconv}{1/(2*\dimn)}
        \path[fill=convexdiffusive]
            (0,\ymin) rectangle (\pconv,\ymax);
        \path[fill=nonconvexdiffusive]
            (\pconv,\ymin) rectangle (\pdiff,\ymax);
        \path[
            pattern=north east lines,
            pattern color=black!45
        ]
            (\pconv,\ymin) rectangle (\pdiff,\ymax);
        \path[fill=nonconvexsuperdiffusive]
            (\pdiff,\ymin) rectangle (1,\ymax);
        \path[
            pattern=dots,
            pattern color=black!45
        ]
            (\pdiff,\ymin) rectangle (1,\ymax);
    \fi
    \draw[
        gray!60,
        line width=0.35pt
    ]
        (0,\ymin) rectangle (1,\ymax);
}
\foreach \dimn in {2,...,7} {
    \pgfmathsetmacro{\ymin}{\dimn-0.36}
    \pgfmathsetmacro{\ymax}{\dimn+0.36}
    \pgfmathsetmacro{\pconv}{1/(2*\dimn)}
    \draw[
        black!75,
        line width=1.1pt
    ]
        (\pconv,\ymin) -- (\pconv,\ymax);
}
\foreach \dimn in {1,...,7} {
    \pgfmathsetmacro{\ymin}{\dimn-0.36}
    \pgfmathsetmacro{\ymax}{\dimn+0.36}
    \pgfmathsetmacro{\pdiff}{(2*\dimn+1)/(4*\dimn)}
    \draw[
        black!75,
        line width=1.1pt,
        densely dashed
    ]
        (\pdiff,\ymin) -- (\pdiff,\ymax);
}
\draw[semithick] (0,0.5) -- (1,0.5);
\draw[semithick] (0,0.5) -- (0,7.5);
\foreach \i/\lab in {
    0/0,
    1/0.1,
    2/0.2,
    3/0.3,
    4/0.4,
    5/0.5,
    6/0.6,
    7/0.7,
    8/0.8,
    9/0.9,
    10/1
} {
    \pgfmathsetmacro{\xpos}{\i/10}
    \draw
        (\xpos,0.5) -- ++(0,-0.08)
        node[below=3pt] {$\lab$};
}
\foreach \dimn in {1,...,7} {
    \draw (-0.012,\dimn) -- (0,\dimn);
    \node[anchor=east]
        at (-0.02,\dimn) {$\dimn$};
}
\node[below=15pt] at (0.5,0.5)
    {Memory parameter $p$};
\node[rotate=90] at (-0.11,4)
    {Dimension};
\fill[convexdiffusive]
    (0.02,-1.36) rectangle (0.08,-1.14);
\draw[gray!60,line width=0.35pt]
    (0.02,-1.36) rectangle (0.08,-1.14);
\node[anchor=west,font=\footnotesize]
    at (0.09,-1.25)
    {Convex and diffusive};
\fill[nonconvexdiffusive]
    (0.58,-1.36) rectangle (0.64,-1.14);
\path[
    pattern=north east lines,
    pattern color=black!45
]
    (0.58,-1.36) rectangle (0.64,-1.14);
\draw[gray!60,line width=0.35pt]
    (0.58,-1.36) rectangle (0.64,-1.14);
\node[anchor=west,font=\footnotesize]
    at (0.65,-1.25)
    {Nonconvex and diffusive};
\fill[nonconvexsuperdiffusive]
    (0.02,-2.04) rectangle (0.08,-1.82);
\path[
    pattern=dots,
    pattern color=black!45
]
    (0.02,-2.04) rectangle (0.08,-1.82);
\draw[gray!60,line width=0.35pt]
    (0.02,-2.04) rectangle (0.08,-1.82);
\node[anchor=west,font=\footnotesize]
    at (0.09,-1.93)
    {Nonconvex and superdiffusive};
\fill[convexsuperdiffusive]
    (0.58,-2.04) rectangle (0.64,-1.82);
\path[
    pattern=crosshatch,
    pattern color=black!45
]
    (0.58,-2.04) rectangle (0.64,-1.82);
\draw[gray!60,line width=0.35pt]
    (0.58,-2.04) rectangle (0.64,-1.82);
\node[anchor=west,font=\footnotesize]
    at (0.65,-1.93)
    {Convex and superdiffusive ($d=1$)};
\draw[
    black!75,
    line width=1.1pt
]
    (0.02,-2.63) -- (0.08,-2.63);
\node[anchor=west,font=\footnotesize]
    at (0.09,-2.63)
    {Loss of convexity of the rate function: $p=1/(2d)$, $d\geq 2$};
\draw[
    black!75,
    line width=1.1pt,
    densely dashed
]
    (0.02,-3.23) -- (0.08,-3.23);
\node[anchor=west,font=\footnotesize]
    at (0.09,-3.23)
    {Diffusive/superdiffusive transition:
     $p=(2d+1)/(4d)$};
\end{tikzpicture}

\caption{Phase diagram of the multi-dimensional elephant random walk.}
\label{fig_phase_diagram}

\end{figure}

\medskip
\noindent\textbf{Organisation of the paper.}
Section~2 introduces the notation and preliminary results used throughout the article. In Section~\ref{section_LDP_existence}, we establish a weak LDP for the empirical measures of the MERW and several related observables. Section~\ref{section_properties_rate_function} is devoted to the properties of the rate function, including its symmetries, its comparison with the Donsker--Varadhan rate function, and its convexity. Finally, Section~\ref{section_proof_of_existence_of_stitching_procedure} provides the technical construction of the stitching procedure used to establish the weak LDP.

\par
\addvspace{22pt}
\noindent{\Large\bfseries Acknowledgements}
\par
\addvspace{15pt}

I thank my supervisors, Tristan Benoist, Noé Cuneo, and Clément Pellegrini, for their guidance, support, and careful reading of the manuscript. I also thank Bernard Bercu, Juan Esaul González Rangel, Lucile Laulin and Kilian Raschel for fruitful discussions, and Bastien Mallein for introducing me to the model. I acknowledge financial support from EUR MINT (ANR-18-EURE-0023).

\par
\addvspace{22pt}
\noindent{\Large\bfseries AI statement}
\par
\addvspace{15pt}

\noindent The main results of this article, with the exception of those presented in Section~\ref{subsection_convex_for_small_p}, were established without the use of generative AI, although some were initially proved using arguments different from those presented here. Generative AI was subsequently used to assist with proofreading, revising the manuscript and producing figures.

\section{Notation and preliminaries}

This section collects the notation and auxiliary results used throughout the paper. A summary of the principal notations is provided in Appendix~\ref{section_appendix_list_of_notations}.

\subsection{Topology on the space of probability
measures}\label{subsection_topology_space_of_proba_measures}

We endow $\ZZ^d$ with the distance induced by the $L^1$ norm. Thus, $\ZZ^d$ is a metric space, and therefore the space of probability measures $\mathcal{P}(\ZZ^d)$ endowed with the topology of weak convergence is itself metrisable \cite[Theorem~3.2.2]{Bogachev}. Consider the Kantorovich--Rubinstein (KR) norm on the space $\mathcal{M}(\ZZ^d)$ of finite signed measures, defined by
\begin{equation}\label{eq_def_BL_norm}
\| \mu \| \eqdef \sup \left\{ \int_{\ZZ^d} f(x) \, \mathrm{d} \mu(x) \; : \; f \in \text{Lip}_1(\ZZ^d,\RR)  \; , \; \sup_{k \in \ZZ^d}|f(k)| \leq 1 \right\},
\end{equation}
where $\text{Lip}_1(\ZZ^d,\RR)$ is the space of functions from $\ZZ^d$ to $\RR$ whose Lipschitz constant is bounded by one. The restriction of the distance induced by the KR norm to the convex subset $\mathcal{P}(\ZZ^d) \subseteq \mathcal{M}(\ZZ^d)$ metrises the topology of weak convergence on $\mathcal{P}(\ZZ^d)$, see \cite[Theorem~3.2.2]{Bogachev}. The corresponding metric is also commonly referred to as the bounded-Lipschitz metric, or Fortet--Mourier metric. Depending on the convention used for the bounded-Lipschitz norm, it either coincides with the KR metric defined above or is equivalent to it. In the rest of the article, when a norm is considered on the space of measures it will always be the KR norm and $B(\mu,\varepsilon)$ stands for the open ball in $\mathcal{P}(\ZZ^d)$ centred at $\mu$ and of radius $\varepsilon$ with respect to the distance induced by the KR norm on $\mathcal{P}(\ZZ^d)$.  We recall the following property of the KR norm, which will be used later. For all $\mu \in \mathcal{P}(\ZZ^d)$,
\begin{equation}\label{eq_lemma_properties_KR_norm_of_proba}
\| \mu \| = 1.
\end{equation}
The verification of this property is straightforward.

\subsection{The space of trajectories}\label{subsection_space_of_trajectories}

While the results are stated in a probabilistic setting, our analysis will primarily be carried out at the level of individual trajectories. For all $x,y \in \ZZ^d$, we write $x \sim y$ if $\|x-y\|_1=1$, in other words $x \sim y$ if and only if they are adjacent in $\ZZ^d$, which is also equivalent to saying that $x-y \in \mathscr{S}_d$. Note that for $(S_n)_{n \in \NN}$ we have $S_{n} \sim S_{n+1}$ for all $n \in \NN$. For $A,B \subseteq \ZZ^d$, we also write $x \sim A$ if there exists $y \in A$ such that $x \sim y$, and $A \sim B$ if there exist $x \in A$ and $y \in B$ such that $x \sim y$. We now introduce several sets of trajectories that will be used in the sequel. For all $n \in \NN$, let
\begin{align}\label{def_Omega_n}
\Omega_n& \eqdef \left\{ (x_1, \dots , x_{n}) \in (\ZZ^d)^n : \, x_i \sim x_{i+1} , \, i =1,\dots,n-1 \right\}, \\
\Omega_n^{(0)} & \eqdef \left\{ w \in \Omega_n  : \, w_1=0 \right\}. \notag
\end{align}
We also consider the empty trajectory $\varnothing$ and set $\Omega_0=\{ \varnothing\}=\Omega_0^{(0)}$. We define the set of finite trajectories
\begin{align*}
\Omega_{\textup{fin}} \eqdef \bigcup_{n \in \NN_0} \Omega_n \qquad \text{and} \qquad \Omega_{\textup{fin}}^{(0)} \eqdef \bigcup_{n \in \NN_0} \Omega_n^{(0)}.
\end{align*}

\noindent In what follows, we will often refer to elements $w \in \Omega_{\text{fin}}$ as \textit{words}. A word $w \in \Omega_{\textup{fin}}$ is said to be of length $n$ if $w \in \Omega_n$, and we denote this by $|w|=n$. In particular, $|\varnothing|=0$. The $k$-th element of $w$, denoted $w_k$, will be called a \textit{letter} of the word $w$. We write $w=w_1w_2 \cdots w_n$ for the decomposition of $w$ into its letters. The last letter of $w$ will be denoted by $w_{-1}$. For $w \in \Omega_n$ and integers $p,q \in \NN$ satisfying $1 \leq p \leq q \leq n$, we define $w_{[p:q]}=w_pw_{p+1} \cdots w_q$, and refer to $w_{[p:q]}$ as a \textit{subword} of $w$. Given two words $w^1,w^2 \in \Omega_{\textup{fin}}$ such that $w^1_{-1} \sim w^2_1$, we can define the \textit{concatenation} of $w^1$ and $w^2$ as 
\begin{align*}
w^1 w^2 \eqdef w^1_1 w^1_2 \cdots w^1_{-1} w^2_1 w^2_2 \cdots w^2_{-1}.
\end{align*}
Concatenation with the empty word $\varnothing$ is always allowed and bypasses the adjacency condition. Note that $w^1 w^2 \in \Omega_{\textup{fin}}$ and both $w^1$ and $w^2$ appear as subwords of their concatenation $w^1 w^2$. Moreover, the length of the concatenation of two words satisfies $|w^1 w^2|=|w^1|+|w^2|$. In addition, for all $w \in \Omega_n$, we have $\varnothing w = w = w \varnothing$. Given $w \in \Omega_{\textup{fin}}$, we define its support by setting
\begin{align*}
\textup{Supp}(w) \eqdef \big\{ x \in \ZZ^d \, : \, \exists k \in \dbl 1,|w| \dbr \, , \, w_k=x \big\}.
\end{align*}
We define the empirical measure of a word $w \in \Omega_{n}$, with $w \neq \varnothing$, by setting
\begin{align*}
\ell(w)=\ell_{n}(w) \eqdef \frac{1}{n}\sum_{j=1}^{n} \delta_{w_j}.
\end{align*}
Finally, we define for all $\mu \in \mathcal{P}(\ZZ^d)$ and $\varepsilon>0$ the set
\begin{equation}\label{eq_def_Omega_mu_varepsilon}
\Omega_n(\mu,\varepsilon) \eqdef \big\{ w \in \Omega_n^{(0)} \, : \, \ell_n(w) \in B(\mu,\varepsilon) \big\},
\end{equation}
which consists of all the trajectories of length $n$ that start at $0$ and whose empirical measure is close to $\mu$. 

For any $U_n \subseteq \Omega_n$, we adopt the shorthand notation $\PP_{d,p}(U_n)$ to denote $\PP_{d,p}(S_{[1:n]} \in U_n)$. With this convention, we will also refer to $U_n$ itself as an event, meaning the event $\{S_{[1:n]} \in U_n\}$.

\medskip

For any given measure $\mu \in \mathcal{P}(\ZZ^d)$, one can find $M \in \NN$ such that most of the mass of $\mu$ is concentrated on the compact set $\dbl-M,M \dbr^d$. Consequently, any trajectory
$w \in \Omega_n(\mu,\varepsilon)$
spends most of its time inside this compact set. Since the MERW is a nearest-neighbour random walk, it follows that the trajectory itself cannot stray too far from this compact set, and therefore from the origin. The following lemmas make this observation precise.

For $A \subseteq \ZZ^d$, let $N_A(w)$ be the occupation time of the trajectory $w \in \Omega_{\textup{fin}}$ within $A$:
\begin{equation}\label{eq_def_occupation_time_N_A}
N_A(w)=\# \big\{ j \in \dbl 1,|w| \dbr \, : \, w_j \in A \big\}.
\end{equation}
We also denote $N_x(w)$ instead of $N_{\{x\}}(w)$ when considering singletons.

\begin{lemma}[Bounds on occupation times]\label{lemma_occupation_meas_with_empirical_meas}
Consider $\mu \in \mathcal{P}(\ZZ^d)$, $\varepsilon>0$ and $n \in \NN$. For all $w \in \Omega_{n}(\mu,\varepsilon)$ and $A \subseteq \ZZ^d$,
\begin{equation*}
\left| N_A(w) - \mu(A)n \right| < \varepsilon n.
\end{equation*}
\end{lemma}

\begin{proof}
Fix $w \in \Omega_n(\mu,\varepsilon)$ and $A \subseteq \ZZ^d$. Since distinct points of $\ZZ^d$ are at distance at least one, the map $f=\mathbbm{1}_A$ is Lipschitz with constant at most $1$. In addition, it satisfies $0\leq f\leq 1$. Since $w \in \Omega_{n}(\mu,\varepsilon)$, the definition of the KR norm applied to $f$ and $-f$ yields
\begin{align*}
\varepsilon & > \|\ell_{n}(w)-\mu\| \\
& \geq \left \vert \int_{\ZZ^d} \mathbbm{1}_A \mathrm{d} \ell_{n}(w) - \int_{\ZZ^d} \mathbbm{1}_A \mathrm{d} \mu \right \vert \\
& = \left \vert \frac{N_A(w)}{n}  - \mu(A) \right \vert.
\end{align*}
Multiplying the preceding inequality by $n$ yields the desired result.
\end{proof}

\begin{lemma}[Time spent outside a finite box]\label{lemma_bound_points_outside_compact}
For all $\mu \in \mathcal{P}(\ZZ^d)$ and $\varepsilon>0$, there exists $R_{\mu}(\varepsilon) \in \NN$ such that for all $n \in \NN$ and $w \in \Omega_n(\mu,\varepsilon)$,
\begin{equation*}\label{eq_lemma_bound_points_outside_compact}
\# \big\{ k \in \dbl 1,n \dbr \, : \, \|w_k\|_{\infty} \geq R_{\mu}(\varepsilon) \big\} < 2\varepsilon n.
\end{equation*}
\end{lemma}

\begin{proof}
Fix $\mu \in \mathcal{P}(\ZZ^d)$ and $\varepsilon>0$. For $R \in \RR_+^*$, let $B_{\infty}(0,R)=\{ x \in \ZZ^d \, : \, \|x\|_{\infty}<R \}$ and set  
\begin{align*}
R_{\mu}(\varepsilon) \eqdef \inf \big\{ k \in \NN \, : \, \mu(B_{\infty}(0,k)) > 1-\varepsilon \big\}.
\end{align*}
By definition of $R_{\mu}(\varepsilon)$, we have $\mu\big(B_{\infty}\big(0,R_{\mu}(\varepsilon)\big)^c\big) < \varepsilon$. The desired result then follows by applying Lemma~\ref{lemma_occupation_meas_with_empirical_meas} with $A=B_{\infty}\big(0,R_{\mu}(\varepsilon)\big)^c$.
\end{proof}

\begin{corollary}[Bound on the maximal displacement]\label{corollary_bound_on_last_point_of_subword}
Fix $\mu \in \mathcal{P}(\ZZ^d)$, $\varepsilon>0$ and let $R_{\mu}(\varepsilon)$ be as in Lemma~\ref{lemma_bound_points_outside_compact}. If $n > \frac{R_{\mu}(\varepsilon)}{\varepsilon}$ and $w \in \Omega_n(\mu,\varepsilon)$, then
\begin{equation*}
\|w_k\|_{\infty} < 3\varepsilon n, \quad k \in \dbl 1,n \dbr.
\end{equation*}
\end{corollary}

\begin{proof}
Suppose for the sake of contradiction that there exists $k \in \dbl 1,n \dbr$ such that $\|w_k\|_{\infty} \geq 3 \varepsilon n$. Since $n > \frac{R_{\mu}(\varepsilon)}{\varepsilon}$, we deduce that $\|w_k\|_{\infty} > R_{\mu}(\varepsilon)+2\varepsilon n$. In addition, since trajectories $w \in \Omega_n(\mu,\varepsilon)$ satisfy $w_i \sim w_{i+1}$ for all $i \in \dbl 1,n-1 \dbr$, we must have $k \geq R_{\mu}(\varepsilon)+2\varepsilon n$ and
\begin{align*}
\|w_{k-j}\|_{\infty} > R_{\mu}(\varepsilon), \quad \forall j\in \dbl 0,\lfloor 2 \varepsilon n \rfloor \dbr.
\end{align*}
Therefore, we have
\begin{align*}
\# \big\{ j \in \dbl 1,n \dbr \, : \, \|w_j\|_{\infty} \geq R_{\mu}(\varepsilon) \big\} > 2 \varepsilon n,
\end{align*}
which contradicts Lemma~\ref{lemma_bound_points_outside_compact}. Thus, such $k$ cannot exist and we deduce that $\|w_k\|_{\infty} < 3\varepsilon n$ for all $k \in \dbl 1,n \dbr$.
\end{proof}

\subsection{Profile of a trajectory}

As mentioned previously, the MERW is in general non-Markovian. The position alone does not, in general, determine the next-step distribution. Indeed, for every $n \geq 2$ and $s \in \mathscr{S}_d$, the conditional probability
\begin{align*}
\PP_{d,p}\big(S_{n+1}=S_n+s \mid S_1 , \dots , S_n \big)
\end{align*}
can be expressed in terms of the proportion of increments taken in the direction $s$ up to time $n$. This dependence is made explicit in Lemma~\ref{lemma_conditional_transition_probability}, and the resulting formula for the probability of a trajectory is given in Corollary~\ref{corollary_product_of_transitions}. This motivates the introduction of the \textit{profile} of a trajectory, which records the cumulative evolution of its increments in each direction.

\medskip

First, consider the simplex over $\mathscr{S}_d$ defined by
\begin{equation*}\label{eq_def_dimplex}
\Delta^{(d)}=\left\{ (p_s)_{s \in \mathscr{S}_d} \in [0,1]^{\mathscr{S}_d} \, : \, \sum_{s \in \mathscr{S}_d} p_s=1 \right\}.
\end{equation*}
We endow $\RR^{\mathscr{S}_d}$ with the $\| \cdot \|_{\infty}$ norm and for an element $x \in \RR^{\mathscr{S}_d}$, we denote by $[x]_s$ the coordinate of $x$ indexed by $s$. We define
\begin{equation}\label{eq_def_H_d_space_of_profiles}
\mathscr{H}^{(d)} \eqdef
\left\{ H \in \mathcal{C}^0\big([0,1],(\RR_+)^{\mathscr{S}_d}\big) \, : \,
\raisebox{-0.5ex}{$\displaystyle
\begin{array}{c}
[H(u)]_s \leq [H(t)]_s, \quad s \in \mathscr{S}_d, \; 0 \leq u \leq t \leq 1, \\[0.3em]
H(0)=0, \qquad \displaystyle\sum_{s \in \mathscr{S}_d} [H(t)]_s = t, \quad \forall t \in [0,1]
\end{array}
$}
\right\}.
\end{equation}
The space $\mathscr{H}^{(d)}$ will serve as the state space for the trajectory profiles. We endow $\mathscr{H}^{(d)}$ with the metric induced by uniform convergence,
\begin{equation}\label{eq_def_uniform_distance_for_profiles}
d_{\mathscr{H}}(H, \widehat{H}) \eqdef \sup_{t \in [0,1]} \| H(t)-\widehat{H}(t) \|_{\infty}, \qquad H,\widehat{H} \in \mathscr{H}^{(d)}.
\end{equation}
We further denote by $B_{\mathscr{H}}(H,\varepsilon)$ the open ball centred at $H$ and of radius $\varepsilon$ for the distance $d_{\mathscr{H}}$.

\medskip

For all $w \in \Omega_{\textup{fin}}$ and $j \in \dbl 1, |w|-1 \dbr$, set
\begin{equation}\label{eq_Delta_word_j}
\Delta_j(w) \eqdef w_{j+1}-w_j \in \mathscr{S}_d,
\end{equation}
which gives the direction of the $j$-th step of the word $w$. For $n \in \NN$, $w \in \Omega_n$, $1 \leq k \leq n$ and $s \in \mathscr{S}_d$, set
\begin{equation}\label{eq_def_step_counter}
C^{(s)}_{k}(w) \eqdef \# \big\{j \in \dbl 1,k-1 \dbr \, : \, \Delta_j(w)=s \big\}.
\end{equation}
In particular, $C_1^{(s)}(w)=0$. The quantity $C^{(s)}_{k}(w)$ counts the number of times the step $s$ has been taken by the trajectory $w$ up to the time $k$. 

\begin{definition}[Profile of a trajectory]\label{definition_profile_H_n}
Let $n \geq 2$ and $w \in \Omega_n^{(0)}$. We define its \textit{profile} $H_n(w,\cdot) \in \mathscr{H}^{(d)}$ by setting, for every $j \in \dbl 1,n \dbr$ and $s \in \mathscr{S}_d$, 
\begin{equation}\label{eq_def_profile_map}
H_n\left(w,\frac{j-1}{n-1} \right) = \left( \frac{C^{(s)}_{j}(w)}{n-1} \right)_{s \in \mathscr{S}_d}, 
\end{equation}
and then extending $H_n(w,\cdot)$ affinely on each interval $\left[\frac{j-1}{n-1},\frac{j}{n-1}\right]$, $j \in \dbl 1,n-1 \dbr$. For $n=1$, we set $H_1(w,t)=t(1/(2d))_{s \in \mathscr{S}_d}$. This convention will be of no importance for all asymptotic statements below.
\end{definition} 

\noindent The profile $H_n(w,\cdot)$ records the cumulative number of increments in every direction, normalised by the total number $n-1$ of increments along the trajectory $w \in \Omega_n^{(0)}$. By \eqref{eq_def_H_d_space_of_profiles}, for every $H \in \mathscr{H}^{(d)}$ and $t \in (0,1]$, the vector $H(t)/t$ belongs to $\Delta^{(d)}$. In particular, when $H=H_n(w,\cdot)$ and $t=(j-1)/(n-1)$, with $2 \leq j \leq n$, this vector is the empirical distribution of the first $j-1$ increments of $w$. Thus, the profile records the information needed to determine the transition probabilities, as shown in Section~\ref{subsection_proba_of_trajectory}. We further write $H_n$ for the map
\begin{equation}
\begin{array}{cccc}
H_n : & \Omega_n^{(0)} & \longrightarrow & \mathscr{H}^{(d)} \\
& w & \longmapsto & H_n(w,\cdot)
\end{array}.
\end{equation}
In what follows, we slightly abuse notation and also write $H_n$ for the random profile $H_n(S_{[1:n]},\cdot)$.

\begin{example}[Two different profiles for $d=3$]\label{example_balanced_unbalanced_profile}
    We identify $\RR^{\mathscr{S}_3}$ with $\RR^6$ through the ordering $\mathscr{S}_3=(e_1,-e_1,e_2,-e_2,e_3,-e_3)$. Thus an element $x \in \RR^{\mathscr{S}_3}$ is represented by
\begin{align*}
    \big(x_{e_1}, x_{-e_1} , x_{e_2}, x_{-e_2}, x_{e_3}, x_{-e_3}\big) \in \RR^6.
\end{align*}
Consider the measure $\widehat{\mu}=1/2 \delta_0 + 1/4 \delta_{e_1} + 1/4 \delta_{e_2} \in \mathcal{P}(\ZZ^3)$. For $n \in \NN$, define the trajectories 
\begin{align*}
    w^{(n),\text{alt}} & \eqdef (0,e_1,0,e_2,0,e_1,0,e_2 \cdots , 0,e_1,0,e_2) \in \Omega_{4n} , \\ 
    w^{(n),\text{seq}} & \eqdef (\underbrace{0,e_1,0,e_1, \cdots , 0,e_1}_{\text{of length } 2n},\underbrace{0,e_2,0,e_2, \cdots,0,e_2}_{\text{of length } 2n}) \in \Omega_{4n}.
\end{align*}
Both trajectories have the same empirical measure, namely
\begin{align*}
    \ell_{4n}\big( w^{(n),\text{alt}} \big)=\widehat{\mu}=\ell_{4n}\big( w^{(n),\text{seq}} \big)
\end{align*}
Nevertheless, both trajectories have very different profiles. The first trajectory alternates visits to $e_1$ and $e_2$, whereas the second makes them sequentially. More precisely, define $H^{\text{alt}}$ and $H^{\text{seq}}$ by
\begin{align*}
    H^{\text{alt}}(t) & = \left( \frac{t}{4} , \frac{t}{4} , \frac{t}{4} , \frac{t}{4} , 0 , 0 \right) \qquad \text{and} \qquad H^{\text{seq}}(t) = 
    \begin{cases}
        \left( \frac{t}{2}, \frac{t}{2},0,0,0,0 \right), & \text{if } t \in \left[ 0,\frac{1}{2} \right] \\
        \left( \frac{1}{4} , \frac{1}{4} , \frac{t}{2}-\frac{1}{4} , \frac{t}{2}-\frac{1}{4} , 0 ,0 \right), & \text{if } t \in \left[ \frac{1}{2},1 \right]
    \end{cases}.
\end{align*}
A direct count on the increments shows that  
\begin{align*}
    H_{4n}\big( w^{(n),\text{alt}} , \cdot \big) \xrightarrow[n \to \infty]{d_{\mathscr{H}}} H^{\text{alt}} \qquad \text{and} \qquad H_{4n}\big( w^{(n),\text{seq}}, \cdot\big) \xrightarrow[n \to \infty]{d_{\mathscr{H}}} H^{\text{seq}}.
\end{align*}
\end{example}

\begin{lemma}[Admissibility of trajectory profiles]\label{lemma_profiles_in_profile_space}
For all $n \in \NN$ and $w \in \Omega_n^{(0)}$,
\begin{equation}\label{eq_lemma_profiles_in_profile_space}
H_n(w,\cdot) \in \mathscr{H}^{(d)}.
\end{equation}
\end{lemma}

\begin{proof}
Fix $n \geq 2$. For every $w \in \Omega_n^{(0)}$ and $s \in \mathscr{S}_d$, the sequence $k \mapsto C_k^{(s)}(w)$ is non-decreasing. Hence, by the definition of $H_n$, the map $t \mapsto [H_n(w,t)]_s$ is non-decreasing on the mesh points $(j-1)/(n-1)$. Since $H_n(w,\cdot)$ is extended affinely between consecutive mesh points, this monotonicity is preserved, and therefore
\begin{align*}
[H_n(w,u)]_s \leq [H_n(w,t)]_s, \qquad \forall \, 0 \leq u \leq t \leq 1.
\end{align*}
Moreover, the condition $H_n(w,0)=0$ follows from \eqref{eq_def_profile_map} with $j=1$. Finally, for every $j \in \dbl 1,n \dbr$, we have
\begin{align*}
\sum_{s \in \mathscr{S}_d} \left[ H_n \left( w, \frac{j-1}{n-1} \right) \right]_s & = \frac{1}{n-1} \sum_{s \in \mathscr{S}_d} C_j^{(s)}(w) \\
& = \frac{j-1}{n-1}.
\end{align*}
Since $H_n(w,\cdot)$ is affine between consecutive mesh points, it follows that, 
\begin{align*}
\sum_{s \in \mathscr{S}_d} [H_n(w,t)]_s=t, \qquad \forall \, t \in [0,1].
\end{align*}
Thus, all three conditions defining $\mathscr{H}^{(d)}$ are satisfied, and consequently $H_n(w,\cdot) \in \mathscr{H}^{(d)}$. The case $n=1$ is immediate from the definition of $H_1$.
\end{proof}

\noindent The following lemma gives useful properties of the space $\mathscr{H}^{(d)}$ and of profiles $H \in \mathscr{H}^{(d)}$. We further write $H_n'(w,t)$ for $\partial_t H_n(w,t)$.

\begin{lemma}[Properties of the profile space]\label{lemma_properties_profile_space}
The following properties hold :
\begin{enumerate}
\item For every $n \geq 2$ and $w \in \Omega_n^{(0)}$, the profile $H_n(w,\cdot)$ satisfies
\begin{equation}\label{eq_lemma_properties_profile_space_derivative}
H_n'(w,t)=\delta_{\Delta_j(w)}, \qquad \forall \, t \in \left(\frac{j-1}{n-1},\frac{j}{n-1} \right), \quad j \in \dbl 1,n-1 \dbr,
\end{equation} 
where $\delta_{s} \in \Delta^{(d)}$ denotes the probability vector concentrated on $s \in \mathscr{S}_d$.
\item Every $H \in \mathscr{H}^{(d)}$ is 1-Lipschitz, absolutely continuous, and
\begin{equation}\label{eq_lemma_properties_profile_space_derivative_everywhere}
H'(t) \in \Delta^{(d)}, \qquad \text{for a.e. } t \in [0,1]. 
\end{equation}
\item The space $\mathscr{H}^{(d)}$ is a compact metric space.
\end{enumerate}
\end{lemma}

\begin{proof}
\begin{enumerate}
\item Equation~\eqref{eq_lemma_properties_profile_space_derivative} follows from the definition of $H_n(w,\cdot)$ at the points $(j-1)/(n-1)$ and $j/(n-1)$, and the fact that it is extended affinely between these two mesh points.

\item Let $H \in \mathscr{H}^{(d)}$ and $0 \leq u \leq t \leq 1$. Since every coordinate is non-decreasing, 
\begin{align*}
0 \leq [H(t)]_s-[H(u)]_s \leq \sum_{r \in \mathscr{S}_d} \big([H(t)]_r-[H(u)]_r\big) = t-u,
\end{align*}
where the last equality follows from $\sum_{r \in \mathscr{S}_d} [H(v)]_r=v$ for every $v \in [0,1]$. Taking the maximum over all $s \in \mathscr{S}_d$, we obtain
\begin{align*}
\| H(t)-H(u) \|_{\infty} \leq |t-u|, \qquad \forall \, t,u \in [0,1].
\end{align*}
In other words, every element of $\mathscr{H}^{(d)}$ is 1-Lipschitz and hence absolutely continuous. In particular, every coordinate function $[H(\cdot)]_s$ is differentiable almost everywhere, see \cite[p.~124, Theorem~10]{Real_analysis}. Since $\mathscr{S}_d$ is finite, there exists a set $A \subseteq [0,1]$ of full Lebesgue measure on which all the coordinate functions are differentiable. Fix $t \in A$. Since $[H(\cdot)]_s$ is non-decreasing, we deduce that $[H'(t)]_s \geq 0$ for all $s \in \mathscr{S}_d$. Moreover, differentiating the identity $\sum_{s \in \mathscr{S}_d} [H(u)]_s=u$ at $u=t$ yields $\sum_{s \in \mathscr{S}_d} [H'(t)]_s=1$. Therefore $H'(t) \in \Delta^{(d)}$ for all $t \in A$, which proves \eqref{eq_lemma_properties_profile_space_derivative_everywhere}.

\item We first show that $\mathscr{H}^{(d)}$ is closed for the uniform topology. Let $(H_k)_{k \in \NN} \subseteq \mathscr{H}^{(d)}$ converge uniformly to some $H \in \mathcal{C}^0\big([0,1],\RR_+^{\mathscr{S}_d}\big)$. Then, by pointwise convergence, we obtain that $H(0)=0$, $\sum_{s \in \mathscr{S}_d} [H(t)]_s=t$ and that each coordinate of $H$ is non-decreasing. Hence $H \in \mathscr{H}^{(d)}$, and $\mathscr{H}^{(d)}$ is closed. 

By the second item of Lemma~\ref{lemma_properties_profile_space}, every $H \in \mathscr{H}^{(d)}$ is 1-Lipschitz. Since $H(0)=0$, we also have $\|H(t)\|_{\infty} \leq 1$ for every $t \in [0,1]$. Thus, $\mathscr{H}^{(d)}$ is uniformly bounded and equicontinuous. The Arzelà--Ascoli theorem therefore implies that $\mathscr{H}^{(d)}$ is relatively compact for the uniform topology. Since it is also closed, it is compact. Equipped with the metric induced by the uniform norm, it is consequently a compact metric space.  \qedhere
\end{enumerate}
\end{proof}

\begin{remark}\label{remark_density_of_attainable_profiles}
It is worth noting that the set of attainable profiles is dense in $\mathscr{H}^{(d)}$ for the uniform topology. More precisely,
\begin{equation*}
\overline{\bigcup_{n \in \NN} H_n\left( \Omega_n^{(0)} \right)} = \mathscr{H}^{(d)},
\end{equation*}
where the closure is taken for the uniform topology. Indeed, partition $[0,1]$ into sufficiently small intervals. On each interval, approximate the average slope of $H$, which belongs to $\Delta^{(d)}$, by a convex combination of the vectors $\{ \delta_s : s \in \mathscr{S}_d\}$. Such a convex combination can be realised by a finite sequence of steps having the prescribed proportions of directions. Concatenating these blocks yields a finite trajectory whose profile is uniformly close to $H$.
\end{remark}

\subsection{Path probabilities}\label{subsection_proba_of_trajectory}

We now establish the probability that a MERW realises a fixed trajectory $w \in \Omega_n^{(0)}$. The following lemma starts by giving the conditional transition probabilities. For every $n \in \NN$, let 
\begin{align*}
\mathcal{F}_n \eqdef \sigma\big(S_1, \dots , S_n \big).
\end{align*}
Thus, $(\mathcal{F}_n)_{n \in \NN}$ is the natural filtration of the MERW $(S_n)_{n \in \NN}$. By the recursive definition of the increments in \eqref{eq_def_ERW_increments}, we have for every $n \geq 2$,
\begin{align*}
\mathcal{F}_n \subseteq \sigma \big( X_1,A_2,\dots,A_{n-1},U_1,\dots,U_{n-2} \big).
\end{align*}
For $n \geq 2$, $w \in \Omega_n$ and $2 \leq k \leq n$, let
\begin{equation}\label{eq_def_empirical_distribution_rho}
\rho_k(w) \eqdef \left(\frac{C_k^{(s)}(w)}{k-1}\right)_{s \in \mathscr{S}_d} \in\Delta^{(d)}
\end{equation}
denote the empirical distribution of the first $k-1$ increments of $w$. For $w\in\Omega_n^{(0)}$, the definition of the profile gives
\begin{equation*}
H_n \left(w,\frac{k-1}{n-1}\right) = \frac{k-1}{n-1}\rho_k(w), \qquad 2 \leq k \leq n.
\end{equation*}
We further define for $r \in \Delta^{(d)}$ and $s \in \mathscr{S}_d$,
\begin{equation}\label{eq_def_Q_p_r_s}
Q_p(s \,|\, r) \eqdef p r_s + \frac{1-p}{2d-1}(1-r_s).
\end{equation}
The quantity $Q_p(s \mid r)$ represents the one-step transition probability of the MERW, namely the probability of selecting the step $s$ when the empirical distribution of past steps is $r \in \Delta^{(d)}$. The following lemma makes this interpretation precise.

\begin{lemma}[Conditional step probability]\label{lemma_conditional_transition_probability}
Fix $n \geq 2$ and $w \in \Omega_n^{(0)}$. Then, for all $k \in \dbl 2,n-1 \dbr$, we have $\PP_{d,p}$-almost surely
\begin{align*}
\PP_{d,p} \big( S_{[1:k+1]}=w_{[1:k+1]} \vert \mathcal{F}_k \big) = \mathbbm{1}_{\{S_{[1:k]}=w_{[1:k]}\}}Q_p \big( \Delta_k(w) \mid \rho_k(w) \big). 
\end{align*}
\end{lemma}

\begin{proof}
Since $(S_j)_{j=1}^k$ is $\mathcal{F}_k$-measurable,
\begin{equation}\label{eq_conditional_proba_step_indicator}
\PP_{d,p} \big( S_{[1:k+1]}=w_{[1:k+1]} \vert \mathcal{F}_k \big) = \mathbbm{1}_{\{S_{[1:k]}=w_{[1:k]}\}}\PP_{d,p}\left(X_{k}=\Delta_k(w) \vert \mathcal{F}_k\right).
\end{equation}
It remains to compute the second factor. Here and throughout, $\EE_{d,p}$ denotes expectation with respect to $\PP_{d,p}$. By the definition of the MERW increments in \eqref{eq_def_ERW_increments}, we have $X_k=A_k X_{U_{k-1}}$. Therefore, conditioning on the value of $U_{k-1}$ and using the tower property, we obtain
\begin{align*}
\PP_{d,p}\left(X_{k}=\Delta_k(w) \vert \mathcal{F}_k\right) & = \sum_{j=1}^{k-1} \EE_{d,p}\left[ \mathbbm{1}_{\{A_k X_j = \Delta_k(w)\}} \mathbbm{1}_{\{U_{k-1}=j\}} \big| \mathcal{F}_k \right] \\
& = \sum_{j=1}^{k-1} \EE_{d,p} \left[ \mathbbm{1}_{\{U_{k-1}=j\}} \EE_{d,p}\left[ \mathbbm{1}_{\{A_k X_j = \Delta_k(w)\}} \big| \mathcal{F}_k \vee \sigma(U_{k-1}) \right] \big\vert \mathcal{F}_k \right].
\end{align*}
By partitioning according to the value of $X_j \in \mathscr{S}_d$ and using the independence of $A_k$ from $\mathcal{F}_k \vee \sigma(U_{k-1})$, we deduce that
\begin{align*}
\EE_{d,p}\left[ \mathbbm{1}_{\{A_k X_j = \Delta_k(w)\}} \big| \mathcal{F}_k \vee \sigma(U_{k-1}) \right] & = \sum_{s \in \mathscr{S}_d} \PP_{d,p}(A_k s = \Delta_k(w)) \mathbbm{1}_{\{X_j=s\}}.
\end{align*}
By definition of the law of $A_k$ and the definition of $Q_p(\cdot \mid \cdot)$ given in \eqref{eq_def_Q_p_r_s}, we obtain
\begin{align*}
\EE_{d,p}&\left[ \mathbbm{1}_{\{A_k X_j = \Delta_k(w)\}} \big| \mathcal{F}_k \vee \sigma(U_{k-1}) \right] \\
& \qquad \qquad \qquad \qquad = \sum_{s \in \mathscr{S}_d} \left( p \mathbbm{1}_{\{s=\Delta_k(w)\}} + \sum_{M \neq \mathbb{I}_d} \frac{1-p}{2d-1}\mathbbm{1}_{\{Ms=\Delta_k(w)\}} \right)\mathbbm{1}_{\{X_j=s\}} \\
&\qquad \qquad \qquad \qquad= p \mathbbm{1}_{\{X_j=\Delta_k(w)\}} + \frac{1-p}{2d-1} \big(1-\mathbbm{1}_{\{X_j = \Delta_k(w)\}}\big).
\end{align*}
Since $X_j=S_{j+1}-S_j$, we deduce that for all $j \in \dbl 1,k-1 \dbr$, the random variable $X_j$ is $\mathcal{F}_k$-measurable. Hence, 
\begin{align*}
\PP_{d,p}\left(X_{k}=\Delta_k(w) \vert \mathcal{F}_k\right) & = \sum_{j=1}^{k-1} \Big( p \mathbbm{1}_{\{X_j=\Delta_k(w)\}} + \frac{1-p}{2d-1} \big(1-\mathbbm{1}_{\{X_j = \Delta_k(w)\}}\big) \Big)\EE_{d,p} \left[ \mathbbm{1}_{\{U_{k-1}=j\}} \big \vert \mathcal{F}_k \right].
\end{align*}
Since $U_{k-1}$ is independent of $\mathcal{F}_k$ and follows the uniform distribution on $\dbl 1,k-1 \dbr$, we deduce that
\begin{align*}
& \PP_{d,p}\left(X_{k}=\Delta_k(w) \vert \mathcal{F}_k\right) \\
& \qquad \qquad = p \left( \frac{1}{k-1} \sum_{j=1}^{k-1} \mathbbm{1}_{\{S_{j+1}-S_j=\Delta_k(w)\}} \right)+\frac{1-p}{2d-1}\left( 1-\left( \frac{1}{k-1} \sum_{j=1}^{k-1} \mathbbm{1}_{\{S_{j+1}-S_j=\Delta_k(w)\}} \right) \right).
\end{align*}
On the event $\{S_{[1:k]}=w_{[1:k]}\}$, we have, by the definition of the step counter $C^{(s)}_k$, given in \eqref{eq_def_step_counter}, that
\begin{align*}
\frac{1}{k-1}\sum_{j=1}^{k-1} \mathbbm{1}_{\{S_{j+1}-S_j=\Delta_k(w)\}} & = \frac{1}{k-1} \sum_{j=1}^{k-1} \mathbbm{1}_{\{w_{j+1}-w_j=\Delta_k(w)\}} \\
& = \frac{C^{(\Delta_k(w))}_k(w)}{k-1}.
\end{align*}
Using the definition of $\rho_k(w)$ given in \eqref{eq_def_empirical_distribution_rho}, we have
\begin{align*}
\frac{1}{k-1}\sum_{j=1}^{k-1} \mathbbm{1}_{\{S_{j+1}-S_j=\Delta_k(w)\}}=[\rho_k(w)]_{\Delta_k(w)}.
\end{align*}
We deduce from \eqref{eq_conditional_proba_step_indicator} and the definition of $Q_p(\cdot \mid \cdot)$, given in \eqref{eq_def_Q_p_r_s}, that
\begin{align*}
\PP_{d,p} \big( S_{[1:k+1]}=w_{[1:k+1]} \vert \mathcal{F}_k \big) & = \mathbbm{1}_{\{S_{[1:k]}=w_{[1:k]}\}}\PP_{d,p}\left(X_{k}=\Delta_k(w) \vert \mathcal{F}_k\right) \\
& =\mathbbm{1}_{\{S_{[1:k]}=w_{[1:k]}\}} Q_p \big( \Delta_k(w) \mid \rho_k(w) \big),
\end{align*}
and the lemma is proved.
\end{proof}

\begin{corollary}[Path probability]\label{corollary_product_of_transitions}
Let $n \geq 2$. For all $w \in \Omega_n^{(0)}$, 
\begin{equation}\label{eq_corollary_product_of_transitions}
\PP_{d,p}\big( S_{[1:n]}=w \big) = \frac{1}{2d} \prod_{j=2}^{n-1} Q_p \big( \Delta_j(w) \mid \rho_j(w) \big).
\end{equation}
\end{corollary}

\begin{proof}
Using Lemma~\ref{lemma_conditional_transition_probability} with $k=n-1$, conditioning on $\mathcal{F}_{n-1}$, we obtain
\begin{align*}
\PP_{d,p}\big( S_{[1:n]}=w \big) & = \EE_{d,p} \left[ \mathbbm{1}_{\{S_{[1:n-1]}=w_{[1:n-1]}\}} Q_p \big( \Delta_{n-1}(w) \mid \rho_{n-1}(w) \big) \right] \\
& = \PP_{d,p} \left( S_{[1:n-1]}=w_{[1:n-1]} \right) Q_p \big( \Delta_{n-1}(w) \mid \rho_{n-1}(w) \big).
\end{align*}
The result follows by induction combined with the fact that the first step is drawn according to the initial distribution $\nu_{\text{unif}}$, so that
\begin{equation*}
\PP_{d,p}(S_{[1:2]}=w_{[1:2]})=\frac{1}{2d}. \qedhere
\end{equation*}
\end{proof}

\begin{remark}\label{remark_formula_different_inital_distrib}
More generally, for $\nu \in \mathcal{P}(\mathscr{S}_d)$, let $\PP_{d,p,\nu}$ denote the law of the MERW whose initial step has distribution $\nu$, not necessarily uniform on $\mathscr{S}_d$, as introduced in Section~\ref{subsection_invariance_of_rate_function}. The same argument gives
\begin{equation}\label{eq_remark_formula_different_inital_distrib}
\PP_{d,p,\nu}\big( S_{[1:n]}=w \big) = \nu\big( \Delta_1(w) \big) \prod_{j=2}^{n-1} Q_p \big( \Delta_j(w) \mid \rho_j(w) \big).
\end{equation}
\end{remark}

We introduce the constant 
\begin{equation}\label{eq_lower_bound_transition_probability}
p_* \eqdef \min \left\{ \frac{1-p}{2d-1} \,  , \, p \right\}>0,
\end{equation}
where the strict positivity follows from the fact that $p \in (0,1)$. 
\begin{lemma}[Uniform bounds on transition probabilities]\label{lemma_lower_bound_transition}
For every $r \in \Delta^{(d)}$ and $s \in \mathscr{S}_d$,
\begin{equation}\label{eq_lemma_lower_bound_transition}
p_* \leq Q_p(s \mid r) \leq 1.
\end{equation}
\end{lemma}

\begin{proof}
By \eqref{eq_def_Q_p_r_s}, the value of $Q_p(s \mid r)$ is the affine function $p r_s + \frac{1-p}{2d-1}(1-r_s)$ of $r_s \in [0,1]$, so its minimum and maximum are reached at $r_s=0$ or $r_s=1$. 
\end{proof} 

\section{Existence of a weak large deviation principle}\label{section_LDP_existence}

The aim of this section is to establish a weak LDP for the sequence of empirical measures of the elephant random walk $(\ell_n)_{n \in \NN}$, defined in \eqref{eq_def_empirical_measure}. To do so, we first establish this for the augmented sequence $(\ell_n,H_n)_{n\in\NN}$ in the product space $ \mathcal{P}(\ZZ^d)\times\mathscr{H}^{(d)}$, where $H_n$ is the profile associated with the walk, defined in \eqref{eq_def_profile_map}.  Our proof relies on a decoupling argument in the spirit of the Ruelle--Lanford approach to large deviations, following the strategy developed in \cite{Cuneo_selective_decoupling} for selectively decoupled processes and, more recently, in \cite{Daures,Daures_2,JLF_RW_compactified} for empirical measures of Markov chains. Its main ingredient is a stitching procedure which combines shorter trajectories whose empirical measures and profiles are close to prescribed targets $(\mu,H) \in \mathcal{P}(\ZZ^d) \times \mathscr{H}^{(d)}$ into a single longer trajectory. This is achieved by concatenating suitable portions of these trajectories and inserting short transition paths between consecutive portions. The construction simultaneously controls the empirical measure and profile of the resulting trajectory, the number of possible preimages under the stitching map, and the loss of probability caused by the concatenation. It thereby yields the local probability estimates required to compare different time scales. Since the construction and the verification of its properties are rather technical, they are deferred to Section~\ref{section_proof_of_existence_of_stitching_procedure}.

The section is organised as follows. In Section~\ref{subsection_decoupling_method}, we introduce the upper and lower Ruelle--Lanford functions and recall that their equality at every point yields a weak LDP. In Section~\ref{subsection_LDP_symmetric_case}, we consider the symmetric parameter $p=1/(2d)$, for which the MERW has the same law as the simple symmetric random walk. Using the stitching procedure, we prove that the Ruelle--Lanford function exists everywhere, thereby obtaining a weak LDP for the augmented sequence under $\PP_{d,1/(2d)}$, with rate function $I_{\textup{A}}^{(1/(2d))}$, where the subscript $\text{A}$ stands for \textit{augmented}. We also show that the corresponding rate function is convex. In Section~\ref{subsection_comparaison_RL_functions}, we extend the weak LDP to every $p \in (0,1)$ by means of an exponential tilting determined by the profile. More precisely, after introducing the profile-energy functional $\mathcal{E}_p$ on $\mathscr{H}^{(d)}$, we compare the trajectory probabilities under $\PP_{d,p}$ and $\PP_{d,1/(2d)}$ and deduce that for every profile $H \in \mathscr{H}^{(d)}$,
\begin{align*}
I_{\textup{A}}^{(p)}(\mu,H) = I_{\textup{A}}^{(1/(2d))}(\mu,H)-\mathcal{E}_p(H).
\end{align*}
\noindent In Section~\ref{subsection_existence_of_a_weak_LDP_measures}, we establish a weak LDP for the augmented sequence $(\ell_n,H_n)_{n \in \NN}$ and apply the contraction principle to obtain a weak LDP for the empirical measures $(\ell_n)_{n \in \NN}$, thereby proving Theorem~\ref{theorem_existence_LDP}. We further show that this weak LDP cannot be strengthened to a full LDP, which proves Proposition~\ref{prop_no_full_LDP}. Finally, as a by-product of the exponential tilting obtained in Section~\ref{subsection_comparaison_RL_functions}, we derive in Section~\ref{subsection_profile_and_path_observables} full large deviation principles for the profile, the rescaled position, and the maximal displacement of the MERW. 

\medskip

We endow $\mathcal{P}(\ZZ^d)\times\mathscr{H}^{(d)}$ with the product topology and its Borel $\sigma$-algebra. For every $n\in\NN$, the pair $(\ell_n,H_n)$ is measurable. Indeed, it depends only on the finite trajectory $S_{[1:n]}$, so the preimage of any Borel set is a finite union of events $\{S_{[1:n]}=w\}$, with $w\in\Omega_n^{(0)}$. Hence the standing measurability assumption used in Dembo--Zeitouni is satisfied, see \cite[Remark p.~5 and p.~116]{DZ_LDP}. We fix throughout this section, once and for all, the dimension $d \in \NN$. 

\subsection{Ruelle--Lanford functions}\label{subsection_decoupling_method}

We start this section by presenting the decoupling method (also called the subadditive method) which we will be using in order to prove the existence of a weak LDP for the augmented sequence $(\ell_n,H_n)_{n \in \NN}$. 

For $ \mu \in \mathcal{P}(\mathbb{Z}^d)$, $H \in \mathscr{H}^{(d)}$ and $\varepsilon>0$, define
\begin{equation}\label{eq_def_Omega_n_mu_f}
\begin{aligned}
\Omega_n(\mu,H,\varepsilon) & \eqdef \big\{ w \in \Omega_n^{(0)} : \|\ell_n(w)-\mu\| < \varepsilon, \; d_{\mathscr{H}}(H_n(w,\cdot),H) < \varepsilon \big\} \\
& = (\ell_n,H_n)^{-1} \big( B(\mu,\varepsilon) \times B_{\mathscr{H}}(H,\varepsilon) \big).
\end{aligned}
\end{equation}
Thus, $\Omega_n(\mu,H,\varepsilon)$ consists of trajectories whose empirical measure and profile are both within a distance $\varepsilon$ of their prescribed targets $\mu$ and $H$.

\begin{definition}[Ruelle--Lanford functions]\label{definition_RL}
For $\mu \in \mathcal{P}(\ZZ^d)$ and $H \in \mathscr{H}^{(d)}$, let
\begin{align*}
\underline{s}_{p}(\mu,H) & \eqdef \inf \limits_{\varepsilon>0} \liminf_{n \to \infty} \frac{1}{n} \log \PP_{d,p} \big( (\ell_n,H_n) \in B(\mu,\varepsilon) \times B_{\mathscr{H}}(H,\varepsilon) \big),\\
\text{and} \qquad \quad \overline{s}_{p}(\mu,H) & \eqdef \inf \limits_{\varepsilon>0} \limsup_{n \to \infty} \frac{1}{n} \log \PP_{d,p} \big( (\ell_n,H_n) \in B(\mu,\varepsilon) \times B_{\mathscr{H}}(H,\varepsilon) \big).
\end{align*}
If $\overline{s}_p(\mu,H)=\underline{s}_p(\mu,H)$, we say that the \textit{Ruelle--Lanford function} of $(\ell_n,H_n)_{n \in \NN}$ exists at $(\mu,H)$ and we set
\begin{equation}\label{eq_def_I_RL_couple}
I_{\textup{A}}^{(p)}(\mu,H) \eqdef -\overline{s}_p(\mu,H)=-\underline{s}_p(\mu,H) \in [0,+\infty].
\end{equation}
\end{definition}
\noindent When writing $\overline{s}_p$, $\underline{s}_p$ and $I_{\textup{A}}^{(p)}$ we omit the dependence on the dimension $d \in \NN$, since it is fixed throughout the section. 

Ruelle--Lanford functions were introduced in \cite{Ruelle} as a notion of entropy and later brought into the theory of large deviations in \cite{Lanford}. This led to the development of subadditive methods, first established for independent and identically distributed random variables \cite{Lanford} and \cite{Bahadur}, and subsequently extended to more general settings in \cite{Pfister}, as well as more recently in \cite{Cuneo_selective_decoupling,Daures,Daures_2}. This method also appears in \cite[Theorem~4.1.11]{DZ_LDP}, although the terminology \textit{Ruelle--Lanford} does not appear there. The significance of the decoupling method is made clear by the following proposition.

\begin{proposition}[Weak LDP for the augmented process]\label{prop_RL_functions_imply_LDP}
Suppose that, under $\PP_{d,p}$, the Ruelle--Lanford function of $(\ell_n,H_n)_{n \in \NN}$ exists for all $(\mu,H) \in \mathcal{P}(\ZZ^d) \times \mathscr{H}^{(d)}$. Then, the augmented sequence $(\ell_n,H_n)_{n \in \NN}$ satisfies a weak LDP with rate function $I_{\textup{A}}^{(p)}$.
\end{proposition}

\begin{proof}
The open balls $B(\mu,\varepsilon)$ and $B_{\mathscr{H}}(H,\eta)$ form bases for the topology of weak convergence on $\mathcal{P}(\ZZ^d)$ and the topology of uniform convergence on $\mathscr{H}^{(d)}$, respectively. The products of such balls form a basis for the product topology on $\mathcal{P}(\ZZ^d) \times \mathscr{H}^{(d)}$, see \cite[Theorem~15.1]{Munkres}. Thus, the family of sets $\big\{ B(\mu,\varepsilon)\times B_{\mathscr{H}}(H,\varepsilon) \, : \, \varepsilon>0 \big\}$ also forms a neighbourhood basis at $(\mu,H)$ for the product topology. Indeed, for any $\varepsilon',\eta'>0$ and $0<\varepsilon \leq \min\{\varepsilon',\eta'\}$,
\begin{align*}
    B(\mu,\varepsilon)\times B_{\mathscr{H}}(H,\varepsilon) \subseteq B(\mu,\varepsilon')\times B_{\mathscr{H}}(H,\eta').
\end{align*}
By \cite[Theorem~4.1.11]{DZ_LDP} and Definition~\ref{definition_RL}, the existence of the Ruelle--Lanford function at every $(\mu,H)$ implies that the sequence $(\ell_n,H_n)_{n \in \NN}$ satisfies a weak LDP under $\PP_{d,p}$ with rate function $I_{\text{A}}^{(p)}$, and furthermore, that $I_{\text{A}}^{(p)}$ is lower semicontinuous.
\end{proof}

\subsection{The symmetric case \texorpdfstring{$p=1/(2d)$}{p=1/(2d)}}%
\label{subsection_LDP_symmetric_case}
The aim of this section is to show the existence, under $\PP_{d,1/(2d)}$, of a weak LDP for the augmented sequence $(\ell_n,H_n)_{n \in \NN}$.

\begin{lemma}[The symmetric MERW is a simple random walk]\label{lemma_MERW_p_half_is_SRW}
Under $\PP_{d,1/(2d)}$, the MERW $(S_n)_{n \in \NN}$ has the same law as a symmetric random walk on $\ZZ^d$ started at the origin. In other words, for every $n \in \NN$ and every $w \in \Omega_n^{(0)}$,
\begin{equation}\label{eq_lemma_MERW_p_half_is_SRW}
\PP_{d,\frac{1}{2d}}\big( S_{[1:n]}=w \big) = \left( \frac{1}{2d} \right)^{n-1}. 
\end{equation}
\end{lemma}

\begin{proof}
The case $n=1$ is immediate. Let $n \geq 2$. When $p=1/(2d)$, the expression for $Q_p(s \mid r)$, given in \eqref{eq_def_Q_p_r_s}, yields
\begin{align*}
Q_p(s \mid r)=\frac{1}{2d}, \qquad \forall (s,r) \in \mathscr{S}_d \times \Delta^{(d)}.
\end{align*}
Using Corollary~\ref{corollary_product_of_transitions} and the fact that the initial distribution is the uniform measure $\nu_{\text{unif}}$ on $\mathscr{S}_d$, we have for all $w \in \Omega_n^{(0)}$,
\begin{align*}
\PP_{d,\frac{1}{2d}}\big( S_{[1:n]}=w \big) & = \frac{1}{2d} \prod_{j=2}^{n-1} Q_p \big( \Delta_j(w) \mid \rho_j(w) \big) \\
& = \left(\frac{1}{2d} \right)^{n-1}.
\end{align*}
Hence, the MERW $(S_n)_{n \in \NN}$ has, under $\PP_{d,1/(2d)}$, the same law as the $d$-dimensional SRW.
\end{proof}

\noindent From now on, we write $\PP_{d,\text{sym}}$, $\overline{s}_{\text{sym}}$, $\underline{s}_{\text{sym}}$, $I_{\text{A}}^{\text{sym}}$ and $I_{d,\text{sym}}$ instead of $\PP_{d,1/(2d)}$, $\overline{s}_{1/(2d)}$, $\underline{s}_{1/(2d)}$, $I_{\text{A}}^{(1/(2d))}$ and $I_{d,1/(2d)}$ to mark this correspondence.

\begin{corollary}\label{corollary_identification_in_symmetric_case}
For all $d \in \NN$, 
\begin{equation}\label{eq_corollary_identification_in_symmetric_case}
I_{d,\textup{sym}}=I_{d,\textup{DV}}.
\end{equation}
\end{corollary}

\begin{proof}
By Lemma~\ref{lemma_MERW_p_half_is_SRW}, under $\PP_{d,\text{sym}}$, the process $(S_n)_{n \in \NN}$ has the same law as the simple symmetric random walk on $\ZZ^d$ started at the origin. By \cite[Lemma~5]{Bryc_Dembo}, the empirical measures of the simple symmetric random walk satisfy a weak LDP with the Donsker--Varadhan functional $I_{d,\text{DV}}$ given in \eqref{Donsker_Varadhan_our_case}. Since $\ZZ^d$, endowed with the discrete topology, is a Polish space, $\mathcal{P}(\ZZ^d)$ endowed with the weak topology is also Polish \cite[Theorem~3.2.2]{Bogachev}. The local characterisation of the rate function for a weak LDP, see the remark following \cite[Theorem~4.1.18]{DZ_LDP}, yields
\begin{align*}
I_{d,\text{DV}}(\mu)=-\inf_{\varepsilon>0} \limsup_{n \to \infty} \frac{1}{n} \log \PP_{d,\text{sym}}\big( \ell_n \in B(\mu,\varepsilon) \big), \qquad \mu \in \mathcal{P}(\ZZ^d).
\end{align*}
By definition of $I_{d,p}$ given in \eqref{expression_rate_function}, we obtain \eqref{eq_corollary_identification_in_symmetric_case}. 
\end{proof}

\noindent The following proposition provides the stitching estimate for
the trajectory sets defined in \eqref{eq_def_Omega_n_mu_f}.

\begin{proposition}[Controlled stitching procedure]\label{prop_good_stitching_procedure}
Fix $(\mu_1,H^{(1)}),(\mu_2,H^{(2)})\in \mathcal{P}(\ZZ^d) \times \mathscr{H}^{(d)}$. For every $\varepsilon>0$ and $\kappa>0$, there exists $\varepsilon'>0$ for which the following holds. There exists $n_0 \in \NN$ and a map $N_0: \{ n\in \NN \, : \, n \geq n_0\} \to \NN$ such that, for every $n \geq n_0$ and $N \geq N_0(n)$, one can choose an integer $k=k(n,N)\in \NN$ satisfying $2kn \leq N$, constants $C_{n,N}>0$ and $D_{n,N}\geq 1$, and a map
\begin{equation*}
\psi : \mathcal{W}_{n,k} \longrightarrow
\Omega_N,
\end{equation*}
where $\mathcal{W}_{n,k}=\big(\Omega_n(\mu_1,H^{(1)},\varepsilon') \times \Omega_n(\mu_2,H^{(2)},\varepsilon')\big)^k$,
so that the following four properties hold.
\begin{enumerate}
\item The image of the smaller target trajectories is contained in the larger target set:
\begin{equation}\label{eq_objectif_inclusion}
\qquad\qquad \qquad \qquad \qquad \psi(\mathcal{W}_{n,k}) \subseteq \Omega_N\left(\frac{\mu_1+\mu_2}{2},\frac{H^{(1)}+H^{(2)}}{2},\varepsilon \right). \qquad \qquad \quad \; \text{\footnotesize[inclusion]}
\end{equation}

\item For every $\overline{w}\in\Omega_N$,
\begin{equation}\label{eq_objectif_cardinality_bound}
\qquad \qquad \qquad\qquad \qquad \qquad \qquad \qquad\#\psi^{-1}(\{\overline{w}\})
\leq
D_{n,N}. \qquad \qquad \quad  \text{\footnotesize[combinatorial bound]}
\end{equation}

\item For all
$\big((w^{(1)},w^{(2)}),\dots,(w^{(2k-1)},w^{(2k)})\big) = (w^{(1)},\dots,w^{(2k)})\in \mathcal{W}_{n,k}$,
\begin{equation}\label{eq_objectif_comp_proba}
\quad \PP_{d,\textup{sym}}\Big(S_{[1:N]} =\psi(w^{(1)},\dots,w^{(2k)})\Big) \geq C_{n,N} \prod_{i=1}^{2k} \PP_{d,\textup{sym}}\big(S_{[1:n]}=w^{(i)}\big). \quad \text{\footnotesize[lower decoupling]}
\end{equation}

\item The constants $C_{n,N}$ and $D_{n,N}$ satisfy 
\begin{equation}\label{eq_objectif_asympto}
\qquad \qquad \qquad \qquad \qquad \qquad \quad \liminf_{n\to\infty} \liminf_{N\to\infty} \frac{1}{N} \log \frac{C_{n,N}}{D_{n,N}} \geq -\kappa.  \qquad \qquad  \text{\footnotesize[asymptotic bound]}
\end{equation}
\end{enumerate}
\end{proposition}

\noindent The proof of Proposition~\ref{prop_good_stitching_procedure}, which contains the construction of $\psi$ and the verification of the four properties above, is relatively long and technical. It is postponed to Section~\ref{section_proof_of_existence_of_stitching_procedure}. The function $\psi$ whose existence is given by Proposition~\ref{prop_good_stitching_procedure} will be called the \textit{stitching procedure}. Note that there are in fact multiple such stitching procedures since $\psi$ depends on the different parameters introduced in this proposition. 

The following lemma uses this stitching procedure to obtain two key inequalities for the existence and convexity of RL functions.

\begin{lemma}[Existence of RL functions and midpoint concavity]\label{lemma_intermediary_inequality_existence}
Under $\PP_{d,\textup{sym}}$, for all $(\mu,H) \in \mathcal{P}(\ZZ^d) \times \mathscr{H}^{(d)}$, 
\begin{equation}\label{eq_lemma_intermediary_inequality_existence}
\underline{s}_{\textup{sym}}(\mu,H)=\overline{s}_{\textup{sym}}(\mu,H).
\end{equation}
In addition, for all $(\mu_1,H^{(1)}),(\mu_2,H^{(2)}) \in \mathcal{P}(\ZZ^d) \times \mathscr{H}^{(d)}$,
\begin{equation}\label{eq_lemma_intermediary_inequality_convexity}
\underline{s}_{\textup{sym}}\left(\frac{\mu_1+\mu_2}{2},\frac{H^{(1)}+H^{(2)}}{2} \right) \geq \frac{1}{2}\overline{s}_{\textup{sym}}(\mu_1,H^{(1)}) + \frac{1}{2}\overline{s}_{\textup{sym}}(\mu_2,H^{(2)}).
\end{equation}
\end{lemma}

\begin{proof}
Fix $(\mu_1,H^{(1)}),(\mu_2,H^{(2)}) \in \mathcal{P}(\ZZ^d) \times \mathscr{H}^{(d)}$, $\varepsilon>0$ and $\kappa>0$. Choose $\varepsilon'>0$, $n_0 \in \NN$ and $N_0(n) \in \NN$ as in Proposition~\ref{prop_good_stitching_procedure}. Let $n \geq n_0$, $N \geq N_0(n)$ and consider $k(n,N) \in \NN$, $C_{n,N}>0$, $D_{n,N} \geq 1$ and $\psi$ again as in Proposition~\ref{prop_good_stitching_procedure}. By the inclusion given in \eqref{eq_objectif_inclusion},
\begin{align*}
\PP_{d,\text{sym}}\left(\Omega_N\left(\frac{\mu_1+\mu_2}{2},\frac{H^{(1)}+H^{(2)}}{2},\varepsilon \right)\right) &\geq \sum_{\overline{w} \in \psi(\mathcal{W}_{n,k})} \PP_{d,\text{sym}}(S_{[1:N]}=\overline{w}).
\end{align*}
Using the combinatorial bound, given in \eqref{eq_objectif_cardinality_bound} and then the lower decoupling property, given in \eqref{eq_objectif_comp_proba}, we obtain for all $\overline{w} \in \psi(\mathcal{W}_{n,k})$,
\begin{align*}
D_{n,N}\PP_{d,\text{sym}}\big(S_{[1:N]}=\overline{w}\big)& \geq \sum_{\substack{(w^{(1)},\dots,w^{(2k)}) \in \mathcal{W}_{n,k} \, : \\ \psi(w^{(1)},\dots,w^{(2k)})=\overline{w}}} \PP_{d,\text{sym}}\big(S_{[1:N]}=\psi(w^{(1)},\dots,w^{(2k)})\big) \\
& \geq C_{n,N} \sum_{\substack{(w^{(1)},\dots,w^{(2k)}) \in \mathcal{W}_{n,k} \, : \\ \psi(w^{(1)},\dots,w^{(2k)})=\overline{w}}} \prod_{i=1}^{2k} \PP_{d,\text{sym}} \big(S_{[1:n]}=w^{(i)}\big).
\end{align*}
By summing the above inequality over all $\overline{w} \in \psi(\mathcal{W}_{n,k})$, we obtain
\begin{align}
& \PP_{d,\text{sym}}\left(\Omega_N\left(\frac{\mu_1+\mu_2}{2},\frac{H^{(1)}+H^{(2)}}{2},\varepsilon\right)\right) \notag\\
& \qquad \qquad \qquad \geq \frac{C_{n,N}}{D_{n,N}} \sum_{(w^{(1)},\dots,w^{(2k)}) \in \mathcal{W}_{n,k} } \prod_{i=1}^{2k} \PP_{d,\text{sym}}(S_{[1:n]}=w^{(i)}) \notag \\
& \qquad \qquad \qquad = \frac{C_{n,N}}{D_{n,N}} \PP_{d,\text{sym}}\big(\Omega_n(\mu_1,H^{(1)},\varepsilon')\big)^k \, \PP_{d,\text{sym}}\big(\Omega_n(\mu_2,H^{(2)},\varepsilon')\big)^k, \label{proof_prop_RL_equality_1}
\end{align}
where the last equality follows from the definition of $\mathcal{W}_{n,k}$. Taking logarithms, dividing by $N$, and using the fact that $2kn\leq N$ with $\log \PP_{d,\text{sym}}(\cdot)\leq 0$, yields
\begin{align*}
&\frac{1}{N}\log \PP_{d,\text{sym}}\left(\Omega_N\left(\frac{\mu_1+\mu_2}{2},\frac{H^{(1)}+H^{(2)}}{2},\varepsilon\right)\right) \\
& \qquad \quad \geq \frac{1}{N}\log\frac{C_{n,N}}{D_{n,N}} +\frac{1}{2n} \log \PP_{d,\text{sym}}\big(\Omega_n(\mu_1,H^{(1)},\varepsilon')\big) + \frac{1}{2n} \log \PP_{d,\text{sym}}\big(\Omega_n(\mu_2,H^{(2)},\varepsilon')\big).
\end{align*}
Since the above inequality holds for all $N \geq N_0(n)$, we can take $\liminf_{N\to\infty}$ in the above equation. Then, using the asymptotic bound given in \eqref{eq_objectif_asympto}, we deduce that for every $\delta>0$ and for all sufficiently large $n$,
\begin{align}
& \liminf_{N\to\infty} \frac{1}{N}\log \PP_{d,\text{sym}}\left(\Omega_N\left(\frac{\mu_1+\mu_2}{2},\frac{H^{(1)}+H^{(2)}}{2},\varepsilon \right)\right) \label{proof_prop_RL_equality_2} \\
& \qquad \qquad \qquad \geq -\kappa -\delta + \frac{1}{2n} \log \PP_{d,\text{sym}}\big(\Omega_n(\mu_1,H^{(1)},\varepsilon')\big) +  \frac{1}{2n} \log \PP_{d,\text{sym}}\big(\Omega_n(\mu_2,H^{(2)},\varepsilon')\big). \notag
\end{align}
Suppose now that $\mu_1=\mu=\mu_2$ and $H^{(1)}=H=H^{(2)}$. Taking $\limsup_{n\to\infty}$ in \eqref{proof_prop_RL_equality_2} and rewriting the events $\Omega_N\left(\mu,H,\varepsilon\right)$ and $\Omega_n(\mu,H,\varepsilon')$ with \eqref{eq_def_Omega_n_mu_f}, we deduce that for every $\delta>0$,
\begin{align*}
& \liminf_{N\to\infty} \frac{1}{N}\log \PP_{d,\text{sym}}\left((\ell_N,H_N) \in B\left(\mu,\varepsilon \right) \times B_{\mathscr{H}}(H,\varepsilon)\right) \\
& \qquad \qquad \qquad \qquad \qquad \qquad \geq -\kappa -\delta + \limsup_{n\to\infty} \frac{1}{n} \log \PP_{d,\text{sym}}\big((\ell_n,H_n) \in B(\mu,\varepsilon') \times B_{\mathscr{H}}(H,\varepsilon')\big).
\end{align*}
Since $B(\mu,\varepsilon') \times B_{\mathscr{H}}(H,\varepsilon')$ is one of the neighbourhoods over which $\overline{s}_{\textup{sym}}(\mu,H)$ is defined, the right-hand side is bounded from below by $-\kappa-\delta+\overline{s}_{\text{sym}}(\mu,H)$. Hence, for every $\varepsilon>0$,
\begin{equation*}
\liminf_{N\to\infty} \frac{1}{N}\log \PP_{d,\text{sym}}\big((\ell_N,H_N) \in B(\mu,\varepsilon) \times B_{\mathscr{H}}(H,\varepsilon)\big) \geq -\kappa-\delta+\overline{s}_{\text{sym}}(\mu,H).
\end{equation*}
Taking the infimum over $\varepsilon>0$, and then letting $\kappa \to 0$ and $\delta \to 0$, we obtain $\underline{s}_{\text{sym}}(\mu,H) \geq \overline{s}_{\text{sym}}(\mu,H)$. 
The reverse inequality is immediate from the definitions of $\overline{s}_{\text{sym}}$ and $\underline{s}_{\text{sym}}$, which proves \eqref{eq_lemma_intermediary_inequality_existence}.

To show \eqref{eq_lemma_intermediary_inequality_convexity}, consider again $(\mu_1,H^{(1)}),(\mu_2,H^{(2)}) \in \mathcal{P}(\ZZ^d) \times \mathscr{H}^{(d)}$. Taking $\liminf_{n \to \infty}$ in \eqref{proof_prop_RL_equality_2}, we obtain
\begin{align*}
& \liminf_{N\to\infty} \frac{1}{N}\log \PP_{d,\text{sym}}\left((\ell_N,H_N) \in B\left(\frac{\mu_1+\mu_2}{2},\varepsilon \right) \times B_{\mathscr{H}} \left( \frac{H^{(1)}+H^{(2)}}{2},\varepsilon \right) \right) \\
& \qquad \qquad \qquad \qquad \quad \geq -\kappa -\delta + \frac{1}{2} \liminf_{n\to\infty} \frac{1}{n} \log \PP_{d,\text{sym}}\big((\ell_n,H_n) \in B(\mu_1,\varepsilon') \times B_{\mathscr{H}}(H^{(1)},\varepsilon')\big) \\
& \qquad \qquad \qquad \qquad \quad \quad + \frac{1}{2} \liminf_{n\to\infty} \frac{1}{n} \log \PP_{d,\text{sym}}\big((\ell_n,H_n) \in B(\mu_2,\varepsilon') \times B_{\mathscr{H}}(H^{(2)},\varepsilon')\big).
\end{align*}
By definition of $\underline{s}_{\text{sym}}$, each local $\liminf_{n \to \infty}$ on the right is greater than or equal to the corresponding $\underline{s}_{\text{sym}}(\mu_i,H^{(i)})$. Using the equality of lower and upper RL functions, given in \eqref{eq_lemma_intermediary_inequality_existence}, we obtain
\begin{align*}
& \liminf_{N\to\infty} \frac{1}{N}\log \PP_{d,\text{sym}}\left((\ell_N,H_N) \in B\left(\frac{\mu_1+\mu_2}{2},\varepsilon \right) \times B_{\mathscr{H}}\left( \frac{H^{(1)}+H^{(2)}}{2},\varepsilon \right)\right) \\
& \qquad \qquad \qquad \qquad \qquad \qquad \qquad \qquad \qquad \qquad  \geq -\kappa - \delta +\frac{1}{2} \overline{s}_{\text{sym}}(\mu_1,H^{(1)}) +\frac{1}{2} \overline{s}_{\text{sym}}(\mu_2,H^{(2)}).
\end{align*}
Taking the infimum over $\varepsilon>0$, and then letting $\kappa \to 0$ and $\delta \to 0$, we obtain \eqref{eq_lemma_intermediary_inequality_convexity} which concludes the proof.
\end{proof}

\begin{proposition}[Weak LDP in the symmetric case]\label{prop_good_stitching_implies_RL_equality}
Under $\PP_{d,\textup{sym}}$, the sequence $(\ell_n,H_n)_{n \in \NN}$ satisfies a weak LDP in $\mathcal{P}(\ZZ^d) \times \mathscr{H}^{(d)}$, with rate function $I_{\textup{A}}^{\textup{sym}}$ given by \eqref{eq_def_I_RL_couple}. Furthermore, the rate function $I_{\textup{A}}^{\textup{sym}}$ is convex on $\mathcal{P}(\ZZ^d) \times \mathscr{H}^{(d)}$.
\end{proposition}

\begin{proof}
By \eqref{eq_lemma_intermediary_inequality_existence} of Lemma~\ref{lemma_intermediary_inequality_existence}, the Ruelle--Lanford function exists at every $(\mu,H) \in \mathcal{P}(\ZZ^d) \times \mathscr{H}^{(d)}$. Proposition~\ref{prop_RL_functions_imply_LDP} therefore yields the weak LDP for the sequence $(\ell_n,H_n)_{n \in \NN}$, with rate function $I_{\text{A}}^{\text{sym}}=-\overline{s}_{\text{sym}}=-\underline{s}_{\text{sym}}$.

Let us now show that $I_{\text{A}}^{\text{sym}}$ is convex. Using \eqref{eq_lemma_intermediary_inequality_convexity} of Lemma~\ref{lemma_intermediary_inequality_existence} and the preceding discussion, 
\begin{align*}
I_{\text{A}}^{\text{sym}}\left( \frac{\mu_1+\mu_2}{2}, \frac{H^{(1)}+H^{(2)}}{2} \right) \leq \frac{1}{2} I_{\text{A}}^{\text{sym}}\big(\mu_1,H^{(1)}\big) + \frac{1}{2} I_{\text{A}}^{\text{sym}}\big(\mu_2,H^{(2)}\big).
\end{align*}
By iterating this inequality, we obtain
\begin{align*}
I_{\text{A}}^{\text{sym}}\big( \lambda \mu_1 + (1-\lambda)\mu_2,\lambda H^{(1)} + (1-\lambda)H^{(2)} \big) \leq \lambda I_{\text{A}}^{\text{sym}}\big(\mu_1,H^{(1)}\big) + (1-\lambda) I_{\text{A}}^{\text{sym}}\big(\mu_2,H^{(2)}\big),
\end{align*}
for all $\lambda=m/2^q$, where $q \in \NN$ and $0 \leq m \leq 2^q$. Using the density of the dyadics and the lower semicontinuity of $I_{\text{A}}^{\text{sym}}$, we extend the above relation to all $\lambda \in [0,1]$ and therefore $I_{\text{A}}^{\text{sym}}$ is convex.
\end{proof}

\noindent The following section extends the existence results obtained in Proposition~\ref{prop_good_stitching_implies_RL_equality} to any $p \in (0,1)$.

\subsection{The general case}\label{subsection_comparaison_RL_functions}

Proposition~\ref{prop_good_stitching_implies_RL_equality} shows that under $\PP_{d,\text{sym}}$, the augmented sequence $(\ell_n,H_n)_{n \in \NN}$ satisfies a weak LDP. In this subsection, we extend this result to all $p \in (0,1)$ by means of exponential tilting. Before stating the main result of this subsection, we introduce the \textit{profile-energy} map $\mathcal{E}_p : \mathscr{H}^{(d)} \longrightarrow \RR$, defined by setting
\begin{equation}\label{eq_def_energy_profile}
\mathcal{E}_p(H) \eqdef \int_0^1 \sum_{s \in \mathscr{S}_d} [H'(t)]_s \log \left( \frac{Q_p(s \mid H(t)/t)}{1/(2d)} \right) \dd t, \qquad H \in \mathscr{H}^{(d)}.
\end{equation}
By definition of $\mathscr{H}^{(d)}$, for each $t \in (0,1]$, the vector $H(t)/t$ belongs to $\Delta^{(d)}$. The integrand may be assigned any value at $t=0$ without affecting the integral. The integral is well defined since, by Lemma~\ref{lemma_properties_profile_space}, $H'(t) \in \Delta^{(d)}$ for a.e. $t \in [0,1]$ and, by Lemma~\ref{lemma_lower_bound_transition}, the logarithm in \eqref{eq_def_energy_profile} is uniformly bounded. 

We will show that the profile-energy map quantifies, at exponential scale, the effect of the memory parameter on trajectory probabilities. More precisely, for an event $A_n \subseteq \Omega_n^{(0)}$ whose trajectories have profiles uniformly close to a given $H \in \mathscr{H}^{(d)}$, we will establish a comparison of the form
\begin{align*}
    \PP_{d,p}\big(A_n\big) \simeq \PP_{d,\text{sym}}\big(A_n\big)e^{n \mathcal{E}_p(H)}.
\end{align*}
Thus, for events localised near $H$, the memory modifies their exponential cost by $-\mathcal{E}_p(H)$. This comparison allows us to transfer the local large-deviation estimates for the augmented sequence $(\ell_n,H_n)_{n \in \NN}$ from the symmetric random walk to the MERW with parameter $p$. A precise statement of this comparison is given in Corollary~\ref{coro_tilting_ineq_for_localised_profile}.

\begin{remark}\label{remark_energy_profile_constant}
Fix $r=(r_s)_{s \in \mathscr{S}_d} \in \Delta^{(d)}$ and the linear profile $H \in \mathscr{H}^{(d)}$ defined by $H(t)=tr$. Then $H'(t)=r$ and $H(t)/t=r$ for all $t \in (0,1]$. Thus,
\begin{equation}\label{eq_remark_energy_profile_constant}
\mathcal{E}_p(H)=\sum_{s \in \mathscr{S}_d} r_s \log \left( \frac{Q_p(s \mid r)}{1/(2d)} \right).
\end{equation}
\end{remark}

\noindent The following is the main proposition of this section.

\begin{proposition}[Weak LDP for general $p$]\label{proposition_expression_RL_rate_function}
Fix $p \in (0,1)$. Under $\PP_{d,p}$, the sequence $(\ell_n,H_n)_{n \in \NN}$ satisfies a weak LDP with rate function $I_{\textup{A}}^{(p)}=-\underline{s}_p=-\overline{s}_p$. Furthermore, for every $(\mu,H) \in \mathcal{P}(\ZZ^d) \times \mathscr{H}^{(d)}$,
\begin{equation}\label{eq_proposition_expression_RL_rate_function}
I_{\textup{A}}^{(p)}(\mu,H) = I_{\textup{A}}^{\textup{sym}}(\mu,H)-\mathcal{E}_p(H).
\end{equation}
\end{proposition}

\begin{remark}\label{remark_I_RL_p_when_p_is_1_slash_2d}
When $p=1/(2d)$, one has $Q_p(s \mid r)=1/(2d)$ for every $(s,r)\in \mathscr{S}_d \times \Delta^{(d)}$. Using \eqref{eq_def_energy_profile}, we obtain $\mathcal{E}_{1/(2d)} \equiv 0$. Equation~\eqref{eq_proposition_expression_RL_rate_function} then yields $I_{\textup{A}}^{(1/(2d))}=I_{\textup{A}}^{\textup{sym}}$ as required.
\end{remark}

\begin{remark}\label{remark_decomposition_E_p_f}
In the expression given in \eqref{eq_proposition_expression_RL_rate_function}, the term $I_{\text{A}}^{\text{sym}}(\mu,H)$ can be seen as a combinatorial term, which reflects the number of trajectories whose empirical measure approximates $\mu$ and whose profile approximates $H$. The term $\mathcal{E}_p(H)$ on the other hand captures the modifications of the probabilities induced by the memory $p$ and only depends on the profile $H$, see Remark~\ref{remark_interpretation_energy_profile}.
\end{remark}

We will use the following elementary result in the proofs of Lemmas~\ref{lemma_exponential_tilting_wrt_SRW} and~\ref{lemma_prefix_does_not_change}.

\begin{lemma}[Uniform bounds and Lipschitz estimates for $Q_p$]\label{lemma_properties_log_Q_p}
There exist constants $M,M'>0$, depending only on $d$ and $p$, such that, for all $s \in \mathscr{S}_d$ and $r,\widetilde{r} \in \Delta^{(d)}$,
\begin{align}
\big|\log \big(2d Q_p(s \mid r)\big) \big| & \leq M, \label{eq_lemma_properties_log_Q_p_1} \\
\big|\log \big(2d Q_p(s \mid r)\big) - \log \big(2d Q_p(s \mid \widetilde{r})\big) \big| & \leq M' \|r - \widetilde{r} \|_{\infty}. \label{eq_lemma_properties_log_Q_p_2}
\end{align}
\end{lemma}

\begin{proof}
By the definition of $Q_p$ and the uniform bound $p_*\leq Q_p(s\mid r)\leq 1$, given in Lemma~\ref{lemma_lower_bound_transition}, the map $r \mapsto  \log(2dQ_p(s \mid r))$ is uniformly bounded on $\Delta^{(d)}$, which shows \eqref{eq_lemma_properties_log_Q_p_1}. Moreover, for $r,\widetilde r\in\Delta^{(d)}$,
\begin{align*}
\big|Q_p(s \mid r)-Q_p(s \mid \widetilde{r})\big|
=\left|p-\frac{1-p}{2d-1}\right|\,|r_s-\widetilde r_s|
\leq \left|p-\frac{1-p}{2d-1}\right|\,\|r-\widetilde r\|_\infty.
\end{align*}
Since $x \mapsto \log x$ is $1/p_*$-Lipschitz on $[p_*,1]$, we obtain \eqref{eq_lemma_properties_log_Q_p_2}.
\end{proof}

We now show that $\mathcal{E}_p$ gives the exponential density of the MERW with respect to the symmetric random walk.

\begin{lemma}[Profile-dependent exponential tilting]\label{lemma_exponential_tilting_wrt_SRW}
There exists a constant $M_{d,p}>0$ such that for every $n \geq 3$ and every $w \in \Omega_n^{(0)}$, 
\begin{equation}\label{eq_lemma_exponential_tilting_wrt_SRW}
\left \vert \frac{1}{n} \log \frac{\PP_{d,p}\left( S_{[1:n]}=w \right)}{\PP_{d,\textup{sym}}\big( S_{[1:n]}=w \big)} - \mathcal{E}_p\big(H_n(w)\big) \right \vert < M_{d,p} \frac{\log n}{n}. 
\end{equation}
\end{lemma}

\begin{proof}
Fix $n \geq 3$ and $w \in \Omega_n^{(0)}$. For $k \in \dbl 0,n-1 \dbr$, let $t_k \eqdef k/(n-1)$. Applying Corollary~\ref{corollary_product_of_transitions} under $\PP_{d,p}$ and $\PP_{d,\textup{sym}}$, we cancel the common first-step probability $1/(2d)$ and obtain
\begin{equation}\label{proof_lemma_exponential_tilting_wrt_SRW_3}
\frac{1}{n} \log \frac{\PP_{d,p}\left( S_{[1:n]}=w \right)}{\PP_{d,\textup{sym}}\big( S_{[1:n]}=w \big)} = \frac{1}{n} \sum_{k=2}^{n-1} \log \Big( 2d Q_p\big(\Delta_k(w) \mid \rho_k(w) \big) \Big).
\end{equation}
Fix $k \in \dbl 2,n-1 \dbr$ and $t \in (t_{k-1},t_k)$. Equation~\eqref{eq_lemma_properties_profile_space_derivative} gives $H_n'(w,t)=\delta_{\Delta_k(w)}$. By construction, $H_n(w,\cdot)$ is affine on $[t_{k-1},t_k]$. Moreover, by \eqref{eq_def_profile_map} and \eqref{eq_def_empirical_distribution_rho}, $H_n(w,t_{k-1})=t_{k-1} \rho_k(w)$. Thus,
\begin{align*}
\frac{H_n(w,t)}{t} & = \frac{t_{k-1}}{t}\frac{H_n(w,t_{k-1})}{t_{k-1}} + H_n'(w,t) \frac{t-t_{k-1}}{t} \\
& = \frac{t_{k-1}}{t} \rho_{k}(w) + \frac{t-t_{k-1}}{t} \delta_{\Delta_k(w)}.
\end{align*}
Using $\| \delta_{\Delta_k(w)} - \rho_k(w) \|_{\infty} \leq 1$ and $t>t_{k-1}=(k-1)/(n-1)$, we obtain
\begin{equation}\label{proof_lemma_exponential_tilting_wrt_SRW_4}
\left\|\frac{H_n(w,t)}{t}- \rho_k(w) \right\|_{\infty} \leq \frac{1}{k-1}.
\end{equation}
By definition of $\mathcal{E}_p$ and the identity $[H_n'(w,t)]_s=\mathbbm{1}_{\{\Delta_k(w)=s\}}$ on $(t_{k-1},t_k)$, 
\begin{align*}
\mathcal{E}_p\big(H_n(w,\cdot)\big)= \sum_{k=1}^{n-1} \int_{t_{k-1}}^{t_k}\log \left( 2dQ_p\left(\Delta_k(w) \biggm| \frac{H_n(w,t)}{t}\right)\right) \dd t.
\end{align*}
By \eqref{eq_lemma_properties_log_Q_p_1}, the contribution of the first interval $(t_0,t_1)$ has absolute value at most $M/(n-1)$. For every $k \geq 2$, applying \eqref{eq_lemma_properties_log_Q_p_2} and \eqref{proof_lemma_exponential_tilting_wrt_SRW_4} gives
\begin{align*}
& \left \vert \mathcal{E}_p\big(H_n(w,\cdot)\big) - \frac{1}{n-1} \sum_{k=2}^{n-1} \log \Big( 2d Q_p\big(\Delta_k(w) \mid \rho_k(w) \big) \Big) \right \vert  \\
& \qquad \qquad \qquad \qquad \qquad \qquad \qquad \qquad \leq \frac{M}{n-1} + \sum_{k=2}^{n-1} \int_{t_{k-1}}^{t_k} M' \left\| \frac{H_n(w,t)}{t}-\rho_k(w) \right\|_{\infty} \dd t \\
& \qquad \qquad \qquad \qquad \qquad \qquad \qquad \qquad \leq \frac{M}{n-1} + \frac{M'}{n-1} \sum_{k=2}^{n-1} \frac{1}{k-1}.
\end{align*}
It remains only to replace the normalisation $1/(n-1)$ by $1/n$. By \eqref{eq_lemma_properties_log_Q_p_1}, 
\begin{align*}
\left \vert \frac{1}{n} \sum_{k=2}^{n-1} \log \Big( 2d Q_p\big(\Delta_k(w) \mid \rho_k(w) \big) \Big) - \frac{1}{n-1} \sum_{k=2}^{n-1} \log \Big( 2d Q_p\big(\Delta_k(w) \mid \rho_k(w) \big) \Big) \right \vert \leq \frac{M}{n-1}.
\end{align*}
Combining this estimate with \eqref{proof_lemma_exponential_tilting_wrt_SRW_3} yields
\begin{align*}
\left \vert \frac{1}{n} \log \frac{\PP_{d,p}\left( S_{[1:n]}=w \right)}{\PP_{d,\textup{sym}}\big( S_{[1:n]}=w \big)} - \mathcal{E}_p\big(H_n(w)\big) \right \vert \leq 2 \frac{M}{n-1} + \frac{M'}{n-1} \sum_{k=1}^{n-2} \frac{1}{k}.
\end{align*}
Finally, for $n \geq 3$, $\sum_{k=1}^{n-2} 1/k \leq 1+\log(n-2) < 2 \log n$ and $1/(n-1)\leq 2/n$. Hence, the right-hand side is at most $4(M+M')\log(n)/n$. Taking $M_{d,p} \eqdef 4(M+M')>0$, which depends only on $d$ and $p$, we obtain \eqref{eq_lemma_exponential_tilting_wrt_SRW}.
\end{proof}

We next derive an alternative, derivative-free representation of $\mathcal{E}_p$, based on an integration-by-parts. This representation will be useful in establishing the continuity of $\mathcal{E}_p$ in Lemma~\ref{lemma_continuity_of_E_p}, and its concavity for $p \leq 1/(2d)$ in Lemma~\ref{lemma_concavity_of_E_p}. We start by introducing the scalar functions $L_p,F_p,G_p:[0,1] \to \RR$ defined by
\begin{align}
L_p(x) & \eqdef \log \left(2d\left(px + \frac{1-p}{2d-1}(1-x)\right)\right), \label{eq_def_L_p}\\
F_p(x) & \eqdef \int_0^x L_p(t) \mathrm{d}t, \label{eq_def_F_p}\\
G_p(x) & \eqdef xL_p(x)-F_p(x). \label{eq_def_G_p}
\end{align}
Note that if $r \in \Delta^{(d)}$,
\begin{align*}
L_p(r_s)=\log \left( \frac{Q_p(s \mid r)}{1/(2d)} \right), \qquad s \in \mathscr{S}_d.
\end{align*}
We know that $0<p_* \leq px + \frac{1-p}{2d-1}(1-x)$ and therefore $L_p$, $F_p$ and $G_p$ are well defined and twice differentiable on $[0,1]$. The following identity gives another representation of the profile-energy map defined in \eqref{eq_def_energy_profile}. 
\begin{lemma}[Derivative-free representation of the profile-energy]\label{lemma_analytic_representation_of_profile_energy_map}
For every $H \in \mathscr{H}^{(d)}$, 
\begin{equation}\label{eq_lemma_analytic_representation_of_profile_energy_map}
\mathcal{E}_p(H)=\sum_{s \in \mathscr{S}_d} \left( F_p \big( [H(1)]_s \big) + \int_0^1 G_p\left( \frac{[H(t)]_s}{t} \right) \mathrm{d}t \right).
\end{equation}
\end{lemma}

\begin{proof}
Fix $s \in \mathscr{S}_d$ and put $r_s(t)=[H(t)]_s/t$ for $t \in (0,1]$. For every $\delta \in (0,1)$, by item~2 of Lemma~\ref{lemma_properties_profile_space}, the function $r_s$ is absolutely continuous on $[\delta,1]$, and
\begin{align*}
[H'(t)]_s=r_s(t)+t r'_s(t), \qquad \text{for a.e. } t \in [\delta,1].
\end{align*}
By definition of $F_p$, we have $F_p'=L_p$. Combining this with an integration by parts for absolutely continuous functions and the definition of $G_p$, we obtain
\begin{align*}
\mathcal{E}_p(H) & = \sum_{s \in \mathscr{S}_d} \left\{ \int_{\delta}^1 [H'(t)]_s L_p \big(r_s(t)\big) \dd t + \int_{0}^{\delta} [H'(t)]_s L_p \big(r_s(t)\big) \dd t \right\} \\
& = \sum_{s \in \mathscr{S}_d} \left\{ \int_{\delta}^1 r_s(t) L_p\big(r_s(t)\big) \dd t + \int_{\delta}^1 t \frac{\dd}{\dd t} F_p\big(r_s(t)\big) \dd t + \int_{0}^{\delta} [H'(t)]_s L_p \big(r_s(t)\big) \dd t \right\} \\
& = \sum_{s \in \mathscr{S}_d} \left\{ F_p \big( [H(1)]_s \big) - \delta F_p\big( r_s(\delta) \big) + \int_{\delta}^1 G_p \big( r_s(t) \big) \dd t + \int_{0}^{\delta} [H'(t)]_s L_p \big(r_s(t)\big) \dd t \right\}.
\end{align*}
Since $L_p$, $F_p$ and $G_p$ are bounded on $[0,1]$, letting $\delta \to 0$ yields \eqref{eq_lemma_analytic_representation_of_profile_energy_map} after summing over $s \in \mathscr{S}_d$.
\end{proof}

\begin{lemma}[Continuity of the profile-energy map]\label{lemma_continuity_of_E_p}
The map $\mathcal{E}_p : \mathscr{H}^{(d)} \longrightarrow \RR$ is continuous for the uniform topology.
\end{lemma}

\begin{proof}
Let $(H_m)_{m \in \NN}$ be a sequence in $\mathscr{H}^{(d)}$ converging uniformly to $H \in \mathscr{H}^{(d)}$. Then, $H_m(1) \to H(1)$, and hence, by continuity of $F_p$,
\begin{equation}\label{proof_lemma_continuity_of_E_p_1}
\sum_{s \in \mathscr{S}_d} F_p\big([H_m(1)]_s\big) \xrightarrow[]{m \to \infty} \sum_{s \in \mathscr{S}_d} F_p\big([H(1)]_s\big).
\end{equation}
For every fixed $t \in (0,1]$, we also have $H_m(t)/t \to H(t)/t$, and hence, by continuity of $G_p$,
\begin{align*}
G_p\left(\frac{[H_m(t)]_s}{t}\right) \xrightarrow[]{m \to \infty} G_p\left(\frac{[H(t)]_s}{t}\right).
\end{align*}
Moreover $H_m(t)/t$ and $H(t)/t$ belong to $\Delta^{(d)}$ for all $t \in (0,1]$, and $G_p$ is bounded on $[0,1]$. The dominated convergence theorem therefore gives, for every $s \in \mathscr{S}_d$,
\begin{equation}\label{proof_lemma_continuity_of_E_p_2}
\int_0^1 G_p\left(\frac{[H_m(t)]_s}{t}\right) \dd t \xrightarrow[]{m \to \infty} \int_0^1 G_p\left(\frac{[H(t)]_s}{t}\right) \dd t.
\end{equation}
Using the second representation of the profile-energy map, given in \eqref{eq_lemma_analytic_representation_of_profile_energy_map}, together with \eqref{proof_lemma_continuity_of_E_p_1} and \eqref{proof_lemma_continuity_of_E_p_2} gives $\mathcal{E}_p(H_m) \to \mathcal{E}_p(H)$ and the lemma is proved.
\end{proof}

\begin{remark}[Interpretation of the profile-energy map] \label{remark_interpretation_energy_profile}
Lemmas~\ref{lemma_exponential_tilting_wrt_SRW} and \ref{lemma_continuity_of_E_p} show that $\mathcal{E}_p(H)$ is the asymptotic log-likelihood advantage of the MERW with memory parameter $p$ over the symmetric random walk along the cumulative increment profile $H$. More precisely, if $w^{(n)} \in \Omega_n^{(0)}$ is any sequence of trajectories such that $H_n(w^{(n)})$ converges uniformly to $H \in \mathscr{H}^{(d)}$, then
\begin{align*}
\lim_{n \to \infty} \frac{1}{n} \log\frac{\PP_{d,p}\big(S_{[1:n]}=w^{(n)}\big)}{\PP_{d,\textup{sym}}\big(S_{[1:n]}=w^{(n)}\big)}=\mathcal{E}_p(H).
\end{align*}
Consequently, if $\mathcal{E}_p(H)>0$, trajectories whose increment profiles are close to $H$ are exponentially more likely under the MERW than under the symmetric random walk. If $\mathcal{E}_p(H)<0$, they are exponentially less likely.

For a.e. $t \in [0,1]$, $H'(t)$ is the instantaneous direction contribution prescribed by the profile $H$, while $H(t)/t$ describes the accumulated increment distribution seen by the memory mechanism. We can therefore rewrite $\mathcal{E}_p(H)$ as
\begin{equation}\label{eq_energy_profile_with_KL_div}
\mathcal{E}_p(H)=\int_0^1 \Big\{D_{\textup{KL}} \left( \nu_{H'(t)} \| \nu_{\textup{unif}} \right) - D_{\textup{KL}} \left( \nu_{H'(t)} \| \nu_{Q_p(\cdot \mid H(t)/t)} \right) \Big\} \, \dd t,
\end{equation}
where $\nu_g$ denotes the probability measure on $\mathscr{S}_d$ associated with $g \in \Delta^{(d)}$ and $D_{\textup{KL}}(\cdot \| \cdot)$ is the Kullback--Leibler divergence. The first divergence of the integrand is the Sanov cost of sustaining the increment composition $H'(t)$ under the symmetric random walk. The second is the corresponding pathwise relative-entropy cost under the MERW. Thus, $\mathcal{E}_p(H)$ measures the time-integrated reduction, or increase when $\mathcal{E}_p(H)$ is negative, in the exponential cost of producing the instantaneous direction law prescribed by $H$ due to memory.
\end{remark}

\begin{remark}\label{remark_energy_integration_along_trajectory}
Proposition~\ref{proposition_expression_RL_rate_function}, together with the expressions for $\mathcal{E}_p$ given in \eqref{eq_def_energy_profile} and \eqref{eq_energy_profile_with_KL_div}, shows that the effect of the memory parameter enters through an integral functional of the entire cumulative profile, rather than through its endpoint alone. This is structurally close to the formula obtained in \cite[Equation~(5)]{Budhiraja_Waterbury} and \cite[Equation~(2.10)]{Zoubouloglou} for reinforced and self-interacting processes.
\end{remark}

The following corollary, which is a consequence of the preceding two lemmas, compares the probabilities of trajectories with similar profiles under the MERW and under the symmetric random walk.
\begin{corollary}[Local exponential tilting]\label{coro_tilting_ineq_for_localised_profile}
Fix $H \in \mathscr{H}^{(d)}$ and $\delta>0$. There exists $\varepsilon_0>0$ such that, for every $n$ sufficiently large and every $A_n \subseteq \Omega_n^{(0)}$ satisfying $H_n(A_n) \subseteq B_{\mathscr{H}}(H,\varepsilon_0)$,
\begin{equation}\label{eq_coro_tilting_ineq_for_localised_profile}
e^{n(\mathcal{E}_p(H)-2\delta)} \PP_{d,\textup{sym}}(A_n) \leq \PP_{d,p}(A_n) \leq e^{n(\mathcal{E}_p(H)+2\delta)} \PP_{d,\textup{sym}}(A_n).
\end{equation}
\end{corollary}

\begin{proof}
By Lemma~\ref{lemma_continuity_of_E_p}, choose $\varepsilon_0>0$ such that $|\mathcal{E}_p(H)-\mathcal{E}_p(\widetilde{H})|<\delta$ whenever $d_{\mathscr{H}}(H,\widetilde{H})<\varepsilon_0$. Since $H_n(A_n) \subseteq B_{\mathscr{H}}(H,\varepsilon_0)$, we deduce that 
\begin{align*}
|\mathcal{E}_p(H_n(w,\cdot))-\mathcal{E}_p(H)|<\delta, \qquad \forall \, w \in A_n.
\end{align*}
Thus, for $n$ sufficiently large so that $M_{d,p} \log (n)/n < \delta$, Lemma~\ref{lemma_exponential_tilting_wrt_SRW} gives, for every $w \in A_n$, 
\begin{align*}
\left \vert \frac{1}{n} \log \frac{\PP_{d,p}\left( S_{[1:n]}=w \right)}{\PP_{d,\textup{sym}}\big( S_{[1:n]}=w \big)} - \mathcal{E}_p(H) \right \vert & < 2 \delta.
\end{align*}
The result follows by exponentiating and summing over $w \in A_n$.
\end{proof}

\noindent We are now ready to give a proof of Proposition~\ref{proposition_expression_RL_rate_function}.
\begin{proof}[Proof of Proposition~\ref{proposition_expression_RL_rate_function}]
Fix $p \in (0,1)$, $(\mu,H) \in \mathcal{P}(\ZZ^d) \times \mathscr{H}^{(d)}$ and $\delta>0$. Let $\varepsilon_0>0$ be given by Corollary~\ref{coro_tilting_ineq_for_localised_profile}. For every $0<\varepsilon\leq\varepsilon_0$ and all sufficiently large $n$, that corollary gives
\begin{align}
    e^{n(\mathcal{E}_p(H)-2\delta)} \PP_{d,\textup{sym}}\big( \Omega_n(\mu,H,\varepsilon) \big) & \leq \PP_{d,p}\big( \Omega_n(\mu,H,\varepsilon) \big) \label{eq_proof_expression_RL_rate_function_1}\\
    & \leq e^{n(\mathcal{E}_p(H)+2\delta)} \PP_{d,\textup{sym}}\big( \Omega_n(\mu,H,\varepsilon) \big). \notag 
\end{align}
Since the sets $\Omega_n(\mu,H,\varepsilon)$ increase with $\varepsilon$, the infima in Definition~\ref{definition_RL} may be restricted to $0 < \varepsilon \leq \varepsilon_0$.

Taking logarithms, dividing by $n$ and then taking $\liminf_{n \to \infty}$ in the first inequality of \eqref{eq_proof_expression_RL_rate_function_1} gives for every $0 < \varepsilon \leq \varepsilon_0$,
\begin{align*}
&\liminf_{n\to\infty}\frac{1}{n} \log \PP_{d,p}\big(
(\ell_n,H_n)\in B(\mu,\varepsilon)\times B_{\mathscr{H}}(H,\varepsilon) \big)\\
&\qquad \qquad \qquad \qquad  \qquad\geq \liminf_{n\to\infty}\frac{1}{n} \log\PP_{d,\textup{sym}}\big((\ell_n,H_n)\in B(\mu,\varepsilon)\times B_{\mathscr{H}}(H,\varepsilon) \big) +\mathcal{E}_p(H)-2\delta\\
&\qquad \qquad \qquad \qquad \qquad\geq \underline{s}_{\textup{sym}}(\mu,H) +\mathcal{E}_p(H)-2\delta,
\end{align*}
where the last inequality follows directly from
Definition~\ref{definition_RL}. Taking the infimum over $\varepsilon \in (0,\varepsilon_0)$, 
\begin{equation}\label{eq_proof_expression_RL_rate_function_2}
\underline{s}_p(\mu,H)
\geq
\underline{s}_{\textup{sym}}(\mu,H)
+\mathcal{E}_p(H)-2\delta.
\end{equation}
Similarly, taking logarithms, dividing by $n$ and then taking $\limsup_{n \to \infty}$ in the second inequality of \eqref{eq_proof_expression_RL_rate_function_1} gives
\begin{equation}\label{eq_proof_expression_RL_rate_function_3}
\overline{s}_p(\mu,H)
\leq
\overline{s}_{\textup{sym}}(\mu,H)
+\mathcal{E}_p(H)+2\delta.
\end{equation}

By Definition~\ref{definition_RL}, $\underline{s}_p(\mu,H) \leq \overline{s}_p(\mu,H)$, while Lemma~\ref{lemma_intermediary_inequality_existence} gives
\begin{align*}
\underline{s}_{\textup{sym}}(\mu,H) = \overline{s}_{\textup{sym}}(\mu,H).
\end{align*}
Combining these facts with
\eqref{eq_proof_expression_RL_rate_function_2} and
\eqref{eq_proof_expression_RL_rate_function_3}, we obtain
\begin{align*}
\underline{s}_{\textup{sym}}(\mu,H)+\mathcal{E}_p(H)-2\delta \leq \underline{s}_p(\mu,H) \leq \overline{s}_p(\mu,H) \leq \underline{s}_{\textup{sym}}(\mu,H)+\mathcal{E}_p(H)+2\delta.
\end{align*}
Letting $\delta \to 0$ yields
\begin{equation}\label{proof_RL_for_p_f_notin_F_att_2}
\underline{s}_p(\mu,H)=\underline{s}_{\textup{sym}}(\mu,H)
+\mathcal{E}_p(H)=\overline{s}_p(\mu,H).
\end{equation}

Thus the Ruelle--Lanford function exists at every $(\mu,H)\in \mathcal{P}(\ZZ^d)\times \mathscr{H}^{(d)}$. Proposition~\ref{prop_RL_functions_imply_LDP} therefore yields a weak LDP under $\PP_{d,p}$ for $(\ell_n,H_n)_{n\in\NN}$, with rate
function
\begin{align*}
I_{\textup{A}}^{(p)}=-\underline{s}_p=-\overline{s}_p.
\end{align*}
Finally, \eqref{proof_RL_for_p_f_notin_F_att_2} gives
\begin{align*}
I_{\textup{A}}^{(p)}(\mu,H)=I_{\textup{A}}^{\textup{sym}}(\mu,H)-\mathcal{E}_p(H),
\end{align*}
which proves \eqref{eq_proposition_expression_RL_rate_function}.
\end{proof}

We conclude this section with a general result on $I_{\text{A}}^{(p)}$ which will be useful in Sections~\ref{subsection_comparison_to_SRW} and \ref{subsection_phase_transition_for_MERW}. As mentioned in Corollary~\ref{corollary_bound_on_last_point_of_subword}, if the empirical measures converge to a probability measure, the walk cannot move a macroscopic distance from the origin, so any limiting profile must balance opposite directions at every time. Consider the set $\mathscr{H}^{(d)}_{\text{bal}}$ defined by 
\begin{equation}\label{eq_def_H_bal}
\mathscr{H}^{(d)}_{\text{bal}} \eqdef \left\{ H \in \mathscr{H}^{(d)} \, : \, [H(t)]_{e_j} = [H(t)]_{-e_j} \, , \forall j \in \dbl 1,d \dbr \, , \, \forall t \in [0,1]\right\}.
\end{equation}
We call the elements of this set balanced profiles.

\begin{lemma}[Infinite cost of unbalanced profiles]\label{lemma_unbalanced_profiles_proba_nul}
For all $p \in (0,1)$, $\mu \in \mathcal{P}(\ZZ^d)$ and $H \notin \mathscr{H}^{(d)}_{\textup{bal}}$,
\begin{equation}\label{eq_lemma_unbalanced_profiles_proba_nul}
I_{\textup{A}}^{(p)}(\mu,H)=+\infty.
\end{equation}
\end{lemma}

\begin{proof}
Fix $p \in (0,1)$, $\mu \in \mathcal{P}(\ZZ^d)$ and $H \notin \mathscr{H}^{(d)}_{\textup{bal}}$. There exist $t^* \in (0,1]$ and $j^* \in \dbl 1,d \dbr$ such that 
\begin{align*}
[H(t^*)]_{e_{j^*}} \neq [H(t^*)]_{-e_{j^*}}.
\end{align*}
Choose $\varepsilon>0$ such that 
\begin{equation}\label{proof_lemma_unbalanced_profiles_proba_nul_1}
8 \varepsilon < \big \vert [H(t^*)]_{e_{j^*}} - [H(t^*)]_{-e_{j^*}} \big \vert.
\end{equation}
Let $R_\mu(\varepsilon)$ be given by Lemma~\ref{lemma_bound_points_outside_compact}, $n > \max \left\{2,R_{\mu}(\varepsilon)/\varepsilon \right\}$ and $w \in \Omega_n(\mu,\varepsilon)$. Corollary~\ref{corollary_bound_on_last_point_of_subword} yields
\begin{equation}\label{proof_lemma_unbalanced_profiles_proba_nul_2}
\|w_k\|_{\infty} < 3 \varepsilon n, \qquad  k \in \dbl 1,n \dbr.
\end{equation}

At each mesh point $t=(k-1)/(n-1)$, since $w_1=0$, the definition of the step counters gives
\begin{align*}
[H_n(w,t)]_{e_{j^*}}-[H_n(w,t)]_{-e_{j^*}} & = \frac{C_k^{(e_{j^*})}(w)}{n-1}-\frac{C_k^{(-e_{j^*})}(w)}{n-1} \\
& =\frac{\pi_{j^*}(w_k)}{n-1},
\end{align*}
where $\pi_{j^*}$ denotes projection onto the $j^*$-th coordinate. Using the fact that $H_n(w,\cdot)$ is affine between the mesh points and $n \geq 2$, together with \eqref{proof_lemma_unbalanced_profiles_proba_nul_2},
\begin{equation}\label{proof_lemma_unbalanced_profiles_proba_nul_3}
\sup_{t \in [0,1]} \big| [H_n(w,t)]_{e_{j^*}}-[H_n(w,t)]_{-e_{j^*}} \big| < \frac{3 \varepsilon n}{n-1} \leq 6 \varepsilon.
\end{equation}
Combining \eqref{proof_lemma_unbalanced_profiles_proba_nul_3} and \eqref{proof_lemma_unbalanced_profiles_proba_nul_1} gives
\begin{align*}
2d_{\mathscr{H}}\big( H_n(w,\cdot),H \big) > 2 \varepsilon.
\end{align*}
Thus, for all sufficiently large $n$ and all $w \in \Omega_n(\mu,\varepsilon)$, $H_n(w,\cdot) \notin B_{\mathscr{H}}(H,\varepsilon)$, and therefore
\begin{align*}
\PP_{d,p} \big( (\ell_n, H_n) \in B(\mu,\varepsilon) \times B_{\mathscr{H}}(H,\varepsilon) \big)=0.
\end{align*}
Taking the logarithm, dividing by $n$, taking the $\limsup_{n \to \infty}$ and taking the infimum over $\varepsilon>0$ yields
\begin{align*}
\overline{s}_p(\mu,H) & = \inf_{\varepsilon>0 } \limsup_{n \to \infty} \frac{1}{n} \log \PP_{d,p} \big( (\ell_n,H_n) \in B(\mu,\varepsilon) \times B_{\mathscr{H}}(H,\varepsilon) \big) \\
& = - \infty.
\end{align*}
Equation~\eqref{eq_lemma_unbalanced_profiles_proba_nul} follows from the fact that $I_{\textup{A}}^{(p)}=-\overline{s}_p$.
\end{proof}

\subsection{Contraction to empirical measures}\label{subsection_existence_of_a_weak_LDP_measures}

The aim of this subsection is to prove Theorem~\ref{theorem_existence_LDP} which states that the sequence of empirical measures $(\ell_n)_{n \in \NN}$ satisfies a weak LDP under $\PP_{d,p}$.

The following proposition is slightly stronger than Theorem~\ref{theorem_existence_LDP}, since it establishes the existence of a weak LDP for the sequence of empirical measures $(\ell_n)_{n \in \NN}$ and provides two alternative expressions for the rate function.

\begin{proposition}[Weak LDP for empirical measures]\label{proposition_existence_weak_LDP_with_different_rate_fct}
For every $d\in\NN$ and $p\in(0,1)$, under $\PP_{d,p}$, the sequence $(\ell_n)_{n\in\NN}$ satisfies a weak LDP with rate function $I_{d,p}$ given, for every $\mu\in\mathcal{P}(\ZZ^d)$, by
\begin{align}
I_{d,p}(\mu) & = \inf \limits_{H \in \mathscr{H}^{(d)}} I_{\textup{A}}^{(p)}(\mu,H) \label{eq_weak_LDP_rate_function_expression_1}\\
& = \inf_{H \in \mathscr{H}^{(d)}} \Big\{ I_{\textup{A}}^{\textup{sym}}(\mu,H) - \mathcal{E}_p(H)\Big\} \label{eq_weak_LDP_rate_function_expression_2}\\
& = -\inf_{\varepsilon>0} \liminf_{n \to \infty} \frac{1}{n} \log \PP_{d,p}\big( \ell_n \in B(\mu,\varepsilon) \big)=-\inf_{\varepsilon>0} \limsup_{n \to \infty} \frac{1}{n} \log \PP_{d,p}\big( \ell_n \in B(\mu,\varepsilon) \big). \label{eq_weak_LDP_rate_function_expression_3}
\end{align}
\end{proposition}

\begin{proof}  
Let $\pi_1:\mathcal{P}(\ZZ^d) \times \mathscr{H}^{(d)} \longrightarrow \mathcal{P}(\ZZ^d)$ denote the projection onto the first coordinate. This map is continuous. Moreover, by the third item of Lemma~\ref{lemma_properties_profile_space}, $\mathscr{H}^{(d)}$ is compact. Thus, for every compact set $\mathcal{K} \subseteq \mathcal{P}(\ZZ^d)$,
\begin{align*}
\pi^{-1}_1(\mathcal{K})=\mathcal{K} \times \mathscr{H}^{(d)}
\end{align*}
is compact and so $\pi_1$ is proper. Applying the weak contraction principle \cite[Lemma~4.1]{Daures} to the weak LDP for $(\ell_n,H_n)_{n\in\NN}$ yields a weak LDP for $(\ell_n)_{n\in\NN}$, with rate function
\begin{align*}
I_{d,p}(\mu) = \inf \limits_{H \in \mathscr{H}^{(d)}} I_{\textup{A}}^{(p)}(\mu,H),
\end{align*}
which is exactly \eqref{eq_weak_LDP_rate_function_expression_1}. The alternative expression \eqref{eq_weak_LDP_rate_function_expression_2} then follows immediately from Proposition~\ref{proposition_expression_RL_rate_function}.

Since $\ZZ^d$ is Polish, the space $\mathcal{P}(\ZZ^d)$, endowed with the weak topology, is itself Polish \cite[Theorem~3.2.2]{Bogachev}. Hence, by the local characterisation of the rate function for a weak LDP, see the remark following \cite[Theorem~4.1.18]{DZ_LDP}, for every $\mu\in\mathcal{P}(\ZZ^d)$,
\begin{align*}
I_{d,p}(\mu)&= -\inf_{\varepsilon>0} \liminf_{n \to \infty} \frac{1}{n} \log \PP_{d,p}\big( \ell_n \in B(\mu,\varepsilon) \big) \\
& = -\inf_{\varepsilon>0} \limsup_{n \to \infty} \frac{1}{n} \log \PP_{d,p}\big( \ell_n \in B(\mu,\varepsilon) \big). 
\end{align*}
This yields the representation \eqref{eq_weak_LDP_rate_function_expression_3}.
\end{proof}

\noindent With Proposition~\ref{proposition_existence_weak_LDP_with_different_rate_fct}, Theorem~\ref{theorem_existence_LDP} is proved.

\begin{remark}\label{remark_weak_LDP_contraction}
Since we only have a weak LDP, the usual contraction principle cannot be applied blindly. Compactness of $\mathscr{H}^{(d)}$ is essential for the argument to work.
\end{remark}

\begin{remark}\label{remark_use_of_different_expressions}
In Section~\ref{subsection_comparison_to_SRW}, Equation~\eqref{eq_weak_LDP_rate_function_expression_2} allows us to compare $I_{d,p}$ to $I_{d,\text{DV}}$, the rate function of the empirical measures of the symmetric random walk.
\end{remark}

\noindent We now show that this weak LDP cannot be strengthened to a full LDP. 

\begin{proof}[Proof of Proposition~\ref{prop_no_full_LDP}]
This proof is a slight modification of an argument given in \cite[pp.~922--923]{Baxter_Jain_Varadhan}. Consider the set
\begin{align*}
\mathcal{C}\eqdef \left\{ \frac{1}{2n+1} \sum_{k=0}^{2n} \delta_{(k,0,\dots,0)} \right\}_{n \in \NN_0} \subseteq \mathcal{P}(\ZZ^d).
\end{align*}
The set $\mathcal{C}$ is closed for the topology of weak convergence. Indeed, for all $m \in \NN_0$, let
\begin{align*}
\mu_m=\frac{1}{2m+1}\sum_{k=0}^{2m}\delta_{(k,0,\ldots,0)}.
\end{align*}
Let $(\mu_{m_n})_{n\in\NN}$ be a sequence in $\mathcal{C}$ converging weakly to some $\mu \in \mathcal{P}(\ZZ^d)$. Suppose that $(m_n)_{n\in\NN}$ is unbounded. Passing to a subsequence, we may assume that $m_n\to+\infty$. For every $x \in \ZZ^d$,
\begin{align*}
0\leq\mu_{m_n}(\{x\})\leq\frac{1}{2m_n+1},
\end{align*}
and hence $\mu_{m_n}(\{x\}) \to 0$. Since $\ZZ^d$ is discrete, the indicator of $\{x\}$ is bounded and continuous. Weak convergence therefore implies
\begin{align*}
\mu(\{x\}) = \lim_{n\to\infty}\mu_{m_n}(\{x\}) = 0, \qquad x\in\ZZ^d.
\end{align*}
This is impossible, since $\ZZ^d$ is countable and $\mu(\ZZ^d)=1$. Thus, $(m_n)_{n\in\NN}$ is bounded. The sequence $(m_n)_{n\in\NN}$ consequently takes values in a finite set, so it admits a constant subsequence, say $m_{n_k}=m$. This subsequence converges both to $\mu$ and to $\mu_m$, and hence $\mu=\mu_m\in\mathcal{C}$. Thus, $\mathcal{C}$ is sequentially closed. Since $\mathcal{P}(\ZZ^d)$ is metrisable for the topology of weak convergence, $\mathcal{C}$ is closed.

Suppose for the sake of contradiction that the sequence $(\ell_n)_{n \in \NN}$ satisfies a full LDP in $\mathcal{P}(\ZZ^d)$ with rate function $\widetilde{I}$. Since $\mathcal{P}(\ZZ^d)$ is a regular topological space, the local characterisation of the rate function in \cite[Theorem~4.1.18]{DZ_LDP} and the definition of $I_{d,p}$, given in \eqref{expression_rate_function}, yield $\widetilde{I}=I_{d,p}$. 

On the one hand, let $n \in \NN$ and consider the trajectory
\begin{align*}
\overline{w}^{(n)} \eqdef \big(0,e_1, \dots , 2n e_1\big) \in \Omega_{2n+1}^{(0)},
\end{align*}
where we recall that $e_1=(1,0,\dots,0)$. We have 
\begin{align*}
\ell_{2n+1}(\overline{w}^{(n)})=\frac{1}{2n+1} \sum_{k=0}^{2n} \delta_{ke_1} \in \mathcal{C}.
\end{align*}
Thus, by inclusion of events and using Corollary~\ref{corollary_product_of_transitions}, we obtain 
\begin{align*}
\PP_{d,p}\big( \ell_{2n+1} \in \mathcal{C} \big)& \geq \PP_{d,p}\left(S_{[1:2n+1]}=\overline{w}^{(n)} \right) \\
& = \frac{1}{2d}\prod_{j=2}^{2n}Q_p\big(e_1 \mid \rho_j(\overline{w}^{(n)}) \big).
\end{align*}
Then, with the lower bound given in Lemma~\ref{lemma_lower_bound_transition}, we deduce that
\begin{align*}
\PP_{d,p}\big( \ell_{2n+1} \in \mathcal{C} \big) \geq \frac{1}{2d}\, p_*^{2n-1}.
\end{align*}
We deduce that
\begin{align}
\log p_* & \leq \limsup_{n \to \infty} \frac{1}{2n+1} \log \PP_{d,p}\big( \ell_{2n+1} \in \mathcal{C} \big) \notag \\
& \leq \limsup_{n \to \infty} \frac{1}{n} \log \PP_{d,p}\big( \ell_{n} \in \mathcal{C} \big). \label{eq_inequality_for_no_full_LDP_1}
\end{align}

On the other hand, since we have assumed that $(\ell_n)_{n \in \NN}$ satisfies a full LDP with rate function $I_{d,p}$, we obtain from the LDP upper bound with the closed set $\mathcal{C}$ that
\begin{align*}
\limsup_{n \to \infty} \frac{1}{n} \log \PP_{d,p}\big( \ell_n \in \mathcal{C}\big) \leq -\inf_{\mu \in \mathcal{C}} I_{d,p}(\mu).
\end{align*}
Consider $\mu \in \mathcal{C}$. There exists $m \in \NN_0$ such that
\begin{align*}
\mu=\mu_m \eqdef \frac{1}{2m+1} \sum_{k=0}^{2m} \delta_{ke_1}.
\end{align*} 
Consider $0<\varepsilon< \frac{1}{8m+4}$ and let $w \in \Omega_n^{(0)}$ such that $\ell_n(w) \in B(\mu_m,\varepsilon)$. Set
\begin{align*}
E^{(d)} \eqdef \left\{ (x_1,\dots,x_d) \in \ZZ^d \, : \, \sum_{i=1}^d x_i \in 2\ZZ \right\}.
\end{align*}
Using the definition of $\mu_m$ and $E^{(d)}$, we have
\begin{align*}
\mu_m\big( E^{(d)} \big) = \frac{1}{2m+1} \sum_{k=0}^{2m} \mathbbm{1}_{\{k \in 2\ZZ \}}=\frac{m+1}{2m+1}.
\end{align*}
Recall that for $A \subseteq \ZZ^d$, we defined $N_A(w)$ in \eqref{eq_def_occupation_time_N_A} as the occupation time of the trajectory $w$ within $A$. Using the fact that $\ell_{n}(w) \in B(\mu_m,\varepsilon)$ and Lemma~\ref{lemma_occupation_meas_with_empirical_meas}, we deduce that
\begin{align*}
\frac{N_{E^{(d)}}(w)}{n} & \geq \mu_m\big(E^{(d)} \big) - \varepsilon > \frac{m+1}{2m+1}- \frac{1}{8m+4}=\frac{1}{2}+\frac{1}{8m+4}.
\end{align*}
Therefore, 
\begin{equation}\label{eq_intermediaire_for_no_full_LDP}
\PP_{d,p} \big( \ell_n \in B(\mu_m,\varepsilon)\big) \leq \PP_{d,p} \left( \frac{N_{E^{(d)}}\big( S_{[1:n]} \big)}{n} > \frac{1}{2} + \frac{1}{8m+4} \right).
\end{equation}
However, since $(S_n)_{n \in \NN}$ is a nearest-neighbour random walk started at 0, we deduce that
\begin{align*}
S_{n} \in E^{(d)}, \quad  \forall n \in 2\NN_0 +1 \qquad \text{and} \qquad S_n \in \ZZ^d \setminus E^{(d)}, \quad \forall n \in 2\NN.
\end{align*}
Thus, we have for all $w \in \Omega_n^{(0)}$,
\begin{align*}
\frac{N_{E^{(d)}}(w)}{n} = 
\begin{cases}
\frac{1}{2}, & \text{if } n \in 2\NN, \\
\frac{1}{2}+\frac{1}{2n}, & \text{if } n \in 1+2\NN. 
\end{cases}
\end{align*}
Therefore, for all $n > 4m+2$ and all $w \in \Omega_n^{(0)}$, we have $\frac{N_{E^{(d)}}(w)}{n}<\frac{1}{2}+\frac{1}{8m+4}$. In other words,
\begin{equation}\label{eq_intermediaire_for_no_full_LDP_2}
\PP_{d,p} \left(\frac{N_{E^{(d)}}\big( S_{[1:n]} \big)}{n} < \frac{1}{2} + \frac{1}{8m+4} \right)=1.
\end{equation}
Combining \eqref{eq_intermediaire_for_no_full_LDP} and \eqref{eq_intermediaire_for_no_full_LDP_2}, we deduce that for all $0<\varepsilon<\frac{1}{8m+4}$ and all $n>4m+2$,
\begin{align*}
\PP_{d,p} \big(\ell_n \in B(\mu_m,\varepsilon) \big)=0.
\end{align*}
In other words, since $\mu_m$ does not assign half of its mass to $E^{(d)}$, it is asymptotically impossible for the empirical measure of a MERW to be $\varepsilon$-close to $\mu_m$ when $\varepsilon$ is sufficiently small. Since this holds for all $m \in \NN_0$, we deduce that for all $\mu \in \mathcal{C}$,
\begin{align*}
\inf_{\varepsilon>0} \limsup_{n \to \infty} \frac{1}{n} \log \PP_{d,p} \big(\ell_n \in B(\mu,\varepsilon) \big) = - \infty.
\end{align*}
Using the expression for $I_{d,p}$ given in \eqref{expression_rate_function} and the above equality, we deduce that for all $\mu \in \mathcal{C}$, we have $I_{d,p}(\mu)=+\infty$. Therefore,
\begin{align*}
\limsup_{n \to \infty} \frac{1}{n} \log \PP_{d,p} \big(\ell_n \in \mathcal{C} \big) & \leq - \inf_{\mu \in \mathcal{C}} I_{d,p}(\mu) \\
& = -\infty,
\end{align*}
which contradicts \eqref{eq_inequality_for_no_full_LDP_1}. Thus, the sequence $(\ell_n)_{n \in \NN}$ cannot satisfy a full LDP in $\mathcal{P}(\ZZ^d)$ and Proposition~\ref{prop_no_full_LDP} is proved.
\end{proof}

\subsection{LDP for profile and path observables}\label{subsection_profile_and_path_observables}

In this subsection, we show that the exponential tilting introduced in Section~\ref{subsection_comparaison_RL_functions} provides a direct route to large deviation principles for other observables of the MERW.

\begin{corollary}[LDPs for the profile and position]
\label{corollary_LDP_profile_position}
Fix $d \in \NN$ and $p \in (0,1)$. Under $\PP_{d,p}$, the sequence $(H_n)_{n\in\NN}$ satisfies a full LDP on $\mathscr{H}^{(d)}$, endowed with the topology of uniform convergence, with good rate function $J_{d,p} : \mathscr{H}^{(d)} \longrightarrow [0,+\infty]$ given by
\begin{equation}\label{eq_corollary_LDP_profile_position_1}
J_{d,p}(H) = \int_0^1 D_{\textup{KL}}  \left( \nu_{H'(t)} \big\| \nu_{Q_p( \cdot \,| \,H(t)/t)} \right)\textup{d} t, \qquad H \in \mathscr{H}^{(d)}.
\end{equation}
Consequently, $(S_n/n)_{n\in\NN}$ satisfies a full LDP on $\RR^d$ with good rate function $K_{d,p} : \RR^d \longrightarrow [0,+\infty]$ given by
\begin{equation}\label{eq_corollary_LDP_profile_position_2}
K_{d,p}(x) = \inf \left\{ J_{d,p}(H) \, : \, H \in \mathscr{H}^{(d)} , \, \sum_{s \in \mathscr{S}_d}s \, [H(1)]_s =x \right\}.
\end{equation}
\end{corollary}

\begin{proof}
Recall that the sequence $(X_j)_{j \in \NN}$ consists of the steps of the MERW. Under $\PP_{d,\mathrm{sym}}$, the variables $\delta_{X_j}\in\RR^{\mathscr{S}_d}$ are i.i.d. and uniformly distributed on the canonical basis. For $n \geq 2$, consider
\begin{align*}
Z_n(t)=\frac{1}{n-1}\sum_{j=1}^{\lfloor (n-1)t \rfloor} \delta_{X_j}, \qquad t \in [0,1].
\end{align*} 
Consider also the associated log-Laplace transform which is finite everywhere and equals
\begin{align*}
\Lambda(\lambda)=\log \EE_{d,\text{sym}}\left[ e^{\langle \lambda,\delta_{X_1} \rangle} \right]=\log \left( \frac{1}{2d} \sum_{s \in \mathscr{S}_d} e^{\lambda_s} \right), \qquad \lambda \in \RR^{\mathscr{S}_d}.
\end{align*}
By Mogulskii's theorem
\cite[Theorem~5.1.2 and Remark~(b), p.~176]{DZ_LDP}, $(Z_n)_{n\in\NN}$ satisfies a full LDP in $\mathcal{D}_0([0,1],\RR^{\mathscr{S}_d})$, the space of càdlàg paths such that $f(0)=0$, equipped with the uniform topology, with good rate function 
\begin{align*}
\widehat{J}_{d,\text{sym}}(\phi)=
\begin{cases}
\int_0^1 \Lambda^*(\phi'(t)) \dd t, & \text{if } \phi \text{ is absolutely continuous and } \phi(0)=0, \\
+\infty, & \text{otherwise},
\end{cases}
\end{align*}
where $\Lambda^*$ denotes the Legendre transform of $\Lambda$.

By definition $H_n$ is the polygonal interpolation of $Z_n$. 
Thus, similarly to \cite[Lemma~5.1.4]{DZ_LDP}, we have $\| H_n - Z_n \|_{\infty} \leq 1/(n-1)$ and so the sequences $(H_n)_{n \in \NN}$ and $(Z_n)_{n \in \NN}$ are exponentially equivalent. By Lemmas~\ref{lemma_profiles_in_profile_space} and \ref{lemma_properties_profile_space}, $H_n$ takes values in $\mathscr{H}^{(d)}$, which is compact and therefore closed in this path space. Thus,  by \cite[Lemma~4.1.5]{DZ_LDP}, the process $(H_n)_{n \in \NN}$ satisfies a full LDP in $\mathscr{H}^{(d)}$.
 
Moreover, by \cite[Corollary~2.2.35 and Exercise~2.2.36]{DZ_LDP}, applied to the uniform distribution on $\mathscr{S}_d$, the Legendre transform of $\Lambda$ satisfies
\begin{align*}
\Lambda^*(q) = \sum_{s \in \mathscr{S}_d} q_s \log(2d q_s) = D_{\textup{KL}}\big(\nu_q \| \nu_{\textup{unif}}\big), \qquad q \in \Delta^{(d)},
\end{align*}
with the convention $0\log 0=0$. Every $H \in \mathscr{H}^{(d)}$ is absolutely continuous, starts at zero, and satisfies $H'(t) \in \Delta^{(d)}$ by Lemma~\ref{lemma_properties_profile_space}. Thus, the symmetric profile rate function is 
\begin{equation}
J_{d,\text{sym}}(H)=\int_0^1 D_{\textup{KL}}\big(\nu_{H'(t)} \| \nu_{\text{unif}} \big) \dd t.
\end{equation}

Corollary~\ref{coro_tilting_ineq_for_localised_profile} then shows that the local upper and lower large-deviation rates under $\PP_{d,p}$ coincide and are given by
\begin{align*}
J_{d,p}(H)=J_{d,\text{sym}}(H)-\mathcal{E}_p(H).
\end{align*}
Using \eqref{eq_energy_profile_with_KL_div}, we obtain precisely \eqref{eq_corollary_LDP_profile_position_1}. The local LDP therefore yields a weak LDP, and this LDP is full since $\mathscr{H}^{(d)}$ is compact by Lemma~\ref{lemma_properties_profile_space}.

Finally, setting
\begin{align*}
\Phi(H) \eqdef
\sum_{s\in\mathscr{S}_d}s \, [H(1)]_s, \qquad H \in \mathscr{H}^{(d)},
\end{align*}
we have $\Phi(H_n)=S_n/(n-1)$. Evaluation at 1 is continuous for the uniform topology, so $\Phi$ is continuous. The contraction principle therefore gives a full LDP for $S_n/(n-1)$ with rate function $K_{d,p}$. Moreover, $(S_n/n)$ and $(S_n/(n-1))$ are exponentially equivalent, which concludes the proof.
\end{proof}

\begin{remark}\label{remark_same_as_franchini}
As mentioned in the introduction, the MERW admits an equivalent formulation as a Pólya urn. In dimension $d=1$, set $h(t)=[H(t)]_{+ e_1}$, so that $[H(t)]_{-e_1}=t-h(t)$. Introducing $\pi_p(x)=1-p+(2p-1)x$, the rate function becomes
\begin{align*}
J_{1,p}(H)=\int_0^1 \left[ h'(t) \log \left( \frac{h'(t)}{\pi_p(h(t)/t)} \right) + (1-h'(t)) \log \left( \frac{1-h'(t)}{1-\pi_p(h(t)/t)} \right) \right] \dd t.
\end{align*} 
This coincides with the rate function obtained in \cite[Equation~(2.7)]{Franchini_urns} for the binary urn associated with the one-dimensional ERW. Hence, in dimension $d=1$, Corollary~\ref{corollary_LDP_profile_position} recovers, in the present notation, the corresponding LDP established in \cite[Theorem~1]{Franchini_urns}.
\end{remark}

The following lemma shows that while the rate function for $(\ell_n)_{n \in \NN}$ is convex in some regimes, the rate functions for the sequence of profiles and for the endpoints $S_n/n$ are always convex.
\begin{lemma}\label{lemma_convexity_rate_profile_rate_function}
    For all $d \in \NN$ and $p \in (0,1)$, the rate functions $J_{d,p}$ and $K_{d,p}$ are convex.
\end{lemma}

\begin{proof}
    Consider $H,\widehat{H} \in \mathscr{H}^{(d)}$ and $\lambda \in [0,1]$. Set $H_\lambda=\lambda H+(1-\lambda)\widehat{H}$. The defining conditions of $\mathscr{H}^{(d)}$ are preserved by convex combinations, so $H_\lambda \in \mathscr{H}^{(d)}$. By \eqref{eq_def_Q_p_r_s}, the map $r \mapsto Q_p(\cdot \mid r)$ is affine on $\Delta^{(d)}$. It follows that
    \begin{align*}
        \nu_{ Q_p\left( \cdot \, | \, H_{\lambda}(t)/t  \right)} = \lambda \nu_{ Q_p\left( \cdot \,|\, H(t)/t \right)} + (1-\lambda) \nu_{ Q_p\left( \cdot \,|\, \widehat{H}(t)/t \right)}.
    \end{align*}
    Moreover, the absolute continuity of $H$ and $\widehat H$ implies that, for almost every $t \in (0,1)$,
    \begin{align*}
        \nu_{H_{\lambda}'(t)}=\lambda \nu_{H'(t)} + (1-\lambda) \nu_{\widehat{H}'(t)}.
    \end{align*}
    The joint convexity of relative entropy, see \cite[Theorem~2.7.2]{Cover_Thomas}, therefore gives, for almost every $t \in (0,1)$,
    \begin{align*}
        & D_{\textup{KL}}\left( \nu_{H_{\lambda}'(t)} \, \| \, \nu_{Q_p(\cdot \, | \,H_{\lambda}(t)/t )} \right) \\
        & \qquad \qquad \qquad \qquad \qquad \quad \leq \lambda D_{\textup{KL}}\left( \nu_{ H'(t)} \, \| \, \nu_{Q_p(\cdot \, | \, H(t)/t)} \right) + (1-\lambda) D_{\textup{KL}}\left( \nu_{\widehat{H}'(t)} \, \| \, \nu_{Q_p(\cdot \, | \, \widehat{H}(t)/t)} \right).
    \end{align*}
    Integrating the inequality and using \eqref{eq_corollary_LDP_profile_position_1}, we deduce that $J_{d,p}$ is convex. 

    The convexity of $K_{d,p}$ follows from that of $J_{d,p}$, since $K_{d,p}$ is obtained by contraction through the linear map $H \mapsto \sum_{s \in \mathscr{S}_d} s[H(1)]_s$, see \cite[Section~3.2.5]{Boyd_Vandenberghe}.
\end{proof}

\medskip

Fix a norm $\| \cdot \|$ on $\RR^d$. The following result follows directly from Corollary~\ref{corollary_LDP_profile_position}. Consider the rescaled maximal displacement
\begin{align*}
M_n \eqdef \frac{1}{n}\max_{1\leq k\leq n}\|S_k\|.
\end{align*}
\begin{corollary}[LDP for the maximal displacement]\label{corollary_maximal_displacement}
Fix $d \in \NN$ and $p \in (0,1)$. Under $\PP_{d,p}$, the sequence $(M_n)_{n \in \NN}$ satisfies a full LDP on $\RR$ with rate function $L_{d,p}:\RR \longrightarrow [0,+\infty]$ given by
\begin{equation}\label{eq_corollary_maximal_displacement}
L_{d,p}(x)=\inf \left\{ J_{d,p}(H) \, : \, H \in \mathscr{H}^{(d)}\, ,\,  \sup_{t \in [0,1]} \Big\| \sum_{s \in \mathscr{S}_d} s \big[H(t)\big]_s \Big\|=x  \right\}
\end{equation}
\end{corollary}

\begin{proof}
Define $\Psi : \mathscr{H}^{(d)} \longrightarrow \RR$ by 
\begin{align*}
\Psi(H) \eqdef \sup_{t \in [0,1]}\Big\| \sum_{s \in \mathscr{S}_d} s \big[H(t)\big]_s \Big\|, \qquad H \in \mathscr{H}^{(d)}.
\end{align*}
Set $C=\max_{s \in \mathscr{S}_d}\|s\|$. For $H,\widehat{H} \in \mathscr{H}^{(d)}$, the triangle inequality gives
\begin{align*}
\big| \Psi(H)-\Psi(\widehat{H}) \big| & \leq \sup_{t \in [0,1]} \Big\| \sum_{s \in \mathscr{S}_d} s \big( [H(t)]_s - [\widehat{H}(t)]_s \big) \Big\| \\
& \leq 2dC d_{\mathscr{H}}(H,\widehat{H}).
\end{align*}
Thus, $\Psi$ is continuous for the topology of uniform convergence. In addition, since $H_n$ is defined on the mesh points and then extended affinely, any maximum must be reached at one of these mesh points. We deduce that
\begin{align*}
\Psi(H_n) & =\max_{1 \leq j \leq n} \left\| \sum_{s \in \mathscr{S}_d} s \left[H_n\left( \frac{j-1}{n-1} \right)\right]_s \right\| \\
& = \frac{n}{n-1} M_n.
\end{align*}
By Corollary~\ref{corollary_LDP_profile_position} and the contraction principle, the sequence $(\Psi(H_n))_{n \in \NN}$ satisfies an LDP with rate function $L_{d,p}$ given in \eqref{eq_corollary_maximal_displacement}. Finally, since every increment has norm at most $C$, we have $\|S_j\| \leq (j-1)C$. Hence,
\begin{align*}
|\Psi(H_n)-M_n|=\frac{M_n}{n-1} \leq \frac{C}{n}, \qquad n \geq 2.
\end{align*}
Thus, $(\Psi(H_n))_{n \in \NN}$ and $(M_n)_{n \in \NN}$ are exponentially equivalent, and $(M_n)_{n \in \NN}$ satisfies the same full LDP with rate function $L_{d,p}$.
\end{proof}

\begin{remark}\label{remark_convexity_question_L_d_p}
    It remains unclear whether the rate function $L_{d,p}$ is convex, since it is obtained by contracting the convex rate function $J_{d,p}$ through a nonlinear map, which does not necessarily preserve convexity.
\end{remark}

\begin{remark}[Possible generalisations]\label{remark_generalisation}
More generally, let $(Y_n)_{n \in \NN}$ be a sequence of $\mathcal{Y}$-valued observables of the MERW trajectory and suppose that, under $\PP_{d,\textup{sym}}$, one can establish an LDP for the pair $(Y_n,H_n)$ in $\mathcal{Y} \times \mathscr{H}^{(d)}$. Since the exponential change of measure from $\PP_{d,\textup{sym}}$ to $\PP_{d,p}$ depends asymptotically only on the profile $H_n$, Lemma~\ref{lemma_exponential_tilting_wrt_SRW} suggests that the same LDP can then be transferred to $\PP_{d,p}$ by replacing the symmetric rate function $I^{\mathrm{sym}}$ with
\begin{align*}
I^{(p)}(y,H)=I^{\mathrm{sym}}(y,H)-\mathcal{E}_p(H).
\end{align*}
By contraction, one would then obtain a weak LDP for the sequence $(Y_n)_{n \in \NN}$ with rate function
\begin{align*}
J^{(p)}(y)=\inf_{H \in \mathscr{H}^{(d)}} \big\{ I^{\mathrm{sym}}(y,H)-\mathcal{E}_p(H) \big\}, \qquad y \in \mathcal{Y}.
\end{align*}
\end{remark}

\section{Properties of the rate function}\label{section_properties_rate_function}

While Proposition~\ref{proposition_existence_weak_LDP_with_different_rate_fct} does not give a particularly explicit expression for the rate function $I_{d,p}$, we establish some of its properties in this section. First, in Section~\ref{subsection_invariance_of_rate_function}, we prove that the rate function is invariant under translations, signed permutations, and changes in the distribution of the initial step. Then, in Section~\ref{subsection_comparison_to_SRW}, we compare this rate function with the Donsker--Varadhan rate function for the simple random walk, first in dimension $d=1$ and then in dimensions $d\geq2$, highlighting the fundamentally different behaviour that arises between the one-dimensional and higher-dimensional cases. Finally, in Section~\ref{subsection_phase_transition_for_MERW}, we study the convexity of the rate function.

\subsection{Symmetries of the rate function}\label{subsection_invariance_of_rate_function} 

In this section only, we consider the MERW under an initial distribution that may be different from $\nu_{\text{unif}}$, the uniform measure on $\mathscr{S}_d$. For $\nu \in \mathcal{P}(\mathscr{S}_d)$, we consider a probability space $(\Omega,\mathcal{F},\PP_{d,p,\nu})$ supporting a MERW $(S_n)_{n \in \NN}$ with memory parameter $p \in (0,1)$ and initial step distribution $\nu$.

In Section~\ref{subsection_existence_of_a_weak_LDP_measures}, we showed that under $\PP_{d,p,\nu_{\text{unif}}}$, the sequence of empirical measures $(\ell_n)_{n \in \NN}$ satisfies a weak LDP with rate function $I_{d,p}$. In this section, we extend the results of Proposition~\ref{proposition_existence_weak_LDP_with_different_rate_fct} to an arbitrary initial step distribution $\nu \in \mathcal{P}(\mathscr{S}_d)$ and show that the choice of $\nu$ is irrelevant. Furthermore, we show that the weak LDP rate function $I_{d,p}$ is invariant under translations and signed permutations of coordinate axes. 

We start with a technical lemma. Choose a nearest-neighbour path $b \in \Omega_{r+1}^{(0)}$ from 0 to $z \in \ZZ^d$,
\begin{align*}
b_1=0,b_2, \dots , b_{r+1}=z,
\end{align*}
where $r \geq \|z\|_1$. For $w \in \Omega_n^{(0)}$, we define
\begin{equation}\label{eq_def_Theta_n_z}
\Theta_{b}^{(n)}(w) \eqdef b_1 b_2 \cdots b_r z (z+w_2) \cdots (z+w_n) \in \Omega_{n+r}^{(0)},
\end{equation}
which is the trajectory $w$ shifted by the prefix $b$. The following lemma shows that applying $\Theta_b^{(n)}$ to a trajectory does not affect, at an exponential scale, the probability of realising this trajectory under a different initial step distribution.

\begin{lemma}[Subexponential cost of a prefix]\label{lemma_prefix_does_not_change}
Fix $\nu$, $\widetilde{\nu} \in \mathcal{P}(\mathscr{S}_d)$, $z \in \ZZ^d$ and a nearest-neighbour path $b \in \Omega_{r+1}^{(0)}$ from 0 to $z$ such that $\nu(\Delta_1(b))>0$. There exists a sequence $(g_{b}(n))_{n \geq 2}$ such that for all $n \geq 2$, $ g_{b}(n)>0$, and
\begin{equation*}
\lim_{n \to \infty} \frac{1}{n} \log g_{b}(n)=0,
\end{equation*}
and such that for every $w \in \Omega_n^{(0)}$,
\begin{equation}\label{eq_lemma_prefix_does_not_change}
\PP_{d,p,\nu}\left( S_{[1:n+r]}=\Theta_{b}^{(n)}(w) \right) \geq g_{b}(n) \, \PP_{d,p,\widetilde{\nu}}\big( S_{[1:n]}=w \big).
\end{equation} 
\end{lemma}

\begin{proof}
Fix a nearest-neighbour path $b \in \Omega_{r+1}^{(0)}$ from 0 to $z \in \ZZ^d$ such that $\nu(\Delta_1(b))>0$. Note that such a path always exists since one can always consider $s^* \in \mathscr{S}_d$ such that $\nu(s^*)>0$ and then complete the path to $z$. Fix $n \geq 2$ and $w \in \Omega_n^{(0)}$. Denote by $\widetilde{w}^b = \Theta_{b}^{(n)}(w)$ the trajectory $w$ shifted by the prefix $b$. Let $M'>0$ be the constant in the logarithmic Lipschitz estimate
\begin{align*}
\big|\log\big(Q_p(s \mid r)\big)-\log\big(Q_p(s \mid r') \big) \big| \leq M'\|r-r'\|_{\infty}, \qquad s \in \mathscr{S}_d, \quad r,r'\in \Delta^{(d)},
\end{align*}
obtained in \eqref{eq_lemma_properties_log_Q_p_2}, and recall that $Q_p(s \mid r) \geq p_* >0$ uniformly in $s$ and $r$. Set
\begin{align*}
c_b \eqdef \nu\big( \Delta_1(b) \big) p_*^r \qquad \text{and} \qquad g_b(n) \eqdef c_b e^{-M'r(1+\log(n))}.
\end{align*}
We prove that this choice of $g_b(n)$ satisfies \eqref{eq_lemma_prefix_does_not_change}.

Using the trajectory-probability formula given in Remark~\ref{remark_formula_different_inital_distrib}, which follows Corollary~\ref{corollary_product_of_transitions}, together with the uniform lower bound on $Q_p(\cdot \mid \cdot)$ given in Lemma~\ref{lemma_lower_bound_transition} and a change of index, we obtain
\begin{align}
\PP_{d,p,\nu}\left( S_{[1:n+r]}=\widetilde{w}^b \right) &  = \nu\big( \Delta_1(b) \big)\prod_{k=2}^{r+1} Q_p \big( \Delta_k(\widetilde{w}^b) \mid \rho_k(\widetilde{w}^b)\big) \times \prod_{j=2}^{n-1} Q_p\big(\Delta_{j+r}(\widetilde{w}^b) \mid \rho_{j+r}(\widetilde{w}^b) \big) \notag \\
& \geq c_b\widetilde{\nu}(\Delta_1(w)) \prod_{j=2}^{n-1} Q_p\big(\Delta_{j+r}(\widetilde{w}^b) \mid \rho_{j+r}(\widetilde{w}^b) \big). \label{proof_lemma_prefix_does_not_change_1}
\end{align}
For every $2 \leq  j \leq n-1$, the first $j+r-1$ increments of $\widetilde{w}^b$ consist of the $r$ increments of $b$ followed by the first $j-1$ increments of $w$. Hence
\begin{align*}
\rho_{j+r}(\widetilde{w}^b)=\frac{r}{j+r-1}\rho_{r+1}(b) + \frac{j-1}{j+r-1}\rho_j(w).
\end{align*}
Since $\rho_j(w)$ and $\rho_{r+1}(b)$ are probability vectors, we deduce that
\begin{align}
\| \rho_{j+r}(\widetilde{w}^b)-\rho_j(w) \|_{\infty} & \leq \frac{r}{j+r-1}. \label{proof_lemma_prefix_does_not_change_2}
\end{align}
Moreover $\Delta_{j+r}(\widetilde{w}^b)=\Delta_j(w)$. Applying the logarithmic Lipschitz estimate given in \eqref{eq_lemma_properties_log_Q_p_2} and then \eqref{proof_lemma_prefix_does_not_change_2}, we obtain 
\begin{align}
Q_p\big( \Delta_{j+r}(\widetilde{w}^b) \mid \rho_{j+r}(\widetilde{w}^b) \big) & = Q_p\big( \Delta_{j}(w) \mid \rho_{j+r}(\widetilde{w}^b) \big) \notag \\
& \geq Q_p\big( \Delta_{j}(w) \mid \rho_j(w) \big) e^{-M' \frac{r}{j+r-1}}. \label{proof_lemma_prefix_does_not_change_3}
\end{align}
Substituting \eqref{proof_lemma_prefix_does_not_change_3} into \eqref{proof_lemma_prefix_does_not_change_1} and using the trajectory-probability formula under $\PP_{d,p,\widetilde{\nu}}$ yields
\begin{align*}
\PP_{d,p,\nu}\left(S_{[1:n+r]} =\widetilde{w}^b \right) & \geq c_b \prod_{j=2}^{n-1}e^{-M'\frac{r}{j+r-1}} \times \widetilde{\nu}(\Delta_1(w))\prod_{j=2}^{n-1} Q_p\big( \Delta_{j}(w) \mid \rho_j(w) \big) \\
&  \geq c_b e^{-M'r(\log(n)+1)} \PP_{d,p,\widetilde{\nu}}\big(S_{[1:n]}=w \big),
\end{align*}
and \eqref{eq_lemma_prefix_does_not_change} follows with the above choice of $g_b(n)$.
\end{proof}

For all $z \in \ZZ^d$, we denote by $T_z : x \mapsto x+z$ the translation by $z$. We also denote by $(T_z)_{\#}(\mu)$ the pushforward of $\mu$ by $T_z$. 

\begin{lemma}[Effect of a prefix on the empirical measure]\label{lemma_invariance_of_empirical_measure}
Consider $\mu \in \mathcal{P}(\ZZ^d)$ and $\varepsilon>0$. Suppose $b \in \Omega_{r+1}^{(0)}$ is a nearest-neighbour path from 0 to $z \in \ZZ^d$. Then, for all $n>2r/\varepsilon$ and $w \in \Omega_n(\mu,\varepsilon)$,
\begin{equation}\label{eq_lemma_invariance_of_empirical_measure}
\Theta_b^{(n)}(w) \in \Omega_{n+r}\big((T_z)_{\#}(\mu),2\varepsilon \big).
\end{equation}
\end{lemma}

\begin{proof}
By definition of $\Theta_b^{(n)}$, we have
\begin{align*}
\big[\Theta_b^{(n)}(w)\big]_j=\begin{cases}
b_j, & \text{if } j \in \dbl 1,r \dbr, \\
T_z(w_{j-r}), & \text{if } j \in \dbl r+1,n+r \dbr.
\end{cases}
\end{align*}
Thus, for all $w \in \Omega_n(\mu,\varepsilon)$, Equation~\eqref{eq_lemma_properties_KR_norm_of_proba} gives
\begin{align*}
\left\| \ell_{n+r}\left( \Theta_b^{(n)}(w) \right) - (T_z)_{\#}(\mu) \right\| & = \left\| \frac{1}{n+r} \sum_{j=1}^r \delta_{b_j} + \frac{1}{n+r} \sum_{j=1}^n \delta_{T_z(w_j)}  - (T_z)_{\#}(\mu) \right\| \\
& \leq \frac{r}{n+r} + \left \vert \frac{n}{n+r}-1 \right \vert + \Big\| (T_z)_{\#}\big(\ell_n(w)\big)-(T_z)_{\#}(\mu) \Big\|.
\end{align*}
The definition of the KR norm, given in \eqref{eq_def_BL_norm}, yields $\| (T_z)_{\#}(\ell_n(w))-(T_z)_{\#}(\mu)\|=\|\ell_n(w)-\mu\|$. Since $w \in \Omega_n(\mu,\varepsilon)$, we deduce that
\begin{align*}
\left\| \ell_{n+r}\left( \Theta_b^{(n)}(w) \right) - (T_z)_{\#}(\mu) \right\| < 2\frac{r}{n+r} + \varepsilon,
\end{align*}
and the result follows from the fact that $n$ is chosen sufficiently large so that $2r/n<\varepsilon$.
\end{proof}

\noindent We now use Lemma~\ref{lemma_prefix_does_not_change} and Lemma~\ref{lemma_invariance_of_empirical_measure} to prove the different symmetries of the rate function $I_{d,p}$.

\subsubsection{Initial-distribution invariance}\label{subsection_initial_distribution}

\begin{proposition}[Independence of the initial step distribution]\label{prop_invariance_under_initial_distribution}
For all $d \in \NN$, $p \in (0,1)$ and $\nu \in \mathcal{P}(\mathscr{S}_d)$, the sequence of empirical measures $(\ell_n)_{n \in \NN}$ satisfies, under $\PP_{d,p,\nu}$, a weak LDP with rate function $I_{d,p}$, given in Proposition~\ref{proposition_existence_weak_LDP_with_different_rate_fct}. In particular, this rate function does not depend on $\nu$.
\end{proposition}

\begin{remark}\label{remark_when_nu_and_nu_prime_full_support}
Suppose $\nu$ has full support and let $\overline{C}=2d\max_s \nu(s)$ and $\underline{C}=2d\min_s \nu(s)$. Then, for all $w \in \Omega_n^{(0)}$, one has
\begin{align*}
\underline{C} \, \PP_{d,p,\nu_{\text{unif}}}(S_{[1:n]}=w) \leq \PP_{d,p,\nu}(S_{[1:n]}=w) \leq \overline{C} \, \PP_{d,p,\nu_{\text{unif}}}(S_{[1:n]}=w),
\end{align*}
and the result of Proposition~\ref{prop_invariance_under_initial_distribution} follows by taking the logarithm, dividing by $n$, taking $\liminf_{n \to \infty}$ in the first inequality and $\limsup_{n \to \infty}$ in the second and finally taking $\varepsilon \to 0^+$. The only difficulty therefore lies in handling the case where the initial distribution is not of full support.
\end{remark}

\begin{proof}[Proof of Proposition~\ref{prop_invariance_under_initial_distribution}]
Fix $\mu \in \mathcal{P}(\ZZ^d)$ and let $\nu,\widetilde{\nu} \in \mathcal{P}(\mathscr{S}_d)$. Choose $s^* \in \mathscr{S}_d$ such that $\nu(s^*)>0$ and set $b=(0,s^*,0) \in \Omega_3^{(0)}$. Fix $\varepsilon>0$ and $n> 4/\varepsilon$. Lemma~\ref{lemma_invariance_of_empirical_measure} gives
\begin{align*}
\Theta_b^{(n)}(w) \in \Omega_{n+2}(\mu,2\varepsilon), \qquad w \in \Omega_n(\mu,\varepsilon).
\end{align*}
Since $\Theta_b^{(n)}$ is injective, we deduce from the above inclusion that
\begin{align*}
\PP_{d,p,\nu}\big(\ell_{n+2} \in B(\mu,2\varepsilon)\big) & \geq \sum_{w \in \Omega_n(\mu,\varepsilon)} \PP_{d,p,\nu}\left(S_{[1:n+2]}=\Theta_b^{(n)}(w) \right) \\
& \geq \sum_{w \in \Omega_n(\mu,\varepsilon)} g_b(n) \, \PP_{d,p,\widetilde{\nu}}\big( S_{[1:n]}=w \big),
\end{align*}
where the last inequality follows from Lemma~\ref{lemma_prefix_does_not_change}. Taking the logarithm, dividing by $n$ and taking $\liminf_{n \to \infty}$, we obtain
\begin{align*}
\liminf_{n \to \infty} \frac{1}{n} \log \PP_{d,p,\nu}\big(\ell_{n+2} \in B(\mu,2\varepsilon)\big) \geq \liminf_{n \to \infty} \frac{1}{n} \log g_b(n) + \liminf_{n \to \infty} \frac{1}{n} \log \PP_{d,p,\widetilde{\nu}}\big(\ell_{n} \in B(\mu,\varepsilon)\big).
\end{align*} 
Since $(n+2)/n \to1$, shifting the index from $n+2$ to $n$ does not change the lower exponential rate. Taking the infimum over $\varepsilon$ and using the property of $g_b(n)$ stated in Lemma~\ref{lemma_prefix_does_not_change}, we obtain
\begin{equation}\label{proof_invariance_under_initial_distribution_1}
\inf_{\varepsilon>0} \liminf_{n \to \infty} \frac{1}{n} \log \PP_{d,p,\nu}\big(\ell_{n} \in B(\mu,\varepsilon)\big) \geq \inf_{\varepsilon>0} \liminf_{n \to \infty} \frac{1}{n} \log \PP_{d,p,\widetilde{\nu}}\big(\ell_{n} \in B(\mu,\varepsilon)\big).
\end{equation}
For the reverse comparison, choose $\widetilde{s}^*\in\mathscr{S}_d$ such that $\widetilde{\nu}(\widetilde{s}^*)>0$ and set $\widetilde b=(0,\widetilde{s}^*,0)$. Applying the same argument with $\nu$ and $\widetilde{\nu}$ interchanged, and using $\limsup$ instead of $\liminf$, yields
\begin{equation}\label{proof_invariance_under_initial_distribution_2}
\inf_{\varepsilon>0} \limsup_{n \to \infty} \frac{1}{n} \log \PP_{d,p,\nu}\big(\ell_{n} \in B(\mu,\varepsilon)\big) \leq \inf_{\varepsilon>0} \limsup_{n \to \infty} \frac{1}{n} \log \PP_{d,p,\widetilde{\nu}}\big(\ell_{n} \in B(\mu,\varepsilon)\big).
\end{equation}
Let $\nu_{\text{unif}}$ be the uniform measure on $\mathscr{S}_d$. By Proposition~\ref{proposition_existence_weak_LDP_with_different_rate_fct} we have
\begin{align*}
\inf_{\varepsilon>0} \limsup_{n \to \infty} \frac{1}{n} \log \PP_{d,p,\nu_{\text{unif}}}\big(\ell_{n} \in B(\mu,\varepsilon)\big) & = -I_{d,p}(\mu) \\
& = \inf_{\varepsilon>0} \liminf_{n \to \infty} \frac{1}{n} \log \PP_{d,p,\nu_{\text{unif}}}\big(\ell_{n} \in B(\mu,\varepsilon)\big).
\end{align*}
Hence, applying \eqref{proof_invariance_under_initial_distribution_1} and \eqref{proof_invariance_under_initial_distribution_2} with $\widetilde{\nu}=\nu_{\text{unif}}$, together with \cite[Theorem~4.1.11]{DZ_LDP}, we deduce that under $\PP_{d,p,\nu}$, the sequence of empirical measures satisfies a weak LDP with rate function $I_{d,p}$ given in Proposition~\ref{proposition_existence_weak_LDP_with_different_rate_fct}.
\end{proof}

\subsubsection{Symmetry under lattice isometries}\label{subsubsection_invariance_under_lattice_isometries}

We now show that the weak LDP rate function $I_{d,p}$ is invariant under translations and signed permutations. By Proposition~\ref{prop_invariance_under_initial_distribution}, the rate function is invariant with respect to the initial step distribution. Thus, from now on, we revert to choosing the uniform distribution $\nu_{\textup{unif}}$ as the initial step distribution.

\begin{lemma}[Symmetry under translations]\label{lemma_invariance_shifts}
Let $z \in \ZZ^d$. Then, for every $\mu \in \mathcal{P}(\ZZ^d)$, 
\begin{equation}\label{eq_lemma_invariance_shifts}
I_{d,p}\left((T_z)_{\#}(\mu)\right)=I_{d,p}(\mu). 
\end{equation}
\end{lemma}

\begin{proof}
The assertion is immediate for $z=0$. For $z \neq 0$, choose a nearest-neighbour path $b \in \Omega_{r+1}^{(0)}$ from 0 to $z$ with $r \geq 1$ increments. Fix $\varepsilon>0$ and
$n>2r/\varepsilon$. Using the injectivity of $\Theta^{(n)}_b$ and Lemma~\ref{lemma_invariance_of_empirical_measure} and then Lemma~\ref{lemma_prefix_does_not_change}, there exists a sequence $(g_b(n))_{n \geq 2}$ such that $\log g_b(n)/n \to 0$ and
\begin{align*}
\PP_{d,p} \Big( \ell_{n+r} \in B\big( (T_z)_{\#}(\mu),2\varepsilon \big) \Big) & \geq \sum_{w \in \Omega_n(\mu,\varepsilon)} \PP_{d,p}\left( S_{[1:n+r]}=\Theta_b^{(n)}(w) \right) \\
& \geq g_b(n) \, \PP_{d,p} \big( \ell_n \in B(\mu,\varepsilon) \big).
\end{align*}
Taking the logarithm, dividing by $n$, taking $\liminf_{n \to \infty}$, taking $\inf_{\varepsilon>0}$ and then using the expression for $I_{d,p}$ given in \eqref{eq_weak_LDP_rate_function_expression_3} of Proposition~\ref{proposition_existence_weak_LDP_with_different_rate_fct}, we obtain
\begin{align*}
I_{d,p}\big( (T_z)_{\#}(\mu) \big) \leq I_{d,p}(\mu), \qquad z \in \ZZ^d.
\end{align*}
Writing $\mu=(T_{-z})_{\#}\big((T_z)_{\#}(\mu)\big)$, we obtain the reverse inequality and hence \eqref{eq_lemma_invariance_shifts} follows.
\end{proof}

We say that $R : \ZZ^d \longrightarrow \ZZ^d$ is a signed permutation of the coordinate axes of $\ZZ^d$ if there exist a permutation $\sigma \in \mathfrak{S}_d$ and signs $(s_1,\dots,s_d) \in \{ \pm 1\}^d$ such that
\begin{align*}
R(z_1,\dots,z_d)=\big( s_1 z_{\sigma(1)}, \dots , s_d z_{\sigma(d)}\big), \qquad (z_1,\dots,z_d) \in \ZZ^d.
\end{align*}

\begin{lemma}[Symmetry under signed permutations]\label{lemma_invariance_rotation}
Let $R$ be a signed permutation of the coordinate axes of $\ZZ^d$. Then, for every $\mu \in \mathcal{P}(\ZZ^d)$,
\begin{equation}\label{eq_lemma_invariance_rotation}
I_{d,p}\left(R_{\#}(\mu)\right)=I_{d,p}(\mu).
\end{equation}
\end{lemma}

\begin{proof}
Let $R$ be a signed permutation of the coordinate axes of $\ZZ^d$. For $w=(w_1,\ldots,w_n)\in\Omega_n^{(0)}$, set
\begin{equation*}
R\cdot w \eqdef \big(R(w_1),\ldots,R(w_n)\big).
\end{equation*}
By definition of $\Delta_j(w)$, given in \eqref{eq_Delta_word_j}, and the definition of $\rho_j$, given in \eqref{eq_def_empirical_distribution_rho}, we have
\begin{align*}
\Delta_j(R \cdot w)=R \big(\Delta_j(w)\big) \qquad \text{and} \qquad \rho_j(R \cdot w)=R_{\#} \rho_j(w).
\end{align*}
Since $[R_{\#}r]_{R(s)}=r_s$, we deduce that
\begin{align*}
Q_p \big( R(s) \mid R_{\#}r\big) = Q_p(s \mid r).
\end{align*}
Thus, Corollary~\ref{corollary_product_of_transitions} yields
\begin{equation}\label{eq_path_probability_rotation}
\PP_{d,p}\big(S_{[1:n]}=R\cdot w\big)=\PP_{d,p}\left(S_{[1:n]}=w\right), \qquad w\in\Omega_n^{(0)}.
\end{equation}
Furthermore,
\begin{align*}
\ell_n(R\cdot w)=\frac{1}{n}\sum_{k=1}^n\delta_{R(w_k)}=R_{\#}\big(\ell_n(w)\big).
\end{align*}
Since $R$ is an isometry of $\ZZ^d$, the KR metric is invariant under
pushforward by $R$. More precisely, for every
$\mu,\nu\in\mathcal{P}(\ZZ^d)$,
\begin{equation*}
\left\|R_{\#}(\mu)-R_{\#}(\nu)\right\|=\|\mu-\nu\|.
\end{equation*}
Indeed, composition with $R$ induces a bijection of the class of real-valued functions on $\ZZ^d$ bounded by $1$ and with Lipschitz constant at most $1$. It follows that
\begin{equation*}
\ell_n(R\cdot w)\in B\big(R_{\#}(\mu),\varepsilon\big) \quad \Longleftrightarrow \quad \ell_n(w)\in B(\mu,\varepsilon).
\end{equation*}
Consequently, using the bijection $w \mapsto R\cdot w$ and \eqref{eq_path_probability_rotation}, we obtain, for every $n\in\NN$ and $\varepsilon>0$,
\begin{equation*}
\PP_{d,p}\left( \ell_n \in B\big(R_{\#}(\mu) , \varepsilon\big)\right)=\PP_{d,p}\left(\ell_n\in B(\mu,\varepsilon) \right).
\end{equation*}
The local characterisation of the rate function given in \eqref{eq_weak_LDP_rate_function_expression_3} yields Equation~\eqref{eq_lemma_invariance_rotation} as claimed.
\end{proof}

\subsection{Comparison with the simple random walk}\label{subsection_comparison_to_SRW}

In this subsection, we show Proposition~\ref{prop_dimensional_dichotomy}. To do so, we compare the rate function $I_{d,p}$ of the empirical measures of the MERW to $I_{d,\text{DV}}$, the rate function of the empirical measures of the simple symmetric random walk on $\ZZ^d$. 

We first establish a general comparison between these two rate functions. 

\begin{proposition}[Uniform comparison with the DV rate function] \label{prop_uniform_comparison_with_DV}
For every $d \in \NN$ and $p \in (0,1)$,  
\begin{equation*}
    \| \mathcal{E}_p\|_{\infty} \eqdef \sup_{H \in \mathscr{H}^{(d)}} \big| \mathcal{E}_p(H) \big|<+\infty.
\end{equation*}
In addition,
\begin{equation}\label{eq_prop_bound_I_d_p}
    I_{d,\textup{DV}}(\mu)-\| \mathcal{E}_p\|_{\infty} \leq I_{d,p}(\mu) \leq I_{d,\textup{DV}}(\mu)+\| \mathcal{E}_p\|_{\infty}, \qquad \mu \in \mathcal{P}(\ZZ^d).
\end{equation}
In particular, the effective domains of $I_{d,p}$ and $I_{d,\textup{DV}}$ coincide.
\end{proposition}

\begin{proof}
By Lemma~\ref{lemma_properties_profile_space}, $\mathscr{H}^{(d)}$ is compact, and by Lemma~\ref{lemma_continuity_of_E_p}, the map $\mathcal{E}_p$ is continuous. Hence, $\| \mathcal{E}_p\|_{\infty}<+\infty$. Furthermore, since $\mathcal{E}_{1/(2d)} \equiv 0$, Equation~\eqref{eq_weak_LDP_rate_function_expression_2}, applied with $p=1/(2d)$, together with Corollary~\ref{corollary_identification_in_symmetric_case}, gives
\begin{equation*}
    \inf_{H \in \mathscr{H}^{(d)}} I_{\textup{A}}^{\textup{sym}}(\mu,H) = I_{d,\textup{sym}}(\mu) = I_{d,\textup{DV}}(\mu).
\end{equation*}
The desired inequalities then follow from \eqref{eq_weak_LDP_rate_function_expression_2} and the definition of $\| \mathcal{E}_p\|_{\infty}$.
\end{proof}

We now investigate whether these two rate functions actually coincide, distinguishing between the one-dimensional and higher-dimensional settings.

\subsubsection{The case \texorpdfstring{$d=1$}{d=1}}
\label{subsubsection_one_dimensional_rate_function}

Let us first consider the one-dimensional setting. We aim to show \eqref{eq_prop_dimensional_dichotomy_dim_1} which states that $I_{1,p}(\mu)=I_{1,\text{DV}}(\mu)$ for all $\mu \in \mathcal{P}(\ZZ)$. 

\begin{proof}[Proof of Equation~\eqref{eq_prop_dimensional_dichotomy_dim_1}]
Consider $H \in \mathscr{H}_{\text{bal}}^{(1)}$. Then, for all $t \in [0,1]$, we have
\begin{align*}
[H(t)]_{e_1}=[H(t)]_{-e_1} \qquad \text{and} \qquad [H(t)]_{e_1}+[H(t)]_{-e_1}=t.
\end{align*}
Thus,
\begin{align*}
\mathscr{H}^{(1)}_{\text{bal}}=\{H^*\}, \qquad H^*(t)=\left( \frac{t}{2},\frac{t}{2} \right), \quad t \in [0,1].
\end{align*}
By Lemma~\ref{lemma_unbalanced_profiles_proba_nul} and the expression for $I_{1,p}$ given in \eqref{eq_weak_LDP_rate_function_expression_1} of Proposition~\ref{proposition_existence_weak_LDP_with_different_rate_fct}, we deduce that
\begin{equation}\label{proof_identification_dim_one_1}
I_{1,p}(\mu)=I_{\text{A}}^{(p)}(\mu,H^*) \qquad \text{and} \qquad I_{1,\text{sym}}(\mu)=I_{\text{A}}^{\text{sym}}(\mu,H^*),
\end{equation}
where the second expression comes from the specialisation of the first to the case $p=1/(2d)$. Combining \eqref{proof_identification_dim_one_1} and \eqref{eq_proposition_expression_RL_rate_function} of  Proposition~\ref{proposition_expression_RL_rate_function}, we obtain for all $\mu \in \mathcal{P}(\ZZ)$,
\begin{align*}
I_{1,p}(\mu)& =I_{\text{A}}^{\text{sym}}(\mu,H^*)-\mathcal{E}_p(H^*) \\
& = I_{1,\text{sym}}(\mu)-\mathcal{E}_p(H^*).
\end{align*}
Since $H^*(t)=tr$, where $r=(1/2,1/2)$, and since, for $d=1$,
\begin{align*}
Q_p(s \mid r)=1/2, \qquad \forall \, p \in (0,1), \quad s \in \{ \pm e_1\},
\end{align*}
the definition of $\mathcal{E}_p$, given in \eqref{eq_def_energy_profile}, yields 
\begin{align*}
\mathcal{E}_p(H^*)=0.
\end{align*}
Therefore, $I_{1,p}(\mu)=I_{1,\text{sym}}(\mu)$. Using Corollary~\ref{corollary_identification_in_symmetric_case}, we deduce that
\begin{align*}
I_{1,p}(\mu)=I_{1,\text{DV}}(\mu), \qquad \forall \mu \in \mathcal{P}(\ZZ),
\end{align*}
and we obtain \eqref{eq_prop_dimensional_dichotomy_dim_1}.
\end{proof}

\begin{remark}\label{remark_explanation_proof_one_dim}
Note that the one-dimensional case is special because the position of the walker at a given time $n$ entirely determines the step distribution. Hence, the one-dimensional ERW is a time-inhomogeneous Markov chain. In particular, one can write the transition probabilities as follows. For all $n \geq 2$, $\sigma \in \{ \pm 1\}$, and $x \in \ZZ$ such that $\PP(S_n=x)>0$,
\begin{align*}
\PP_{1,p}\big( S_{n+1}=x + \sigma \big \vert S_n=x \big) = \frac{1}{2} + \sigma \frac{(2p-1)x}{2(n-1)}.
\end{align*}
As shown in  Corollary~\ref{corollary_bound_on_last_point_of_subword}, if $w \in \Omega_n(\mu,\varepsilon)$, then for large enough $n \in \NN$, one has $|x|<3\varepsilon n$. Thus, 
\begin{align*}
\left \vert \PP_{1,p}\big( S_{n+1}=x + \sigma \big \vert S_n=x \big) - \frac{1}{2} \right \vert \leq \frac{|2p-1|3 \varepsilon n}{2 (n-1)} < \frac{3}{2} \varepsilon.
\end{align*}
Thus, for $n$ large enough, the transitions resemble those of the simple random walk, which explains the identity $I_{1,p}=I_{1,\text{sym}}$.
\end{remark}

\subsubsection{The case $d \geq 2$}\label{subsubsection_quali_study_in_higher_dimension}

While no explicit expression for the rate function $I_{d,p}$ is available in  dimension $d \geq 2$, one can nevertheless show that the behaviour of $I_{d,p}$  differs markedly from the one-dimensional case. This is the content of Equation~\eqref{eq_prop_dimensional_dichotomy_dim_geq_1}, which asserts that  for $d \geq 2$ and $p \neq \frac{1}{2d}$,
\begin{equation*}
I_{d,p} \neq I_{d,\textup{DV}}.
\end{equation*}

We start by proving some preliminary lemmas. Let $e_1=(1,0,\dots,0) \in \ZZ^d$ and denote by $\pi_1:\ZZ^d\longrightarrow \ZZ$ the projection on the first coordinate. We consider the measure
\begin{equation*}\label{eq_def_mu_1}
\mu_1=\frac{1}{2} \delta_0 + \frac{1}{2} \delta_{e_1} \in \mathcal{P}(\ZZ^d).
\end{equation*}

\begin{lemma}\label{lemma_I_DV_for_mu_1}
The Donsker--Varadhan rate function evaluated at $\mu_1$ is given by
\begin{equation}\label{eq_lemma_I_DV_for_mu_1}
I_{d,\textup{DV}}(\mu_1)=\log(2d).
\end{equation}
\end{lemma}

\begin{proof}
Since the variational expression in \eqref{Donsker_Varadhan_our_case} is invariant under multiplication of $u$ by a positive constant and every admissible $u$ satisfies $0<\inf_{k \in \ZZ^d} u_k$, we may normalise $u$ so that $u_k \geq 1$ for all $k \in \ZZ^d$. Then, using the lower bound $u_k >0$ for all $k \in \ZZ^d \setminus \{0,e_1\}$, we deduce that
\begin{align*}
\sum_{k \in \ZZ^d} \mu_1(k) \log \left( \frac{u_k}{\frac{1}{2d} \sum_{s \in \mathscr{S}_d} u_{k+s}} \right) & < \log(2d) + \frac{1}{2} \log \left( \frac{u_0}{u_{e_1}} \right) + \frac{1}{2} \log \left( \frac{u_{e_1}}{u_{0}} \right) \\
& = \log(2d).
\end{align*}
Taking the supremum over all admissible $(u_k)_{k \in \ZZ^d}$ in the variational formula \eqref{Donsker_Varadhan_our_case}, we obtain $I_{d,\text{DV}}(\mu_1) \leq \log(2d)$. 

For the reverse inequality, fix $m \in \NN$ and define $(u_k^{(m)})_{k \in \ZZ^d}$ by 
\begin{align*}
u_0^{(m)}=m=u_{e_1}^{(m)} \qquad \text{and} \qquad u_k^{(m)}=1, \quad \forall k \in \ZZ^d \backslash \{0,e_1\}.
\end{align*}
This choice is admissible in \eqref{Donsker_Varadhan_our_case}, and gives
\begin{align*}
\log(2d) + \log \left( \frac{m}{m+2d-1} \right) \leq I_{d,\text{DV}}(\mu_1).
\end{align*}
Letting $m \to +\infty$ yields $\log(2d) \leq I_{d,\text{DV}}(\mu_1)$ and the lemma is proved.
\end{proof}

Define the profile $\widehat{H} \in \mathscr{H}^{(d)}_{\textup{bal}} \subseteq \mathscr{H}^{(d)}$ by
\begin{equation*}
[\widehat{H}(t)]_s=
\begin{cases}
\frac{t}{2}, & \text{if } s \in \{ \pm e_1\}, \\
0, & \text{otherwise},
\end{cases} \qquad t \in [0,1].
\end{equation*}
\noindent The following lemma shows that, asymptotically, the profile $\widehat{H}$ is the only profile that can be realised by trajectories whose empirical measures are close to $\mu_1$.

\begin{lemma}[Profile localisation at $\mu_1$]\label{lemma_controlling_profile_for_mu_1}
Fix $\varepsilon>0$. For every $n>\max \left\{1/\varepsilon,2 \right\}$ and $w \in \Omega_n(\mu_1,\varepsilon)$,
\begin{equation}\label{eq_lemma_controlling_profile_for_mu_1}
d_{\mathscr{H}}\big(H_n(w,\cdot),\widehat{H} \big)<4 \varepsilon.
\end{equation}
\end{lemma}

\begin{proof}
Fix $\varepsilon>0$. Consider $n > \max\left\{ \frac{1}{\varepsilon},2 \right\}$, $w \in \Omega_n(\mu_1,\varepsilon)$ and let 
\begin{align*}
A \eqdef \{0,e_1\}^c=\text{Supp}(\mu_1)^c.
\end{align*}
Using Lemma~\ref{lemma_occupation_meas_with_empirical_meas} and the fact that $\mu_1(A)=0$, we obtain 
\begin{equation}\label{eq_lemma_controlling_profile_for_mu_1_1}
\#\big\{j \in \dbl 1,n \dbr \, : \, w_j \in A \big\}<\varepsilon n.
\end{equation}
Suppose for the sake of contradiction that there exists $k \in \dbl 1,n \dbr$ such that $\|w_k\|_{\infty} \geq 2 \varepsilon n$. Since the trajectory is nearest-neighbour and $\varepsilon n >1$, this would imply that at least $\varepsilon n$ consecutive points of the trajectory lie outside $\{0,e_1\}$, contradicting \eqref{eq_lemma_controlling_profile_for_mu_1_1}. Thus,
\begin{equation}\label{eq_lemma_controlling_profile_for_mu_1_2}
\|w_k\|_{\infty}<2 \varepsilon n, \qquad k \in \dbl 1,n \dbr.
\end{equation}
If $\Delta_j(w) \notin \{ \pm e_1\}$, then either $w_j$ or $w_{j+1}$ belongs to $(\ZZ \times \{0_{d-1}\})^c \subseteq A$. Combining this with \eqref{eq_lemma_controlling_profile_for_mu_1_1}, we obtain
\begin{equation*} 
\# \big\{j \in \dbl 1,n-1 \dbr \, : \, \Delta_j(w) \notin \{ \pm e_1\} \big\} < 2\varepsilon n.
\end{equation*} 
Fix a mesh point $t_j=j/(n-1)$, with $j \in \dbl 0,n-1 \dbr$. Equation~\eqref{eq_def_profile_map} and the inequality just above give
\begin{equation}\label{eq_lemma_controlling_profile_for_mu_1_3}
0 \leq \sum_{s \notin \{ \pm e_1\}}[H_n(w,t_j)]_s < \frac{2 \varepsilon n}{n-1} \leq 4 \varepsilon, 
\end{equation}
where the last inequality comes from the fact that $n \geq 2$. Moreover,
\begin{align*}
&[H_n(w,t_j)]_{e_1}+[H_n(w,t_j)]_{-e_1}=t_j-\sum_{s \notin \{ \pm e_1\}} [H_n(w,t_j)]_{s}, \\
\text{and} \qquad & [H_n(w,t_j)]_{e_1}-[H_n(w,t_j)]_{-e_1}=\frac{\pi_1(w_{j+1})}{n-1}.
\end{align*}
Combining these two equations and then using \eqref{eq_lemma_controlling_profile_for_mu_1_2} and \eqref{eq_lemma_controlling_profile_for_mu_1_3} yields, for $\sigma \in \{ \pm 1\}$,
\begin{align*}
\left| [H_n(w,t_j)]_{\sigma e_1}-\frac{t_j}{2} \right| & \leq \frac{1}{2} \left| \frac{\|w_{j+1}\|_{\infty}}{n-1} + \sum_{s \notin \{ \pm e_1\}} [H_n(w,t_j)]_s \right| \\
& < 4 \varepsilon.
\end{align*}
The profile $H_n(w,\cdot)$ and the target $\widehat{H}$ are affine between consecutive mesh points. Hence, the same coordinatewise bounds hold for every $t \in [0,1]$, proving \eqref{eq_lemma_controlling_profile_for_mu_1}.
\end{proof}

\begin{lemma}[MERW--SRW comparison at $\mu_1$]\label{lemma_value_of_I_mu_1}
For all $p \in (0,1)$,
\begin{align}
I_{d,p}(\mu_1) & = -\log \left( p \frac{2d(d-1)}{2d-1} + \frac{d}{2d-1} \right) + I_{d,\textup{DV}}(\mu_1) \label{eq_lemma_value_of_I_mu_1} \\
& = -\log \left(p\frac{(d-1)}{2d-1}+\frac{1}{2(2d-1)}\right). \label{eq_other_value_of_I_mu_1}
\end{align}
\end{lemma}

\begin{proof}
Consider $H \in \mathscr{H}^{(d)}$ such that $H \neq \widehat{H}$.  Choose $\varepsilon>0$ sufficiently small so that $5 \varepsilon <d_{\mathscr{H}}(H,\widehat{H})$. By Lemma~\ref{lemma_controlling_profile_for_mu_1}, for all $n > \max \left\{ 1/\varepsilon,2 \right\}$ and $w \in \Omega_n(\mu_1,\varepsilon)$
\begin{align*}
d_{\mathscr{H}}\big( H_n(w,\cdot),H \big) & \geq d_{\mathscr{H}}\big( \widehat{H},H \big) - d_{\mathscr{H}}\big( H_n(w,\cdot),\widehat{H} \big) \\
& > \varepsilon.
\end{align*}
Thus, for all sufficiently large $n$,
\begin{align*}
\PP_{d,\text{sym}} \big( (\ell_n,H_n) \in B(\mu_1,\varepsilon) \times B_{\mathscr{H}}(H,\varepsilon) \big) =0.
\end{align*}
Taking the logarithm, dividing by $n$, taking the $\limsup_{n \to \infty}$ and then the infimum over $\varepsilon>0$ yields
\begin{align*}
I_{\text{A}}^{\text{sym}}(\mu_1,H)=+\infty, \qquad H \in \mathscr{H}^{(d)} \setminus \{ \widehat{H} \}.
\end{align*}
Thus, using the expression for $I_{d,p}$ given in \eqref{eq_weak_LDP_rate_function_expression_2}, we obtain for all $p \in (0,1)$,
\begin{equation}\label{eq_lemma_value_of_I_mu_1_1}
I_{d,p}(\mu_1)=I_{\text{A}}^{\text{sym}}(\mu_1,\widehat{H})-\mathcal{E}_p(\widehat{H}).
\end{equation}
As mentioned in Remark~\ref{remark_I_RL_p_when_p_is_1_slash_2d}, $\mathcal{E}_{1/(2d)} \equiv 0$. Applying the above equality to $p=1/(2d)$ yields
\begin{equation*}
I_{\text{A}}^{\text{sym}}(\mu_1,\widehat{H})=I_{d,\text{sym}}(\mu_1).
\end{equation*}
By Corollary~\ref{corollary_identification_in_symmetric_case},
$I_{d,\textup{sym}}=I_{d,\textup{DV}}$. Therefore,
\begin{equation}\label{eq_lemma_value_of_I_mu_1_2}
I_{\textup{A}}^{\textup{sym}}(\mu_1,\widehat{H})
=I_{d,\textup{DV}}(\mu_1).
\end{equation}
Since $\widehat{H}(t)=tr$ with $r_s=1/2$ for $s \in \{ \pm e_1\}$ and $r_s=0$ otherwise, Equation~\eqref{eq_remark_energy_profile_constant} gives
\begin{equation}\label{eq_lemma_value_of_I_mu_1_3}
\mathcal{E}_p(\widehat{H})=\log\left(2d\left( \frac{1}{2}p+\frac{1}{2}\frac{1-p}{2d-1} \right)\right)=\log \left( p \frac{2d(d-1)}{2d-1} + \frac{d}{2d-1} \right).
\end{equation}
Equation~\eqref{eq_lemma_value_of_I_mu_1} then follows by substituting \eqref{eq_lemma_value_of_I_mu_1_2} and \eqref{eq_lemma_value_of_I_mu_1_3} into \eqref{eq_lemma_value_of_I_mu_1_1}. Equation~\eqref{eq_other_value_of_I_mu_1} follows by substituting in the value of $I_{d,\text{DV}}(\mu_1)$ obtained in \eqref{eq_lemma_I_DV_for_mu_1} into \eqref{eq_lemma_value_of_I_mu_1}.
\end{proof}

\noindent The proof of \eqref{eq_prop_dimensional_dichotomy_dim_geq_1} is a direct consequence of Lemma~\ref{lemma_value_of_I_mu_1}.

\begin{proof}[Proof of Equation~\ref{eq_prop_dimensional_dichotomy_dim_geq_1}]
For all $p \in (0,1)$, we have
\begin{align*}
p \frac{2d(d-1)}{2d-1} + \frac{d}{2d-1} =1 \Longleftrightarrow p=\frac{1}{2d}.
\end{align*}
Thus, by \eqref{eq_lemma_value_of_I_mu_1} of Lemma~\ref{lemma_value_of_I_mu_1}, we deduce that $I_{d,p} \neq I_{d,\textup{DV}}$ as soon as $p \neq \frac{1}{2d}$.
\end{proof}

\begin{remark}\label{remark_extensions_not_equality_when_affine_hyperplane}
The mechanism used in this section is not specific to the measure $\mu_1$ considered here and should, in principle, extend to probability measures supported on an affine line $x + \mathrm{Vect}_{\ZZ}(e_i)$ of $\ZZ^d$, for $x \in \ZZ^d$ and $i \in \dbl 1,d \dbr$. Approximating such thinly supported measures is costly for a nearest-neighbour random walk. Indeed, to remain close to the line, the walker must take most of its steps in the direction of $e_i$, leaving it with essentially only $2$ effective choices out of $2d$ at each time step. For the symmetric random walk, this results in an effective transition probability of $\frac{1}{d}$. In the reinforcing regime $p>1/(2d)$, the elephant random walk, by contrast, benefits from its memory. With probability $p$, it repeats a past step, which is likely to be aligned with $e_i$ if the trajectory is already concentrated along the line. This self-reinforcing mechanism reduces the large deviation cost of approximating such measures. The phenomenon is illustrated by Lemma~\ref{lemma_value_of_I_mu_1}. Indeed, when $p>1/(2d)$, we have
\begin{equation*}
p \frac{2d(d-1)}{2d-1}+\frac{d}{2d-1}>1.
\end{equation*}
It follows from Equation~\eqref{eq_lemma_value_of_I_mu_1} that $I_{d,p}(\mu_1) < I_{d,\textup{DV}}(\mu_1)$. In other words, producing an empirical measure close to $\mu_1$ is less costly for the MERW than for the simple random walk.
\end{remark}

The previous remark states that the MERW uses its memory to reduce the cost of approximating thinly supported measures. This observation will be at the heart of the next section. 

\subsection{Convexity of the rate function}\label{subsection_phase_transition_for_MERW}

The aim of this section is to prove Proposition~\ref{prop_change_of_regime} which states that for $d \geq 2$, the rate function $I_{d,p}$ is convex if and only if $p \leq \frac{1}{2d}$. In Section~\ref{subsection_convex_for_small_p}, we show using the exponential tilting presented in Proposition~\ref{proposition_expression_RL_rate_function} that the rate function is convex for all $p \leq \frac{1}{2d}$. In Section~\ref{subsection_convexity_breaks}, we show that for two well-chosen measures $\mu_1$ and $\mu_2$, the rate function $I_{d,p}$ fails to be convex on the segment $[\mu_1,\mu_2]$ for all $p>1/(2d)$. Throughout this section, we always assume $d \geq 2$.

\subsubsection{\texorpdfstring{Convexity for $p \leq \frac{1}{2d}$}{Convexity for p}}\label{subsection_convex_for_small_p}

To show the convexity of $I_{d,p}$, we first show that the Ruelle--Lanford function $I_{\text{A}}^{(p)}$ for the augmented sequence $(\ell_n,H_n)_{n \in \NN}$ is convex. To do so we make use of the identity 
\begin{align*}
I_{\text{A}}^{(p)}(\mu,H)=I_{\text{A}}^{\text{sym}}(\mu,H)-\mathcal{E}_p(H), \qquad \forall (\mu, H) \in \mathcal{P}(\ZZ^d) \times \mathscr{H}^{(d)},
\end{align*}
established in Proposition~\ref{proposition_expression_RL_rate_function}. 

\begin{lemma}[Concavity of the energy map]\label{lemma_concavity_of_E_p}
Fix $p \in \left(0,\frac{1}{2d}\right]$. For all $H^{(1)},H^{(2)} \in \mathscr{H}^{(d)}$ and $\lambda \in [0,1]$, setting $H_{\lambda} \eqdef \lambda H^{(1)} + (1-\lambda) H^{(2)} \in \mathscr{H}^{(d)}$, we have
\begin{equation}\label{eq_lemma_concavity_of_E_p}
\mathcal{E}_p(H_{\lambda}) \geq \lambda \mathcal{E}_p\big(H^{(1)}\big) + (1-\lambda) \mathcal{E}_p\big(H^{(2)}\big).
\end{equation}
\end{lemma}

\begin{proof}
Let $C_1=p-(1-p)/(2d-1)$ and $C_2=(1-p)/(2d-1)$. For all $x \in [0,1]$, one obtains directly from \eqref{eq_def_F_p} and \eqref{eq_def_G_p} that
\begin{align*}
F_p''(x)=\frac{C_1}{C_1x+C_2} \qquad \text{and} \qquad G_p''(x)=\frac{C_1C_2}{(C_1x+C_2)^2}.
\end{align*}
Since $p \leq \frac{1}{2d}$, we have $C_1 \leq 0$. In addition, we have $C_2 \geq 0$ and $C_1x+C_2>0$. Therefore, $F_p$ and $G_p$ are concave on $[0,1]$. Since $\mathscr{H}^{(d)}$ is convex and for every $t>0$,
\begin{align*}
\frac{H_{\lambda}(t)}{t} = \lambda \frac{H^{(1)}(t)}{t} + (1-\lambda) \frac{H^{(2)}(t)}{t}, 
\end{align*}
Equation~\eqref{eq_lemma_concavity_of_E_p} then follows from the representation of $\mathcal{E}_p$, given in \eqref{eq_lemma_analytic_representation_of_profile_energy_map}, and the concavity of $F_p$ and $G_p$.
\end{proof}

\begin{lemma}[Convexity of the augmented rate function]\label{lemma_convexity_of_I_RL_low_p}
For all $p \in \left( 0,\frac{1}{2d} \right]$, the map $I_{\textup{A}}^{(p)}$ is convex on $\mathcal{P}(\ZZ^d) \times \mathscr{H}^{(d)}$.
\end{lemma}

\begin{proof}
Fix $p \leq \frac{1}{2d}$, $(\mu_1,H^{(1)}),(\mu_2,H^{(2)}) \in \mathcal{P}(\ZZ^d) \times \mathscr{H}^{(d)}$ and $\lambda \in [0,1]$. Set
\begin{align*}
\mu_{\lambda} \eqdef \lambda \mu_1 + (1-\lambda) \mu_2 \qquad \text{and} \qquad H_{\lambda} \eqdef \lambda H^{(1)} + (1-\lambda)H^{(2)}.
\end{align*}
If $I_{\text{A}}^{(p)}(\mu_i,H^{(i)})=+\infty$ for some $i \in \{1,2\}$, then there is nothing to prove. By Proposition~\ref{prop_good_stitching_implies_RL_equality},  $I_{\text{A}}^{\text{sym}}$ is convex, while Lemma~\ref{lemma_concavity_of_E_p} shows that $\mathcal{E}_p$ is concave on $\mathscr{H}^{(d)}$. Combining these facts with Proposition~\ref{proposition_expression_RL_rate_function}, we obtain that
\begin{align*}
I_{\text{A}}^{(p)}(\mu_{\lambda},H_{\lambda}) & = I_{\text{A}}^{\text{sym}}(\mu_{\lambda},H_{\lambda})-\mathcal{E}_p(H_{\lambda}) \\
& \leq  \lambda I_{\text{A}}^{\text{sym}}(\mu_1,H^{(1)}) + (1-\lambda) I_{\text{A}}^{\text{sym}}(\mu_2,H^{(2)}) - \big(\lambda \mathcal{E}_p(H^{(1)}) + (1-\lambda) \mathcal{E}_p(H^{(2)})\big) \\
& = \lambda I_{\text{A}}^{(p)}(\mu_1,H^{(1)}) + (1-\lambda) I_{\text{A}}^{(p)}(\mu_2,H^{(2)}). 
\end{align*}
Thus, the map $I_{\text{A}}^{(p)}$ is convex on $\mathcal{P}(\ZZ^d) \times \mathscr{H}^{(d)}$.
\end{proof}

\begin{proposition}[Convexity of the rate function]\label{prop_convexity_small_p}
For all $p \in \left(0,\frac{1}{2d}\right]$, the rate function $I_{d,p}$ is convex.
\end{proposition}

\begin{proof}
Fix $\mu_1,\mu_2 \in \mathcal{P}(\ZZ^d)$ and $\lambda \in [0,1]$, and set $\mu_{\lambda}=\lambda \mu_1 + (1-\lambda)\mu_2$. Since $I_{\textup{A}}^{(p)}$ is lower semicontinuous and $\mathscr{H}^{(d)}$ is compact, the infimum in the definition of $I_{d,p}(\mu)$ is attained \cite[Theorem~2.40]{Aliprantis}. Thus, there exist $\widehat{H}^{(1)},\widehat{H}^{(2)} \in \mathscr{H}^{(d)}$ such that
\begin{align*}
I_{\textup{A}}^{(p)}\big(\mu_i,\widehat{H}^{(i)}\big)=I_{d,p}(\mu_i), \qquad i \in \{1,2\}.
\end{align*}
Using the expression for $I_{d,p}$ given in \eqref{eq_weak_LDP_rate_function_expression_1} and then the convexity of $I^{(p)}_{\text{A}}$, given by Lemma~\ref{lemma_convexity_of_I_RL_low_p}, we obtain
\begin{align*}
I_{d,p}(\mu_{\lambda}) & \leq I_{\textup{A}}^{(p)}\big(\mu_{\lambda},\lambda \widehat{H}^{(1)} + (1-\lambda)\widehat{H}^{(2)}\big) \\
& \leq \lambda I_{\textup{A}}^{(p)}\big(\mu_1,\widehat{H}^{(1)}\big) + (1-\lambda) I_{\textup{A}}^{(p)}\big(\mu_2,\widehat{H}^{(2)}\big) \\
& = \lambda I_{d,p}(\mu_1) + (1-\lambda) I_{d,p}(\mu_2).
\end{align*} 
Thus, $I_{d,p}$ is convex on $\mathcal{P}(\ZZ^d)$.
\end{proof}

\subsubsection{\texorpdfstring{Non-convexity for $p>\frac{1}{2d}$}{Non-convexity for p>1/(2d)}}\label{subsection_convexity_breaks}

The aim of this section is to show the following result.
\begin{proposition}[Non-convexity in the reinforcing case]\label{prop_non_convexity_large_p}
For all $p \in \left( \frac{1}{2d},1 \right)$, the rate function $I_{d,p}$ is not convex on $\mathcal{P}(\ZZ^d)$.
\end{proposition}

To prove this result, we show that the rate function $I_{d,p}$ is non-convex on a specific segment. In addition to the measure $\mu_1$ defined above, we consider 
\begin{align*}
\mu_2 = \frac{1}{2} \delta_{2e_2} + \frac{1}{2} \delta_{3e_2} \in \mathcal{P}(\ZZ^d),
\end{align*}
where we recall that $e_2=(0,1,0,\dots,0)$. 

\begin{lemma}[Rate function at $\mu_2$]\label{lemma_value_of_I_d_p_mu_2}
For all $p \in (0,1)$,
\begin{equation}\label{eq_lemma_value_of_I_d_p_mu_2}
I_{d,p}(\mu_2)=I_{d,p}(\mu_1)=-\log  \left(p\frac{d-1}{2d-1}+\frac{1}{2(2d-1)}\right).
\end{equation}
\end{lemma}

\begin{proof}
The first equality in \eqref{eq_lemma_value_of_I_d_p_mu_2} is a direct consequence of Lemmas~\ref{lemma_invariance_shifts} and \ref{lemma_invariance_rotation} which state that the rate function $I_{d,p}$ is invariant under translations and signed permutations. The second equality follows from \eqref{eq_other_value_of_I_mu_1}.
\end{proof}

Consider the measure $\mu^* \in \mathcal{P}(\ZZ^d)$ defined by
\begin{equation*}\label{eq_def_mu_star}
\mu^* \eqdef \frac{1}{2} \mu_1 + \frac{1}{2} \mu_2 = \frac{1}{4} \delta_0 + \frac{1}{4} \delta_{e_1} + \frac{1}{4} \delta_{2e_2} + \frac{1}{4} \delta_{3e_2}.
\end{equation*}
To show that $I_{d,p}$ fails to be convex when $p>\frac{1}{2d}$, it is enough to show that for all $p \in \left( \frac{1}{2d},1 \right)$,
\begin{equation}\label{eq_aim_to_show_non_convexity}
I_{d,p}\left( \mu^* \right) > - \log  \left(p\frac{d-1}{2d-1}+\frac{1}{2(2d-1)}\right).
\end{equation}
Thus, we turn to deriving a lower bound on $I_{d,p}(\mu^*)$. Recall that by Proposition~\ref{proposition_expression_RL_rate_function}, we have
\begin{align*}
I_{\text{A}}^{(p)}(\mu^*,H) = I_{\text{A}}^{\text{sym}}(\mu^*,H)-\mathcal{E}_p(H), \qquad H \in \mathscr{H}^{(d)}.
\end{align*}
Thus, in the following, we find a lower bound for $I_{\text{A}}^{\text{sym}}(\mu^*,H)$ and an upper bound for $\mathcal{E}_p(H)$ uniformly in $H \in \mathscr{H}^{(d)}$. Consider 
\begin{align*}
A_1 \eqdef \{0,e_1\} \qquad \text{,} \qquad A_2 \eqdef \{2e_2,3e_2\} \qquad \text{and} \qquad A \eqdef A_1 \cup A_2 =\text{Supp}(\mu^*).
\end{align*} 

\begin{lemma}[Symmetric rate-function lower bound at $\mu^*$]\label{lemma_lower_bound_I_RL_sym_mu_star}
For all $H \in \mathscr{H}^{(d)}$,
\begin{equation}\label{eq_lemma_lower_bound_I_RL_sym_mu_star}
I_{\textup{A}}^{\textup{sym}}(\mu^*,H) \geq \log(2d).
\end{equation}
\end{lemma}

\begin{proof}
By contraction, see \eqref{eq_weak_LDP_rate_function_expression_1}, and Corollary~\ref{corollary_identification_in_symmetric_case},
\begin{equation}\label{proof_lemma_lower_bound_I_RL_sym_mu_star}
I_{\text{A}}^{\text{sym}}(\mu^*,H) \geq I_{d,\text{sym}}(\mu^*) = I_{d,\text{DV}}(\mu^*).
\end{equation}
Let $M \geq 1$ and consider $(\widehat{u}_k)_{k \in \ZZ^d}$ defined by
\begin{align*}
\widehat{u}_k=
\begin{cases}
M, & \text{if } k \in A, \\
1, & \text{otherwise}.
\end{cases}
\end{align*}
Using the definition of $I_{d,\text{DV}}$, given in \eqref{Donsker_Varadhan_our_case}, and the fact that each vertex of $A$ has exactly one neighbour in $A$, one obtains
\begin{align*}
I_{d,\text{DV}}(\mu^*) & \geq \sum \limits_{k \in \ZZ^d} \mu^*(k) \log \left( \frac{ \widehat{u}_k}{\frac{1}{2d}\sum_{s \in \mathscr{S}_d}\widehat{u}_{k+s}} \right) \\
& = \log(2d) + \frac{1}{4} \sum_{k \in A} \log \left( \frac{M}{M+2d-1} \right).
\end{align*}
Taking $M$ to infinity gives $I_{d,\text{DV}}(\mu^*) \geq \log(2d)$ and \eqref{eq_lemma_lower_bound_I_RL_sym_mu_star} follows from \eqref{proof_lemma_lower_bound_I_RL_sym_mu_star}.
\end{proof}

\noindent Define the point $x_{\text{end}} \in \Delta^{(d)}$ by 
\begin{equation}\label{eq_def_x_end}
[x_{\text{end}}]_s=
\begin{cases}
1/4, & \text{if } s \in \{ \pm e_1, \pm e_2 \} \\
0, &\text{otherwise}
\end{cases}.
\end{equation}
\begin{lemma}[End-point constraint]\label{lemma_condition_on_profiles}
Let $p \in (0,1)$. For all $H \in \mathscr{H}^{(d)}$ such that $H(1) \neq x_{\textup{end}}$,
\begin{equation}\label{eq_lemma_condition_on_profiles}
I_{\textup{A}}^{(p)}(\mu^*,H)=+\infty.
\end{equation}
\end{lemma}

\begin{proof}
To show this result, we introduce for every $n \in \NN$, every $w \in \Omega_{n}$ and every $x,y \in \ZZ^d$, the quantity
\begin{align*}
N_{(x,y)}^{\text{tr}}(w) \eqdef \# \big\{ j \in \dbl 1,n-1 \dbr \, : \, (w_j,w_{j+1})=(x,y) \big\},
\end{align*}
which counts the number of times the trajectory $w$ makes the transition from $x$ to $y$. Note that if $x \not \sim y$, then $N_{(x,y)}^{\text{tr}}(w)=0$. Then, let 
\begin{align*}
\mathcal{G} \eqdef \big\{ (0,e_1),(e_1,0),(2e_2,3e_2),(3e_2,2e_2) \big\},
\end{align*}
be the set of good transitions, meaning a transition along one of the two edges of $A_1$ and $A_2$. Fix $\varepsilon>0$ and $n>3/\varepsilon$. If $(w_j,w_{j+1}) \notin \mathcal{G}$, then either $w_j \notin A$ or $w_{j+1} \notin A$. Thus, for all $w \in \Omega_n(\mu^*,\varepsilon)$, Lemma~\ref{lemma_occupation_meas_with_empirical_meas} yields
\begin{equation}\label{proof_lemma_condition_on_profiles_1}
\sum_{x,y \in \ZZ^d \, : \, (x,y) \notin \mathcal{G}} N_{(x,y)}^{\text{tr}}(w) \leq 2 N_{A^c}(w) < 2 \varepsilon n.
\end{equation} 
For all $s \notin \{ \pm e_1, \pm e_2\}$ and all $x \in \ZZ^d$, note that $(x,x+s) \notin \mathcal{G}$. Thus, for all $s \notin \{ \pm e_1, \pm e_2\}$, by decomposing the number of $s$-steps according to the starting point of the step, we obtain
\begin{align*}
C_n^{(s)}(w) & = \sum_{x \in \ZZ^d} N^{\text{tr}}_{(x,x+s)}(w) \\
& \leq \sum_{\substack{x,y \in \ZZ^d \, : \, (x,y) \notin \mathcal{G}}} N_{(x,y)}^{\text{tr}}(w) \\
& <2 \varepsilon n.
\end{align*}
Thus, we deduce that for all $w \in \Omega_n(\mu^*,\varepsilon)$ and all $s \notin \{ \pm e_1, \pm e_2\}$, 
\begin{equation}\label{proof_lemma_condition_on_profiles_2}
[H_n(w,1)]_s = \frac{C_n^{(s)}(w)}{n-1} < \frac{2 \varepsilon n}{n-1} \leq 4 \varepsilon.
\end{equation}
We now turn towards the case $s \in \{ \pm e_1, \pm e_2\}$. We start with $s=e_1$. First, since $w_1=0$, note that
\begin{align}\label{proof_lemma_condition_on_profiles_extra}
N_{e_1}(w)=N_{(0,e_1)}^{\text{tr}}(w) + \sum_{s \in \mathscr{S}_d \backslash \{+e_1\}} N_{(e_1-s,e_1)}^{\text{tr}}(w).
\end{align}
For all $s \in \mathscr{S}_d \backslash \{+e_1\}$, $(e_1-s,e_1) \notin \mathcal{G}$. Thus, using Lemma~\ref{lemma_occupation_meas_with_empirical_meas} together with the fact that $\mu^*(\{e_1\})=1/4$ and then \eqref{proof_lemma_condition_on_profiles_1}, we obtain
\begin{align}
\left|\frac{N^{\text{tr}}_{(0,e_1)}(w)}{n}-\frac{1}{4} \right| & < \frac{1}{n} \sum_{s \in \mathscr{S}_d \backslash \{+e_1\}} N_{(e_1-s,e_1)}^{\text{tr}}(w) + \varepsilon  \notag \\
& <3 \varepsilon. \label{proof_lemma_condition_on_profiles_3}
\end{align}
Thus, we deduce that for all $w \in \Omega_n(\mu^*,\varepsilon)$,
\begin{equation*}
\left \vert \big[ H_n(w,1) \big]_{+e_1} - \frac{1}{4} \right \vert = \left \vert \frac{N_{(0,e_1)}^{\text{tr}}(w)}{n-1} + \frac{1}{n-1} \sum_{x \in \ZZ^d \backslash \{0\}} N_{(x,x+e_1)}^{\text{tr}}(w) - \frac{1}{4} \right \vert.
\end{equation*}
Using \eqref{proof_lemma_condition_on_profiles_3} and the fact that $(x,x+e_1) \notin \mathcal{G}$ for all $x \in \ZZ^d \backslash \{0\}$, we obtain
\begin{align*}
\left \vert \big[ H_n(w,1) \big]_{+e_1} - \frac{1}{4} \right \vert & \leq \frac{1}{n} + 3 \varepsilon + \frac{1}{n-1}\sum_{x,y \in \ZZ^d \, : \, (x,y) \notin \mathcal{G}} N_{(x,y)}^{\text{tr}}(w) \notag \\
& < \frac{1}{n} + 3 \varepsilon + \frac{2 \varepsilon n}{n-1} \\
& < 8 \varepsilon, 
\end{align*}
where the second inequality follows from \eqref{proof_lemma_condition_on_profiles_1} and the last inequality from the fact that $n$ is chosen sufficiently large so that $n>3/\varepsilon$. By the same arguments, one can show that for all $s \in \{ \pm e_1,\pm e_2\}$, 
\begin{equation}\label{proof_lemma_condition_on_profiles_5}
\left \vert \big[ H_n(w,1) \big]_{s} - \frac{1}{4} \right \vert < 8 \varepsilon.
\end{equation}
The only modification that is required is that for $s=-e_1$, instead of \eqref{proof_lemma_condition_on_profiles_extra}, we use
\begin{align*}
N_{0}(w)-1=N_{(e_1,0)}^{\text{tr}}(w) + \sum_{s \in \mathscr{S}_d \backslash \{-e_1\}} N_{(0-s,0)}^{\text{tr}}(w),
\end{align*}
since $w_1=0$. This does not, however, modify \eqref{proof_lemma_condition_on_profiles_5}. Combining \eqref{proof_lemma_condition_on_profiles_2} with \eqref{proof_lemma_condition_on_profiles_5}, we deduce that for all $\varepsilon>0$, $n>3/\varepsilon$ and $w \in \Omega_n(\mu^*,\varepsilon)$,
\begin{equation}\label{proof_lemma_condition_on_profiles_6}
\big\| H_n(w,1) - x_{\text{end}} \big\|_{\infty} <8 \varepsilon.
\end{equation}

Consider now $H \in \mathscr{H}^{(d)}$ such that $H(1) \neq x_{\textup{end}}$. Choose $\varepsilon>0$ such that
\begin{align*}
9 \varepsilon < \big\|x_{\text{end}}-H(1)\big\|_{\infty}.
\end{align*}
By \eqref{proof_lemma_condition_on_profiles_6}, for every $n>3/ \varepsilon$ and $w \in \Omega_n(\mu^*,\varepsilon)$,
\begin{align*}
d_{\mathscr{H}}(H_n(w,\cdot),H) & \geq \big\|H_n(w,1)-H(1)\big\|_{\infty} \\
 & \geq  \big\|x_{\text{end}}-H(1)\big\|_{\infty} - \big\| x_{\text{end}} - H_n(w,1) \big\|_{\infty} \\
& > \varepsilon,
\end{align*}
Thus, 
\begin{align*}
\PP_{d,p}\big( (\ell_n,H_n) \in B(\mu^*,\varepsilon) \times B_{\mathscr{H}}(H,\varepsilon) \big)=0, \qquad n>\frac{3}{\varepsilon}.
\end{align*}
Taking the logarithm, dividing by $n$, taking the $\limsup_{n \to \infty}$ and then the infimum over $\varepsilon>0$, we deduce that
\begin{align*}
I_{\text{A}}^{(p)}(\mu^*,H)=-\overline{s}_p(\mu^*,H)=+\infty,
\end{align*}
and the lemma is proved.
\end{proof}

\noindent We introduce the set of profiles
\begin{equation}\label{eq_def_H_star}
\mathscr{H}^* \eqdef \Big\{ H \in \mathscr{H}_{\text{bal}}^{(d)}  \, : \, H(1)=x_{\text{end}} \Big\} \subseteq \mathscr{H}^{(d)},
\end{equation}
which consists of profiles that are balanced, as defined in \eqref{eq_def_H_bal}, and have the prescribed endpoint $x_{\text{end}} \in \Delta^{(d)}$. The following lemma gives an upper bound on the energy of such profiles.

\begin{lemma}[Energy bound for admissible profiles]\label{lemma_bound_energy_profiles_in_F_star}
For all $p \in \left(\frac{1}{2d},1 \right)$, there exists $\delta_{p}>0$ such that for all $H \in \mathscr{H}^*$,
\begin{equation}\label{eq_lemma_bound_energy_profiles_in_F_star}
\mathcal{E}_p(H) \leq \log \left( 2d \left( p \frac{d-1}{2d-1} + \frac{1}{2(2d-1)} \right) \right)-\delta_p.
\end{equation}
\end{lemma}

\begin{proof}
Consider $H \in \mathscr{H}^*$. Since $H(1)=x_{\text{end}}$, we have
\begin{align*}
\sum_{s \in \mathscr{S}_d} F_p \big( [H(1)]_s \big) = 4 F_p\left( \frac{1}{4} \right) + (2d-4)F_p(0)=4 F_p\left( \frac{1}{4} \right),
\end{align*}
where the last equality uses the fact that $F_p(0)=0$. For $t>0$, let
\begin{align*}
\alpha_j(t) \eqdef \frac{[H(t)]_{e_j}}{t}, \qquad j \in \dbl 1, d \dbr.
\end{align*}
Since $H$ is balanced and $\sum_{s}[H(t)]_s=t$,
\begin{align*}
\big( \alpha_j(t) \big)_{j=1}^d \in \Delta^{(d)}_{1/2} \eqdef \left\{ (x_1, \dots , x_d) \in \RR_+^d \, : \sum_{i=1}^d x_i= \frac{1}{2} \right\}.
\end{align*}
Moreover,
\begin{equation}\label{proof_bound_energy_profiles_in_F_star_1}
\sum_{s \in \mathscr{S}_d} G_p\left(\frac{[H(t)]_s}{t}\right)=2 \sum_{j=1}^d G_p\big(\alpha_j(t)\big).
\end{equation}
Let $C_1=p-(1-p)/(2d-1)$ and $C_2=(1-p)/(2d-1)$. Since $p>1/(2d)$, we have $C_1>0$. Thus, for all $x \in [0,1]$,
\begin{align*}
F_p''(x)=\frac{C_1}{C_1x+C_2} \qquad \text{and} \qquad G_p''(x)=\frac{C_1C_2}{(C_1x+C_2)^2} > 0,
\end{align*}
and $F_p$ and $G_p$ are strictly convex on $[0,1]$. Since $F_p(0)=0$, the strict convexity of $F_p$ yields the existence of a $\delta_p>0$ such that
\begin{equation}\label{proof_bound_energy_profiles_in_F_star_2}
F_p \left( \frac{1}{4} \right) \leq \frac{1}{2} F_p \left( \frac{1}{2} \right)-\delta_p.
\end{equation} 
In addition, $(x_1,\dots,x_d) \mapsto \sum_{j=1}^d G_p(x_j)$ is strictly convex on $\Delta^{(d)}_{1/2}$ and its maximum is attained on the extremal points of the set, see \cite[Theorem~3.4.7]{Niculescu_convexity}. Hence, for all $t \in (0,1]$,
\begin{equation}\label{proof_bound_energy_profiles_in_F_star_3}
\sum_{j=1}^d G_p\big(\alpha_j(t)\big) \leq G_p\left( \frac{1}{2} \right).
\end{equation}
Combining \eqref{proof_bound_energy_profiles_in_F_star_1} and \eqref{proof_bound_energy_profiles_in_F_star_3}, we deduce that
\begin{align*}
\int_0^1 \sum_{s \in \mathscr{S}_d} G_p\left(\frac{[H(t)]_s}{t}\right) \mathrm{d}t \leq 2G_p\left( \frac{1}{2} \right).
\end{align*}
Using the second representation of $\mathcal{E}_p$ in Lemma~\ref{lemma_analytic_representation_of_profile_energy_map}, we obtain for all $H \in \mathscr{H}^*$,
\begin{align*}
\mathcal{E}_p(H) \leq 4F_p \left( \frac{1}{4} \right) +2G_p\left( \frac{1}{2} \right).
\end{align*}
Using \eqref{proof_bound_energy_profiles_in_F_star_2} and the relation linking $L_p$, $F_p$ and $G_p$, given in \eqref{eq_def_G_p}, we deduce that
\begin{align*}
\mathcal{E}_p(H) \leq 2F_p \left( \frac{1}{2} \right) + 2G_p \left( \frac{1}{2} \right) - \delta_p = L_p \left( \frac{1}{2} \right) - \delta_p,
\end{align*}
and Equation~\eqref{eq_lemma_bound_energy_profiles_in_F_star} follows from the definition of $L_p$, given in \eqref{eq_def_L_p}. 
\end{proof}

\noindent We can now turn to the proof of Proposition~\ref{prop_non_convexity_large_p} which states that the rate function $I_{d,p}$ fails to be convex when $p>1/(2d)$.
\begin{proof}[Proof of Proposition~\ref{prop_non_convexity_large_p}]
Fix $p > 1/(2d)$. By combining Proposition~\ref{proposition_expression_RL_rate_function}, Lemma~\ref{lemma_unbalanced_profiles_proba_nul} and Lemma~\ref{lemma_condition_on_profiles}, we obtain that for all $H \notin \mathscr{H}^*$,
\begin{align*}
I_{\text{A}}^{(p)}(\mu^*,H)=+\infty.
\end{align*}
Thus, by Proposition~\ref{proposition_existence_weak_LDP_with_different_rate_fct}, we deduce that
\begin{align*}
I_{d,p}(\mu^*)=\inf_{H \in \mathscr{H}^*} \Big\{ I_{\text{A}}^{\text{sym}}(\mu^*,H) - \mathcal{E}_p(H) \Big\}.
\end{align*}
Using the bounds obtained in Lemma~\ref{lemma_lower_bound_I_RL_sym_mu_star} and Lemma~\ref{lemma_bound_energy_profiles_in_F_star}, we deduce that
\begin{align*}
I_{d,p}(\mu^*) & \geq \log(2d) - \left( \log \left( 2d \left( p\frac{d-1}{2d-1} + \frac{1}{2(2d-1)} \right) \right)-\delta_p \right) \\
& > -\log \left( p\frac{d-1}{2d-1} + \frac{1}{2(2d-1)} \right).
\end{align*}
Combining the above inequality with Lemma~\ref{lemma_value_of_I_d_p_mu_2}, we deduce that
\begin{align*}
I_{d,p}(\mu^*)>\frac{1}{2} I_{d,p}(\mu_1) + \frac{1}{2} I_{d,p}(\mu_2).
\end{align*}
Since $\mu^* = \frac{1}{2} \mu_1 + \frac{1}{2} \mu_2$, we deduce that $I_{d,p}$ fails to be convex on the segment $[\mu_1,\mu_2]$.
\end{proof}

\begin{remark}\label{remark_transition_convexity_to_concavity}
By Corollary~\ref{corollary_identification_in_symmetric_case}, we have the identity $I_{d,1/(2d)}=I_{d,\text{DV}}$. Using \eqref{Donsker_Varadhan_our_case} and the fact that $\text{Supp}(\mu_1) \not\sim \text{Supp}(\mu_2)$, one can show that $I_{d,\text{DV}}$ is affine on the segment $[\mu_1,\mu_2]$. Thus, on  the segment $[\mu_1,\mu_2]$, the rate function is convex when $p<1/(2d)$, affine at the crossover $p=1/(2d)$ and fails to be convex when $p>1/(2d)$. Note that one can use the invariance of $I_{d,p}$ under translations and signed permutations, shown in Lemmas~\ref{lemma_invariance_shifts} and \ref{lemma_invariance_rotation}, to obtain an infinite number of segments where this crossover appears. 
\end{remark}

\begin{remark}\label{remark_interpretation_proof}
Let us give a probabilistic interpretation of the proof of Proposition~\ref{prop_non_convexity_large_p}. By Lemma~\ref{lemma_exponential_tilting_wrt_SRW}, for all $w \in \Omega_n(\mu^*,\varepsilon)$,
\begin{align*}
\PP_{d,p}\big(S_{[1:n]}=w\big) \lesssim \PP_{d,\text{sym}}\big(S_{[1:n]}=w\big)e^{n \mathcal{E}_p\left( H_n(w) \right)}.
\end{align*}
Using Lemma~\ref{lemma_bound_energy_profiles_in_F_star}, which is an asymptotic bound on $\mathcal{E}_p\big(H_n(w)\big)$, we obtain
\begin{align*}
& \frac{1}{n} \log \PP_{d,p}\big( \ell_n \in B(\mu^*,\varepsilon) \big) \\
& \qquad \qquad \quad \lesssim \left( \frac{1}{n} \log \# \Omega_n(\mu^*,\varepsilon) + \log \left( \frac{1}{2d} \right) \right) + \left( \log \left( 2d \left(p \frac{d-1}{2d-1} + \frac{1}{2(2d-1)} \right) \right)-\delta_p \right).
\end{align*}
Lemma~\ref{lemma_lower_bound_I_RL_sym_mu_star} could have also been shown by giving a combinatorial bound on $\#\Omega_n(\mu^*,\varepsilon)$. Indeed, the support of $\mu^*$ is partitioned into two connected components
\begin{align*}
\text{Supp}(\mu^*) = \text{Supp}(\mu_1) \cup \text{Supp}(\mu_2) \qquad \text{with} \qquad \text{Supp}(\mu_1) \not\sim \text{Supp}(\mu_2).
\end{align*}
Thus, every transition between these two connected components forces the trajectory to visit the set $\ZZ^d \backslash \text{Supp}(\mu^*)$, which cannot happen too frequently for a trajectory in $\Omega_n(\mu^*,\varepsilon)$. Hence, one can bound the number of such transitions and thereby show that $\#\Omega_n(\mu^*,\varepsilon)$ is subexponential. We thereby obtain
\begin{align*}
\frac{1}{n} \log \PP_{d,p}\big( \ell_n \in B(\mu^*,\varepsilon) \big) \lesssim \log \left(p \frac{d-1}{2d-1} + \frac{1}{2(2d-1)} \right)-\delta_p,
\end{align*}
which gives us the required strict lower bound on $I_{d,p}(\mu^*)$.
\end{remark}

\begin{proof}[Proof of Proposition~\ref{prop_change_of_regime}] Combining Propositions~\ref{prop_convexity_small_p} and \ref{prop_non_convexity_large_p}, one obtains that $I_{d,p}$ is convex if and only if $p \in \left(0,\frac{1}{2d} \right]$. Thus, Proposition~\ref{prop_change_of_regime} is proved.
\end{proof}

\section{The stitching procedure}\label{section_proof_of_existence_of_stitching_procedure}

This section proves Proposition~\ref{prop_good_stitching_procedure}, the decoupling criterion used in Section~\ref{subsection_decoupling_method}. Our approach shares some similarities with the construction in \cite{Daures}, where a stitching procedure is introduced for the empirical measures of discrete-space Markov processes. In the present setting, however, we need to construct a stitching procedure that simultaneously controls the empirical measure and the profile, that is, the pair $(\ell_n,H_n)$. 

We begin by choosing suitable parameters and defining the stitching map $\psi$ in Section~\ref{subsection_defining_stitching_map}. Then, in Section~\ref{subsection_verifying_prop_of_stitching}, we show that $\psi$ satisfies the properties required by Proposition~\ref{prop_good_stitching_procedure}. More precisely, in Section~\ref{subsubsection_inclusion_property}, we show that the trajectory resulting from this stitching procedure has empirical measure and profile arbitrarily close to the prescribed targets $(\mu_1+\mu_2)/2$ and $(H^{(1)}+H^{(2)})/2$, provided the parameters are chosen appropriately. We then show in Section~\ref{subsubsection_combinatorial_bound} that each preimage $\psi^{-1}(\{\overline{w}\})$ is not too large. We derive subadditive-type estimates in Section~\ref{subsubsection_probability_comparaison} by comparing the probability of realising the longer stitched trajectory against the product of the probabilities of its constituent components. Finally, we prove Proposition~\ref{prop_good_stitching_procedure} in Section~\ref{subsection_existence_RL_function_for_couple}.

In the rest of this section, we fix $\mu_1,\mu_2 \in \mathcal{P}(\ZZ^d)$ and $H^{(1)},H^{(2)} \in \mathscr{H}^{(d)}$.  

\subsection{Construction of the stitching map}\label{subsection_defining_stitching_map}

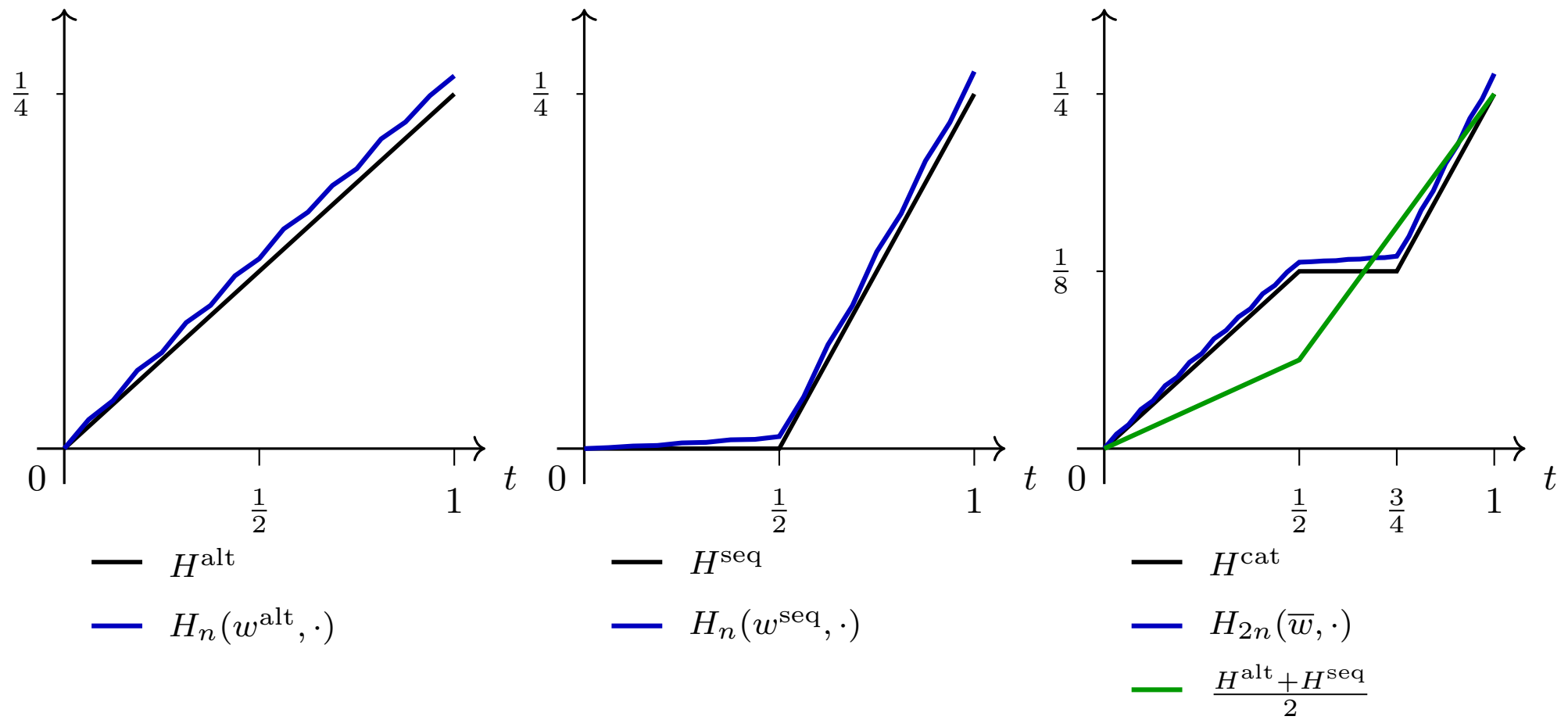
\begin{figure}[htbp]
\centering
\begin{tikzpicture}[
    x=3.3cm,
    y=12cm,
    font=\small,
    axis/.style={->, semithick},
    original/.style={black, line width=1.05pt},
    approx/.style={blue!75!black, line width=1.15pt},
    average/.style={green!60!black, line width=1.15pt}
]
\pgfmathsetmacro{\bA}{0}
\pgfmathsetmacro{\bB}{0.0048}
\pgfmathsetmacro{\bC}{0.0026}
\pgfmathsetmacro{\bD}{0.0083}
\pgfmathsetmacro{\bE}{0.0051}
\pgfmathsetmacro{\bF}{0.0107}
\pgfmathsetmacro{\bG}{0.0072}
\pgfmathsetmacro{\bH}{0.0124}
\pgfmathsetmacro{\bI}{0.0089}
\pgfmathsetmacro{\bJ}{0.0142}
\pgfmathsetmacro{\bK}{0.0105}
\pgfmathsetmacro{\bL}{0.0136}
\pgfmathsetmacro{\bM}{0.0098}
\pgfmathsetmacro{\bN}{0.0152}
\pgfmathsetmacro{\bO}{0.0115}
\pgfmathsetmacro{\bP}{0.0144}
\pgfmathsetmacro{\bQ}{0.0128}

\pgfmathsetmacro{\uA}{0}
\pgfmathsetmacro{\uB}{0.0007}
\pgfmathsetmacro{\uC}{0.0018}
\pgfmathsetmacro{\uD}{0.0022}
\pgfmathsetmacro{\uE}{0.0040}
\pgfmathsetmacro{\uF}{0.0044}
\pgfmathsetmacro{\uG}{0.0062}
\pgfmathsetmacro{\uH}{0.0065}
\pgfmathsetmacro{\uI}{0.0085}

\pgfmathsetmacro{\uJ}{0.0052}
\pgfmathsetmacro{\uK}{0.0108}
\pgfmathsetmacro{\uL}{0.0070}
\pgfmathsetmacro{\uM}{0.0138}
\pgfmathsetmacro{\uN}{0.0094}
\pgfmathsetmacro{\uO}{0.0154}
\pgfmathsetmacro{\uP}{0.0112}
\pgfmathsetmacro{\uQ}{0.0158}

\pgfmathsetmacro{\fA}{\bA}
\pgfmathsetmacro{\fB}{0.015625+\bB}
\pgfmathsetmacro{\fC}{0.03125+\bC}
\pgfmathsetmacro{\fD}{0.046875+\bD}
\pgfmathsetmacro{\fE}{0.0625+\bE}
\pgfmathsetmacro{\fF}{0.078125+\bF}
\pgfmathsetmacro{\fG}{0.09375+\bG}
\pgfmathsetmacro{\fH}{0.109375+\bH}
\pgfmathsetmacro{\fI}{0.125+\bI}
\pgfmathsetmacro{\fJ}{0.140625+\bJ}
\pgfmathsetmacro{\fK}{0.15625+\bK}
\pgfmathsetmacro{\fL}{0.171875+\bL}
\pgfmathsetmacro{\fM}{0.1875+\bM}
\pgfmathsetmacro{\fN}{0.203125+\bN}
\pgfmathsetmacro{\fO}{0.21875+\bO}
\pgfmathsetmacro{\fP}{0.234375+\bP}
\pgfmathsetmacro{\fQ}{0.25+\bQ}

\pgfmathsetmacro{\gA}{\uA}
\pgfmathsetmacro{\gB}{\uB}
\pgfmathsetmacro{\gC}{\uC}
\pgfmathsetmacro{\gD}{\uD}
\pgfmathsetmacro{\gE}{\uE}
\pgfmathsetmacro{\gF}{\uF}
\pgfmathsetmacro{\gG}{\uG}
\pgfmathsetmacro{\gH}{\uH}
\pgfmathsetmacro{\gI}{\uI}
\pgfmathsetmacro{\gJ}{0.03125+\uJ}
\pgfmathsetmacro{\gK}{0.0625+\uK}
\pgfmathsetmacro{\gL}{0.09375+\uL}
\pgfmathsetmacro{\gM}{0.125+\uM}
\pgfmathsetmacro{\gN}{0.15625+\uN}
\pgfmathsetmacro{\gO}{0.1875+\uO}
\pgfmathsetmacro{\gP}{0.21875+\uP}
\pgfmathsetmacro{\gQ}{0.25+\uQ}

\pgfmathsetmacro{\hA}{0.5*\bA}
\pgfmathsetmacro{\hB}{0.0078125+0.5*\bB}
\pgfmathsetmacro{\hC}{0.015625+0.5*\bC}
\pgfmathsetmacro{\hD}{0.0234375+0.5*\bD}
\pgfmathsetmacro{\hE}{0.03125+0.5*\bE}
\pgfmathsetmacro{\hF}{0.0390625+0.5*\bF}
\pgfmathsetmacro{\hG}{0.046875+0.5*\bG}
\pgfmathsetmacro{\hH}{0.0546875+0.5*\bH}
\pgfmathsetmacro{\hI}{0.0625+0.5*\bI}
\pgfmathsetmacro{\hJ}{0.0703125+0.5*\bJ}
\pgfmathsetmacro{\hK}{0.078125+0.5*\bK}
\pgfmathsetmacro{\hL}{0.0859375+0.5*\bL}
\pgfmathsetmacro{\hM}{0.09375+0.5*\bM}
\pgfmathsetmacro{\hN}{0.1015625+0.5*\bN}
\pgfmathsetmacro{\hO}{0.109375+0.5*\bO}
\pgfmathsetmacro{\hP}{0.1171875+0.5*\bP}
\pgfmathsetmacro{\hQ}{0.125+0.5*\bQ}

\pgfmathsetmacro{\jB}{\hQ+0.5*\gB}
\pgfmathsetmacro{\jC}{\hQ+0.5*\gC}
\pgfmathsetmacro{\jD}{\hQ+0.5*\gD}
\pgfmathsetmacro{\jE}{\hQ+0.5*\gE}
\pgfmathsetmacro{\jF}{\hQ+0.5*\gF}
\pgfmathsetmacro{\jG}{\hQ+0.5*\gG}
\pgfmathsetmacro{\jH}{\hQ+0.5*\gH}
\pgfmathsetmacro{\jI}{\hQ+0.5*\gI}
\pgfmathsetmacro{\jJ}{\hQ+0.5*\gJ}
\pgfmathsetmacro{\jK}{\hQ+0.5*\gK}
\pgfmathsetmacro{\jL}{\hQ+0.5*\gL}
\pgfmathsetmacro{\jM}{\hQ+0.5*\gM}
\pgfmathsetmacro{\jN}{\hQ+0.5*\gN}
\pgfmathsetmacro{\jO}{\hQ+0.5*\gO}
\pgfmathsetmacro{\jP}{\hQ+0.5*\gP}
\pgfmathsetmacro{\jQ}{\hQ+0.5*\gQ}

\begin{scope}
    \draw[axis] (-0.07,0) -- (1.08,0) node[below right=1pt] {$t$};
    \draw[axis] (0,-0.025) -- (0,0.31);
    \draw (0,0.25) -- (-0.02,0.25) node[left=2pt] {$\frac14$};
    \draw (0.5,0) -- (0.5,-0.012) node[below=2pt] {$\frac12$};
    \draw (1,0) -- (1,-0.012) node[below=2pt] {$1$};
    \node[below left=1pt] at (0,0) {$0$};

    \draw[original] (0,0) -- (1,0.25);

    \draw[approx]
        (0,\fA)
        -- (0.0625,\fB) -- (0.125,\fC)
        -- (0.1875,\fD) -- (0.25,\fE)
        -- (0.3125,\fF) -- (0.375,\fG)
        -- (0.4375,\fH) -- (0.5,\fI)
        -- (0.5625,\fJ) -- (0.625,\fK)
        -- (0.6875,\fL) -- (0.75,\fM)
        -- (0.8125,\fN) -- (0.875,\fO)
        -- (0.9375,\fP) -- (1,\fQ);

    \draw[original] (0.07,-0.08) -- (0.20,-0.08);
    \node[anchor=west,font=\footnotesize] at (0.23,-0.08)
        {$H^{\text{alt}}$};

    \draw[approx] (0.07,-0.125) -- (0.20,-0.125);
    \node[anchor=west,font=\footnotesize] at (0.23,-0.125)
        {$H_n(w^{\text{alt}},\cdot)$};
\end{scope}

\begin{scope}[xshift=4.4cm]
    \draw[axis] (-0.07,0) -- (1.08,0) node[below right=1pt] {$t$};
    \draw[axis] (0,-0.025) -- (0,0.31);
    \draw (0,0.25) -- (-0.02,0.25) node[left=2pt] {$\frac14$};
    \draw (0.5,0) -- (0.5,-0.012) node[below=2pt] {$\frac12$};
    \draw (1,0) -- (1,-0.012) node[below=2pt] {$1$};
    \node[below left=1pt] at (0,0) {$0$};

    \draw[original] (0,0) -- (0.5,0) -- (1,0.25);

    \draw[approx]
        (0,\gA)
        -- (0.0625,\gB) -- (0.125,\gC)
        -- (0.1875,\gD) -- (0.25,\gE)
        -- (0.3125,\gF) -- (0.375,\gG)
        -- (0.4375,\gH) -- (0.5,\gI)
        -- (0.5625,\gJ) -- (0.625,\gK)
        -- (0.6875,\gL) -- (0.75,\gM)
        -- (0.8125,\gN) -- (0.875,\gO)
        -- (0.9375,\gP) -- (1,\gQ);

    \draw[original] (0.07,-0.08) -- (0.20,-0.08);
    \node[anchor=west,font=\footnotesize] at (0.23,-0.08)
        {$H^{\text{seq}}$};

    \draw[approx] (0.07,-0.125) -- (0.20,-0.125);
    \node[anchor=west,font=\footnotesize] at (0.23,-0.125)
        {$H_n(w^{\text{seq}},\cdot)$};
\end{scope}

\begin{scope}[xshift=8.8cm]
    \draw[axis] (-0.07,0) -- (1.08,0) node[below right=1pt] {$t$};
    \draw[axis] (0,-0.025) -- (0,0.31);
    \draw (0,0.125) -- (-0.02,0.125) node[left=2pt] {$\frac18$};
    \draw (0,0.25) -- (-0.02,0.25) node[left=2pt] {$\frac14$};
    \draw (0.5,0) -- (0.5,-0.012) node[below=2pt] {$\frac12$};
    \draw (0.75,0) -- (0.75,-0.012) node[below=2pt] {$\frac34$};
    \draw (1,0) -- (1,-0.012) node[below=2pt] {$1$};
    \node[below left=1pt] at (0,0) {$0$};

    \draw[original]
        (0,0) -- (0.5,0.125) -- (0.75,0.125) -- (1,0.25);

    \draw[approx]
        (0,\hA)
        -- (0.03125,\hB) -- (0.0625,\hC)
        -- (0.09375,\hD) -- (0.125,\hE)
        -- (0.15625,\hF) -- (0.1875,\hG)
        -- (0.21875,\hH) -- (0.25,\hI)
        -- (0.28125,\hJ) -- (0.3125,\hK)
        -- (0.34375,\hL) -- (0.375,\hM)
        -- (0.40625,\hN) -- (0.4375,\hO)
        -- (0.46875,\hP) -- (0.5,\hQ)
        -- (0.53125,\jB) -- (0.5625,\jC)
        -- (0.59375,\jD) -- (0.625,\jE)
        -- (0.65625,\jF) -- (0.6875,\jG)
        -- (0.71875,\jH) -- (0.75,\jI)
        -- (0.78125,\jJ) -- (0.8125,\jK)
        -- (0.84375,\jL) -- (0.875,\jM)
        -- (0.90625,\jN) -- (0.9375,\jO)
        -- (0.96875,\jP) -- (1,\jQ);

    \draw[average]
        (0,0) -- (0.5,0.0625) -- (1,0.25);

    
    \draw[original] (0.07,-0.08) -- (0.20,-0.08);
    \node[anchor=west,font=\footnotesize] at (0.23,-0.08)
        {$H^{\text{cat}}$};

    \draw[approx] (0.07,-0.125) -- (0.20,-0.125);
    \node[anchor=west,font=\footnotesize] at (0.23,-0.125)
        {$H_{2n}(\overline{w},\cdot)$};

    \draw[average] (0.07,-0.17) -- (0.20,-0.17);
    \node[anchor=west,font=\footnotesize] at (0.23,-0.17)
        {$\frac{H^{\text{alt}}+H^{\text{seq}}}{2}$};

\end{scope}

\end{tikzpicture}
\captionsetup{justification=centering}
\caption{
Schematic evolution of the $e_2$-coordinate of the two target
profiles, their concatenation, and their arithmetic mean.
}
\label{figure_evolution_e_2_coordinate}
\end{figure}

In Example~\ref{example_balanced_unbalanced_profile},  we constructed two trajectories with the same empirical measure $\widehat{\mu}=1/2 \delta_0+ 1/4 \delta_{e_1} + 1/4 \delta_{e_2}$, but with different limiting profiles $H^{\text{alt}}$ and $H^{\text{seq}}$. Consider two trajectories 
\begin{align*}
    w^{\text{alt}} \in \Omega_n \big(\widehat{\mu},H^{\text{alt}},\varepsilon \big) \qquad \text{and} \qquad w^{\text{seq}} \in \Omega_n \big(\widehat{\mu},H^{\text{seq}},\varepsilon \big).
\end{align*}
When their endpoints are compatible, we can concatenate them to obtain $\overline{w}=w^{\text{alt}}w^{\text{seq}} \in \Omega_{2n}$. At the macroscopic time scale, the profile of the concatenation is naturally described by rescaling and translating the profiles of the constituent trajectories. This leads us to consider the function
\begin{equation}\label{eq_concatenated_profile}
H^{\text{cat}}(t)=
\begin{cases}
\displaystyle \frac12 H^{\text{alt}}(2t),
& t\in[0,\frac12],\\[0.6em]
\displaystyle \frac12 H^{\text{alt}}(1)
+\frac12 H^{\text{seq}}(2t-1),
& t\in[\frac12,1].
\end{cases}
\end{equation}
Since the profiles of the constituent trajectories are $\varepsilon$-close to $H^{\text{alt}}$ and $H^{\text{seq}}$, we expect $H_{2n}(\overline{w},\cdot)$ to be uniformly close to $H^{\text{cat}}$. This target profile is generally different from the arithmetic mean $(H^{\text{alt}}+H^{\text{seq}})/2$. This is illustrated in Figure~\ref{figure_evolution_e_2_coordinate}.

\noindent This example shows that concatenation need not preserve the arithmetic mean of the profiles: unlike the empirical measure, which is insensitive to the order of visits, the profile records the temporal ordering of increments. Thus, the stitching construction of \cite{Daures}, which controls only the empirical measure, does not directly extend to our setting, where the profile must also be controlled.

\subsubsection{Choice of parameters}\label{subsubsection_paramters_for_construction}

Before formally defining the stitching procedure $\psi$, we introduce the parameters and scaling choices required for the construction. The order in which the parameters are chosen is important, and their roles are summarised in Table~\ref{table_order_parameter_choices}.

For $m \in \NN$, define the regular mesh of the segment $[0,1]$ by 
\begin{equation}\label{eq_def_overline_x_m}
\overline{x}(m) \eqdef  \left\{ x_h = \frac{h-1}{m} \, : \, h \in \dbl 1, m+1 \dbr \right\}.
\end{equation}

Fix $\varepsilon>0$ and $\kappa>0$. We first choose $m \in \NN$ sufficiently large so that 
\begin{equation}\label{eq_cdt_m_sufficiently_large}
\frac{1}{m}<\frac{\varepsilon}{2}.
\end{equation} 
For a word $w \in \Omega_n^{(0)}$, we decompose $w$ according to $\overline{x}(m)$ as follows. For $h \in \dbl 1,m \dbr$, set
\begin{equation}\label{eq_def_subwords_for_profile_partition}
u^h(w)=u^h \eqdef w_{[ \lfloor x_h n \rfloor +1 : \lfloor x_{h+1} n \rfloor ]}.
\end{equation}
Then $w=u^1 u^2 \cdots u^{m}$ and
\begin{align}
 (x_{h+1}-x_h)n-1 \leq |u^{h}| & \leq (x_{h+1}-x_h)n+1, \qquad \quad \label{eq_bound_length_u_h_with_x_h}\\
\sum_{l=1}^{h} |u^l| &= \lfloor x_{h+1}n \rfloor. \label{eq_bound_length_sum_u_h_with_x_h}
\end{align}
In what follows, we will impose conditions on $n$ which ensure that each block $u^h$ is non-empty, see the condition given in \eqref{eq_condition_n}. We define 
\begin{equation}\label{eq_def_R_varepsilon}
R(\varepsilon) \eqdef \max \left\{ R_{\mu_1}(\varepsilon) \, , \, R_{\mu_2}(\varepsilon) \right\},
\end{equation}
where $R_{\mu_i}(\varepsilon)$ is given by Lemma~\ref{lemma_bound_points_outside_compact}.

\begin{lemma}[Uniformly short connectors]\label{lemma_existence_of_xi_with_right_property}
Let $n > R(\varepsilon) / \varepsilon$. For all $w,w' \in \Omega_n(\mu_1,\varepsilon) \cup \Omega_n(\mu_2,\varepsilon)$ and $i,j \in \dbl 1,n \dbr$, there exists a word $\xi \in \Omega_{\textup{fin}}$ such that 
\begin{enumerate}
\item $w_i \sim \xi_1$ and $\xi_{-1} \sim w_j'$,
\item $|\xi| < 6d\varepsilon n$,
\end{enumerate}
with the convention that the first item is satisfied by $\xi=\varnothing$ when $w_i \sim w_j'$.
\end{lemma}

\begin{proof}
Fix $w,w' \in \Omega_n(\mu_1,\varepsilon) \cup \Omega_n(\mu_2,\varepsilon)$ and $i,j \in \dbl 1,n \dbr$. Since
$\ell_n(w),\ell_n(w') \in B(\mu_1,\varepsilon) \cup B(\mu_2,\varepsilon)$, Corollary~\ref{corollary_bound_on_last_point_of_subword}, the definition of $R(\varepsilon)$ and the fact that $n>R(\varepsilon)/\varepsilon$ yield
\begin{align*}
\|w_i\|_{\infty} < 3 \varepsilon n
\qquad \text{and} \qquad
\|w'_j\|_{\infty} < 3 \varepsilon n.
\end{align*}
It follows that the $L^1$-distance between $w_i$ and $w'_j$ satisfies
\begin{align*}
\|w_i-w'_j\|_1 < 6d\varepsilon n.
\end{align*}
Hence, one can construct a nearest-neighbour path from a site adjacent to $w_i$
to a site adjacent to $w'_j$ with at most $6d\varepsilon n$ steps. This defines
a word $\xi \in \Omega_{\textup{fin}}$ such that
$w_i \sim \xi_1$, $\xi_{-1} \sim w'_j$, and
$|\xi| < 6d\varepsilon n$. In the case where $w_i \sim w'_j$, one can simply take $\xi$ to be the empty word.
\end{proof}

\medskip

We then set
\begin{equation}\label{eq_def_C_d_m}
C_{d,m} \eqdef 6d (m+1) \in \NN.
\end{equation}
\textbf{Condition on $\varepsilon'$} : Fix $\varepsilon>0$, $\kappa>0$ and $m \in \NN$ satisfying \eqref{eq_cdt_m_sufficiently_large}. Consider $\varepsilon_0'(\varepsilon,\kappa,m)=\varepsilon_0'$ sufficiently small such that
\begin{align}
6d\varepsilon_0' & < 1,  \label{eq_condition_varepsilon_eta_2} \\
\frac{1}{m} + \left(20C_{d,m}+m+20 \right) \varepsilon_0' & < \varepsilon, \label{eq_condition_varepsilon_eta_3}\\
(2C_{d,m}+2) \log(2d) \, \varepsilon_0'& <\kappa. \label{eq_condition_varepsilon_eta_4} 
\end{align}
Because of \eqref{eq_cdt_m_sufficiently_large}, such an $\varepsilon'_0$ exists. The constants are not intended to be sharp. The conditions simply require $\varepsilon'_0$ to be sufficiently small
once $\varepsilon$, $\kappa$ and $m$ have been fixed. Note that Equations \eqref{eq_condition_varepsilon_eta_2}--\eqref{eq_condition_varepsilon_eta_4} hold for all $\varepsilon'<\varepsilon_0'$. We now introduce an additional set of conditions on $n$ and $N$.

\begin{hyp}\label{condition_n_N}
Fix $\varepsilon'>0$ and $m \in \NN$. We say that $(n,N) \in \NN^2$ satisfies Condition~\ref{condition_n_N} if 
\begin{align}
n & > \max \Bigg\{ \frac{4}{\varepsilon'} , \frac{R(\varepsilon')}{\varepsilon'}, 4m\Bigg\}, \label{eq_condition_n}\\
\text{and} \qquad N & > 2n(1+C_{d,m}\varepsilon')\left(\frac{4}{\varepsilon'}+1\right). \label{eq_condition_N}
\end{align}
\end{hyp}

\noindent As before, the exact form of these conditions is not important. Some of them are redundant, but we keep them for the sake of clarity. What matters is that $N$ and $n$ are sufficiently large. For every $\varepsilon'>0$ and $m \in \NN$, the condition is satisfied by first taking $n$ and then $N$ sufficiently large.

For each $m \in \NN$, $\varepsilon'>0$, $n \in \NN$ and $N \in \NN$, define
\begin{equation}\label{eq_def_of_k}
k(m,\varepsilon',n,N)=k \eqdef \left \lfloor \frac{N}{2n(1+C_{d,m}\varepsilon')} \right \rfloor. 
\end{equation}

\begin{lemma}[Block-count estimates]\label{lemma_N_large_enough}
If the pair $(n,N)$ satisfies Condition~\ref{condition_n_N}, then
\begin{align}
\frac{4}{\varepsilon'} & <k, \label{eq_lemma_N_large_enough} \\
\text{and} \qquad \left| 1 - \frac{2kn}{N} \right| & < \big(2C_{d,m}+1 \big)\varepsilon'. \label{eq_bound_on_distance_n_N_k}
\end{align}
\end{lemma}

\begin{proof}
Using the definition of $k$ and \eqref{eq_condition_N}, we deduce that
\begin{align*}
k+1 \geq \frac{N}{2n(1+C_{d,m}\varepsilon')}>\frac{4}{\varepsilon'}+1,
\end{align*}
which proves \eqref{eq_lemma_N_large_enough}. Moreover, the definition of $k$ gives
\begin{align*}
N-2nk \leq 2knC_{d,m} \varepsilon' + 2n(1+C_{d,m} \varepsilon').
\end{align*}
Since $2kn \leq N$, we obtain
\begin{align*}
0 \leq 1-\frac{2nk}{N} \leq C_{d,m} \varepsilon'  + \frac{1+C_{d,m} \varepsilon'}{k}.
\end{align*}
By \eqref{eq_lemma_N_large_enough}, we know that $k \geq 1$ and so
\begin{align*}
\left| 1 - \frac{2kn}{N} \right| & \leq 2 C_{d,m} \varepsilon' + \frac{1}{k},
\end{align*}
and Equation~\eqref{eq_bound_on_distance_n_N_k} follows by applying \eqref{eq_lemma_N_large_enough} again, which gives $1/k < \varepsilon'/4<\varepsilon'$.
\end{proof}

\noindent We summarise the dependence of the parameters in Table~\ref{table_order_parameter_choices}.

\begin{table}[htbp]
\centering
\caption{Order of parameter choices and their role in the stitching argument.}
\label{table_order_parameter_choices}
\renewcommand{\arraystretch}{1.25}
\begin{tabular}{|c|p{0.13\textwidth}|p{0.72\textwidth}|}
\hline
Step & Parameters & Role of the choice \\
\hline
1
& $\varepsilon>0$
& Determines the target neighbourhood of the final event
\newline
$\displaystyle
B\!\left(\frac{\mu_1+\mu_2}{2},\varepsilon\right)
\times
B_{\mathscr{H}}\!\left(
\frac{H^{(1)}+H^{(2)}}{2},\varepsilon
\right).$
\arraybackslash
\\
\hline
1'
& $\kappa>0$
& Fixes the auxiliary tolerance for the exponential error term. \\
\hline
2
& $m \in \NN$
& Fixes the size of the mesh on $[0,1]$.\\
\hline
3
& $0<\varepsilon'<\varepsilon_0'$
& Determines the smaller input neighbourhood
$B(\mu_i,\varepsilon')\times
B_{\mathscr{H}}(H^{(i)},\varepsilon')$. \\
\hline
4
& $n\in\mathbb{N}$
& Length of the smaller input words. Chosen large enough to control boundary effects and ensure that every block in the decomposition is non-empty. \\
\hline
5
& $N\in\mathbb{N}$
& Length of the output of the stitching procedure. Determines the number $k=k(m,\varepsilon',n,N)$ of pairs of input trajectories, and hence
the total number $2k$ of input trajectories to be stitched. \\
\hline
\end{tabular}
\end{table}

\medskip

\noindent Note that the parameter $\kappa$ is chosen at the same time as the parameter $\varepsilon$. Its role is to measure how much loss we allow in the logarithmic probability estimates when proving that the stitching procedure satisfies the required lower bound.

\subsubsection{Definition of the stitching map}\label{subsubsection_constructing_the_stitching_map}

Set
\begin{equation}\label{eq_def_W_n_k}
\mathcal{W}_{n,k} \eqdef \Big( \Omega_n\big(\mu_1,H^{(1)},\varepsilon'\big)\times \Omega_n\big(\mu_2,H^{(2)},\varepsilon'\big) \Big)^k.
\end{equation}
We identify each element of this set with the corresponding $2k$-tuple 
\begin{equation}\label{eq_convention_W_n_k_couples}
\begin{aligned}
\big((w^{(1)},w^{(2)}),\ldots, (w^{(2k-1)},w^{(2k)})\big)\longleftrightarrow (w^{(1)},\ldots,w^{(2k)}).
\end{aligned}
\end{equation}
Thus, whenever $(w^{(1)},\ldots,w^{(2k)})\in\mathcal{W}_{n,k}$, it is understood that, for every $i\in\dbl 1,k\dbr$,
\begin{equation*}
w^{(2i-1)} \in \Omega_n\big(\mu_1,H^{(1)},\varepsilon'\big) \qquad \text{and} \qquad w^{(2i)} \in \Omega_n\big(\mu_2,H^{(2)},\varepsilon'\big).
\end{equation*}

The stitching map $\psi$ introduced below is the main ingredient of the decoupling argument. It constructs a trajectory of length $N$ from $2k$ trajectories of length $n$. For each $i \in \dbl 1,k\dbr$, one trajectory is localised near $(\mu_1,H^{(1)})$ and the other near $(\mu_2,H^{(2)})$. 

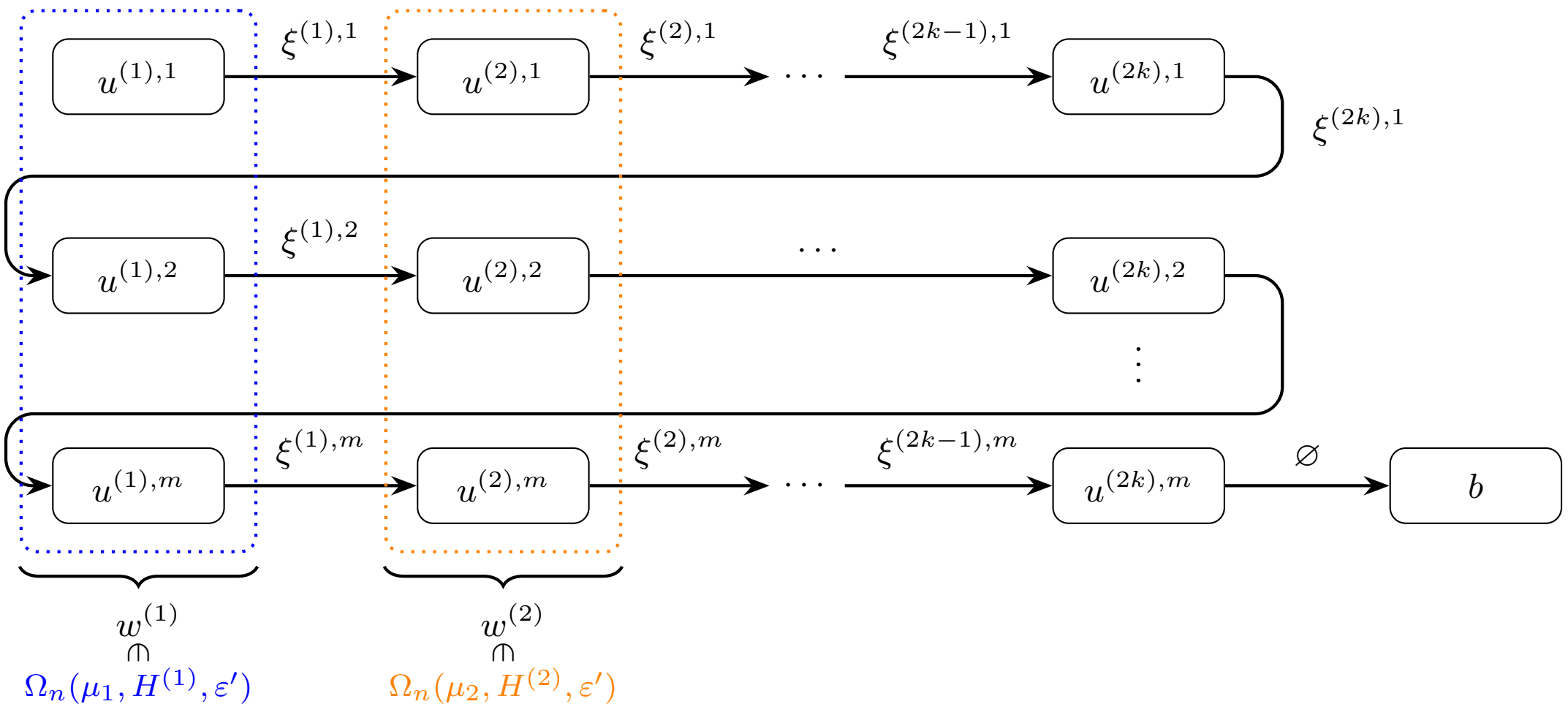
\begin{figure}[htbp]
\centering
\begin{tikzpicture}[
    x=1cm,y=1cm,
    every node/.style={font=\small},
    block/.style={
        draw=black,
        rounded corners,
        minimum width=1.55cm,
        minimum height=.68cm
    },
    conn/.style={-{Stealth},draw=black,thick},
    lab/.style={font=\small,inner sep=1pt}
]

\node[block] (u11) at (0,0) {$u^{(1),1}$};
\node[block] (u21) at (3.30,0) {$u^{(2),1}$};
\node        (d1)  at (6.05,0) {$\cdots$};
\node[block] (uk1) at (9.05,0) {$u^{(2k),1}$};

\node[block] (u12) at (0,-1.8) {$u^{(1),2}$};
\node[block] (u22) at (3.30,-1.8) {$u^{(2),2}$};
\node[block] (uk2) at (9.05,-1.8) {$u^{(2k),2}$};

\node at (9.05,-2.50) {$\vdots$};

\node[block] (u1m) at (0,-3.7) {$u^{(1),m}$};
\node[block] (u2m) at (3.30,-3.7) {$u^{(2),m}$};
\node        (dm)  at (6.05,-3.7) {$\cdots$};
\node[block] (ukm) at (9.05,-3.7) {$u^{(2k),m}$};
\node[block] (b)   at (12.1,-3.7) {$b$};

\draw[conn]
    (u11) -- node[lab,above=3pt] {$\xi^{(1),1}$} (u21);

\draw[conn]
    (u21) -- node[lab,above=3pt] {$\xi^{(2),1}$} (d1);

\draw[conn]
    (d1) -- node[lab,above=3pt] {$\xi^{(2k-1),1}$} (uk1);

\draw[conn]
    (u12) -- node[lab,above=3pt] {$\xi^{(1),2}$} (u22);

\draw[conn]
    (u22) -- node[lab,above=3pt] {$\cdots$} (uk2);

\draw[conn]
    (u1m) -- node[lab,above=3pt] {$\xi^{(1),m}$} (u2m);

\draw[conn]
    (u2m) -- node[lab,above=3pt] {$\xi^{(2),m}$} (dm);

\draw[conn]
    (dm) -- node[lab,above=3pt] {$\xi^{(2k-1),m}$} (ukm);

\draw[conn]
    (ukm) -- node[lab,above=3pt] {$\varnothing$} (b);

\draw[conn,rounded corners=7pt]
    (uk1.east)
    -- (10.35,0)
    -- node[midway,right=4pt] {$\xi^{(2k),1}$} (10.35,-0.90)
    -- (-1.20,-0.90)
    -- (-1.20,-1.80)
    -- (u12.west);

\draw[conn,rounded corners=7pt]
    (uk2.east)
    -- (10.35,-1.80)
    -- (10.35,-3.05)
    -- (-1.20,-3.05)
    -- (-1.20,-3.70)
    -- (u1m.west);

\node[
    draw=blue,
    dotted,
    line width=0.8pt,
    rounded corners,
    fit=(u11)(u12)(u1m),
    inner xsep=0.28cm,
    inner ysep=0.25cm
] (w1box) {};

\node[
    draw=orange,
    dotted,
    line width=0.8pt,
    rounded corners,
    fit=(u21)(u22)(u2m),
    inner xsep=0.28cm,
    inner ysep=0.25cm
] (w2box) {};

\draw[
    decorate,
    decoration={brace,mirror,amplitude=5pt},
    thick
]
    ($(w1box.south west)+(0,-0.12)$) --
    ($(w1box.south east)+(0,-0.12)$);

\draw[
    decorate,
    decoration={brace,mirror,amplitude=5pt},
    thick
]
    ($(w2box.south west)+(0,-0.12)$) --
    ($(w2box.south east)+(0,-0.12)$);

\node at ($(w1box.south)+(0,-0.60)$)
    {$\; \; w^{(1)}$};

\node[rotate=270] at ($(w1box.south)+(0,-0.90)$)
    {$\in$};

\node[font=\footnotesize] at ($(w1box.south)+(0,-1.2)$)
    {\textcolor{blue}{$\Omega_n(\mu_1,H^{(1)},\varepsilon')$}};

\node at ($(w2box.south)+(0,-0.60)$)
    {$\; \; w^{(2)}$};

\node[rotate=270] at ($(w2box.south)+(0,-0.90)$)
    {$\in$};

\node[font=\footnotesize] at ($(w2box.south)+(0,-1.2)$)
    {\textcolor{orange}{$\Omega_n(\mu_2,H^{(2)},\varepsilon')$}};

\end{tikzpicture}

\caption{Schematic representation of the stitching procedure $\psi$.}
\label{fig_stitching_procedure}
\end{figure}

\begin{definition}[Stitching procedure]\label{def_psi_stitching_map}
Consider $m \in \NN$ and $0<\varepsilon'<\varepsilon_0'$, and suppose that $(n,N) \in \NN^2$ satisfies Condition~\ref{condition_n_N}. For $(w^{(1)},\dots,w^{(2k)}) \in \mathcal{W}_{n,k}$, decompose each input word according to \eqref{eq_def_subwords_for_profile_partition}, writing $u^{(i),h}=u^h(w^{(i)})$. We define the stitching map
\begin{align*}
\psi : \mathcal{W}_{n,k} \longrightarrow  \Omega_N,
\end{align*}
as follows. For every $(w^{(1)},\dots,w^{(2k)}) \in \mathcal{W}_{n,k}$, set	 
\begin{equation}\label{eq_def_psi_stitching_map}
\begin{array}{cccc}
\psi \left(w^{(1)},\dots,w^{(2k)}\right) = u^{(1),1} \xi^{(1),1} u^{(2),1} \cdots u^{(2k),1} \xi^{(2k),1} u^{(1),2} \cdots u^{(2k),m} b, 
\end{array}
\end{equation}
where the components of the image are constructed as follows.
\begin{enumerate}
\item \textit{Transitions ($\xi$):} For each $(i,h) \neq (2k,m)$, $\xi^{(i),h}$ is a transition word of length $|\xi^{(i),h}| < 6d \varepsilon' n$ ensuring compatibility between the subwords:
\begin{equation*}
u^{(i),h}_{-1} \sim \xi^{(i),h}_1 \quad \text{and} \quad \xi^{(i),h}_{-1} \sim u^{(i+1),h}_1,
\end{equation*}
where, for every $h \in \dbl 1,m-1 \dbr$, we use the cyclic convention  $u^{(2k+1),h} \eqdef u^{(1),h+1}$. For notational convenience, we set $\xi^{(2k),m}=\varnothing$ and insert this empty word between $u^{(2k),m}$ and $b$.
    
\item \textit{Padding ($b$):} The word $b \in \Omega_{\textup{fin}}$ pads the sequence to total length $N$, and satisfies the final compatibility condition $u^{(2k),m}_{-1} \sim b_1$.
\end{enumerate}
Whenever several connecting words satisfy the required conditions, we choose the first one with respect to a fixed total order on $\Omega_{\textup{fin}}$. With this convention, $\psi$ is a well-defined map. It is schematically represented in Figure~\ref{fig_stitching_procedure}.
\end{definition}

\noindent Note that the map $\psi$ also depends on the parameters $n,N,m,\mu_1,\mu_2,H^{(1)},H^{(2)}$ and $\varepsilon'$. By the conventions introduced in Section~\ref{subsection_space_of_trajectories}, the empty word $\xi^{(2k),m}=\varnothing$ has length zero, contributes no increments, and leaves concatenations unchanged. In particular, in all sums involving connecting words, the term corresponding to $(i,h)=(2k,m)$ contributes zero.

\begin{remark}\label{remark_psi_well_defined}
The map $\psi$ is well defined. Indeed, \eqref{eq_condition_n} permits the use of Lemma~\ref{lemma_existence_of_xi_with_right_property}. Hence one can always find transition words $\xi^{(i),h}$ which satisfy the desired properties. Moreover,
\begin{align*}
\sum_{i=1}^{2k} \sum_{h=1}^{m}(|u^{(i),h}| + |\xi^{(i),h}|) & < \sum_{i=1}^{2k} |w^{(i)}| + \sum_{i=1}^{2k} \sum_{h=1}^{m} 6d \varepsilon' n \\
& = 2kn + 2km 6d \varepsilon' n \\
& \leq N,
\end{align*}
where the last inequality comes from the definition of $C_{d,m}$ and \eqref{eq_def_of_k}. Hence a non-empty padding word $b$ of the remaining length can be appended to the last subword.
\end{remark}

\subsection{Properties of the stitching map}\label{subsection_verifying_prop_of_stitching}

In this section, we show that the stitching procedure defined above satisfies the properties \eqref{eq_objectif_inclusion}, \eqref{eq_objectif_cardinality_bound}, \eqref{eq_objectif_comp_proba} and \eqref{eq_objectif_asympto} required by Proposition~\ref{prop_good_stitching_procedure}.

\subsubsection{Target-set inclusion}\label{subsubsection_inclusion_property}

Throughout this section, we fix $\varepsilon>0$, $\kappa>0$ and $m \in \NN$ satisfying \eqref{eq_cdt_m_sufficiently_large}. Our aim is to show that $\psi$ satisfies the inclusion property stated in \eqref{eq_objectif_inclusion}. 

The following lemma shows that the empirical measure of the trajectory produced by the stitching procedure is close to the average empirical measure of its constituent trajectories.

\begin{lemma}[Empirical measure stability under stitching]\label{lemma_property_stitching_map_measure_close}
Suppose that $\varepsilon'<\varepsilon_0'$ and that $n \in \NN$ and $N \in \NN$ satisfy Condition~\ref{condition_n_N}. Then, for all $(w^{(1)},\dots,w^{(2k)}) \in \mathcal{W}_{n,k}$, 
\begin{equation*}
\ell_N\big( \psi(w^{(1)},\dots,w^{(2k)}) \big) \in  B\left( \frac{\mu_1+\mu_2}{2}, \varepsilon \right).
\end{equation*}
\end{lemma}

\begin{proof}
Fix $(w^{(1)},\dots,w^{(2k)}) \in \mathcal{W}_{n,k}$ and let $\overline{w}=\psi\big(w^{(1)},\dots,w^{(2k)}\big)$. By decomposing $\overline{w}$ as in \eqref{eq_def_psi_stitching_map} and using the properties of the KR norm given in \eqref{eq_lemma_properties_KR_norm_of_proba}, we have
\begin{align*}
\left\| \ell_N(\overline{w}) - \frac{\mu_1+\mu_2}{2} \right\| & = \left\| \frac{1}{N} \sum_{i=1}^{2k} \sum_{h=1}^{ m} \sum_{j=1}^{|u^{(i),h}|} \delta_{u^{(i),h}_j} + \frac{1}{N} \sum_{i=1}^{2k} \sum_{h=1}^{m} \sum_{j=1}^{|\xi^{(i),h}|} \delta_{\xi^{(i),h}_j} + \frac{1}{N} \sum_{j=1}^{|b|} \delta_{b_j} -\frac{\mu_1+\mu_2}{2} \right\| \\
&  \leq \left| \frac{2nk}{N} - 1 \right| + \frac{1}{2k} \sum_{i=1}^{k} \left\| \frac{1}{n} \sum_{j=1}^{n} \delta_{w^{(2i-1)}_j} - \mu_1 \right\| + \frac{1}{2k} \sum_{i=1}^{k} \left\| \frac{1}{n} \sum_{j=1}^{n} \delta_{w^{(2i)}_j} - \mu_2 \right\| \\
& \quad + \frac{1}{N} \left(\sum_{i=1}^{2k} \sum_{h=1}^{m}|\xi^{(i),h}|+|b| \right).
\end{align*}
Since $\sum_{i=1}^{2k} \sum_{h=1}^{m}|\xi^{(i),h}|+|b|=N-2nk$, we deduce that
\begin{align*}
\left\| \ell_N(\overline{w}) - \frac{\mu_1+\mu_2}{2} \right\|  \leq 2 \left| \frac{2nk}{N} - 1 \right| + \frac{1}{2k} \sum_{i=1}^k \left\| \ell_n \big(w^{(2i-1)}\big) - \mu_1 \right\|+\frac{1}{2k} \sum_{i=1}^k \left\| \ell_n \big(w^{(2i)}\big) - \mu_2 \right\|.
\end{align*}
By definition of $\mathcal{W}_{n,k}$, given in \eqref{eq_def_W_n_k}, we have $\ell_n(w^{(2i-1)}) \in B(\mu_1,\varepsilon')$ and $\ell_n(w^{(2i)}) \in B(\mu_2,\varepsilon')$, for all $i \in \dbl 1,k \dbr$. Combining this with \eqref{eq_bound_on_distance_n_N_k}, we obtain 
\begin{align*}
\left\| \ell_N(\overline{w}) - \frac{\mu_1+\mu_2}{2} \right\| < 2 \big(2C_{d,m}+1\big)\varepsilon' + \varepsilon'.
\end{align*}
Since $4C_{d,m}+3<20C_{d,m}+m+20$ and $\varepsilon'<\varepsilon_0'$, Equation~\eqref{eq_condition_varepsilon_eta_3} implies $(4C_{d,m}+3)\varepsilon'<\varepsilon$. The desired result follows.
\end{proof}

Before showing that the stitching procedure also preserves the profile, we first state a technical lemma. Fix $(w^{(1)},\dots,w^{(2k)}) \in \mathcal{W}_{n,k}$. For $h \in \dbl 2,m+1 \dbr$, set
\begin{equation}\label{eq_def_N_h_and_n_h}
n_h \eqdef \sum_{r=1}^{h-1}|u^{(1),r}| \qquad \text{and} \qquad N_h \eqdef \sum_{i=1}^{2k} \sum_{r=1}^{h-1}(|u^{(i),r}| + |\xi^{(i),r}|).
\end{equation}
Note that for all $h \in \dbl 2, m \dbr$, $N_h+1$ is the index of $u^{(1),h}_1$ in $\psi(w^{(1)},\dots,w^{(2k)})$ and $n_h+1$ is the index of $u^{(i),h}_1$ in $w^{(i)}$ for all $i \in \dbl 1,2k \dbr$. We also define the corresponding profile times
\begin{equation}
\sigma_h \eqdef \frac{n_h-1}{n-1} \qquad \text{and} \qquad \tau_h = \frac{N_h-1}{N-1},
\end{equation}
and for $h=1$, we set $\sigma_1=\tau_1=0$.

\begin{lemma}[Location of stitched block boundaries]\label{lemma_sigma_h_and_tau_h}
Suppose that $0<\varepsilon'<\varepsilon_0'$ and that $n \in \NN$ and $N \in \NN$ satisfy Condition~\ref{condition_n_N}. Then, for all $h \in \dbl 1, m+1 \dbr$,
\begin{equation}\label{eq_lemma_sigma_h_and_tau_h_1}
\left \vert \sigma_h-x_h \right \vert <  \varepsilon' \qquad \text{and} \qquad \left| \tau_h-x_h \right| < \left(4C_{d,m}+4\right) \varepsilon'.
\end{equation}
In particular,
\begin{equation}\label{eq_lemma_sigma_h_and_tau_h_2}
\left| \sigma_h - \tau_h \right| < \left(4C_{d,m}+5\right) \varepsilon'.
\end{equation}
\end{lemma}

\begin{proof}
The claim is immediate for $h=1$ because in this case $\sigma_1=\tau_1=x_1=0$. For $h \in \dbl 2,m \dbr$, given \eqref{eq_bound_length_sum_u_h_with_x_h}, we have $n_h=\lfloor x_h n \rfloor$. Thus,
\begin{align*}
|\sigma_h-x_h| \leq \frac{| \lfloor x_h n \rfloor - x_h n| + |1-x_h|}{n-1} \leq \frac{2}{n-1}< \varepsilon',
\end{align*}
where the last inequality comes from \eqref{eq_condition_n}. For $h=m+1$, we have $\sigma_h=1=x_h$ and the first inequality in \eqref{eq_lemma_sigma_h_and_tau_h_1} holds for all $h \in \dbl 1,m+1 \dbr$.

By the definition of $N_h$, given in \eqref{eq_def_N_h_and_n_h}, and \eqref{eq_bound_length_sum_u_h_with_x_h}, we have
\begin{equation}\label{eq_lemma_expression_N_h}
N_h=2k\lfloor x_{h} n \rfloor + \sum_{i=1}^{2k} \sum_{r=1}^{h-1} |\xi^{(i),r}|.
\end{equation}
Using the fact that 
\begin{align*}
\sum_{i=1}^{2k} \sum_{r=1}^{h-1} |\xi^{(i),r}| \leq \sum_{i=1}^{2k} \sum_{r=1}^{m} |\xi^{(i),r}| + |b| =  N-2kn,
\end{align*}
we deduce that
\begin{align*}
\left| \frac{N_h}{N}-x_h \right| & \leq \left| \frac{2k\lfloor x_{h} n \rfloor + \sum_{i=1}^{2k} \sum_{r=1}^{h-1} |\xi^{(i),r}|}{N} - \frac{2knx_h}{N} \right| +  \left| 1-\frac{2kn}{N} \right| x_h  \\
& \leq \frac{2k}{N} + 2 \left| 1- \frac{2kn}{N} \right| \\
& < \left(4C_{d,m}+3\right) \varepsilon',
\end{align*}
where the last inequality follows from \eqref{eq_bound_on_distance_n_N_k}, which can be applied since $n$ and $N$ satisfy Condition~\ref{condition_n_N}. Finally,
\begin{align*}
\left| \tau_h - \frac{N_h}{N} \right| = \frac{N-N_h}{N(N-1)} \leq \frac{1}{N-1}<\frac{1}{n}<\varepsilon',
\end{align*}
and we thereby obtain the second inequality in \eqref{eq_lemma_sigma_h_and_tau_h_1}. Equation~\eqref{eq_lemma_sigma_h_and_tau_h_2} follows directly from the two inequalities in \eqref{eq_lemma_sigma_h_and_tau_h_1} and the triangle inequality.
\end{proof}

\begin{lemma}[Profile stability under stitching]\label{lemma_preservation_of_profile}
Suppose that $0<\varepsilon'<\varepsilon_0'$ and that $n \in \NN$ and $N \in \NN$ satisfy Condition~\ref{condition_n_N}. Then, for all $(w^{(1)},\dots,w^{(2k)}) \in \mathcal{W}_{n,k}$, we have
\begin{equation}\label{eq_lemma_preservation_of_profile}
H_N\big( \psi(w^{(1)},\dots,w^{(2k)}) \, , \, \cdot \, \big) \in B_{\mathscr{H}} \left(\frac{H^{(1)}+H^{(2)}}{2},\varepsilon \right).
\end{equation}
\end{lemma}

\begin{proof}
Fix $(w^{(1)},\dots,w^{(2k)}) \in \mathcal{W}_{n,k}$ and let $\overline{w}=\psi(w^{(1)}, \dots , w^{(2k)})$ and $\overline{H}=(H^{(1)}+H^{(2)})/2$. We have to verify that  for all $t \in [0,1]$,
\begin{equation}\label{proof_lemma_preservation_of_profile_1}
\left\| H_N \big( \overline{w},t \big) - \overline{H}(t) \right\|_{\infty} < \varepsilon.
\end{equation}
We first work at time $\tau_h$. The case $h=1$ is immediate because by Lemma~\ref{lemma_profiles_in_profile_space}, $H_N(\overline{w},\cdot) \in \mathscr{H}^{(d)}$ and therefore
\begin{align*}
H_N(\overline{w},0)=0=\overline{H}(0).
\end{align*}
Fix $h \in \dbl 2,m+1 \dbr$ and $s \in \mathscr{S}_d$. By \eqref{eq_def_profile_map}, the choices of $\sigma_h$ and $\tau_h$ give
\begin{equation}\label{proof_lemma_preservation_of_profile_2}
\big[H_n(w^{(i)},\sigma_h)\big]_s=\frac{C^{(s)}_{n_h}(w^{(i)})}{n-1} \qquad \text{and} \qquad \big[H_N(\overline{w},\tau_h)\big]_s=\frac{C^{(s)}_{N_h}(\overline{w})}{N-1}
\end{equation}
By the construction of $\psi$, given in \eqref{eq_def_psi_stitching_map} and the choice of $N_h$, the stitched path $\overline w$ contains up to time $N_h$, precisely the first $h-1$ blocks of each of the $2k$ input words, together with the connectors inserted between these blocks. Away from the endpoints of the blocks and from the connectors, the increments of $\overline w$ coincide with the corresponding increments of the original words. Hence,
\begin{equation}\label{proof_lemma_preservation_of_profile_3}
C^{(s)}_{N_h}(\overline{w})=\sum_{i=1}^{2k} C_{n_h}^{(s)}\bigl(w^{(i)}\bigr)+\mathcal R_{h,s},
\end{equation}
for some remainder $\mathcal R_{h,s}$. We now bound this remainder. The connectors contribute at most $\sum_{i=1}^{2k}\sum_{r=1}^{h-1} |\xi^{(i),r}|$ increments. Moreover, for each of the $2k(h-1)$ block boundaries involved in the construction, the stitching may remove at most one increment from the original word and create at most one new joining increment. Thus the total contribution of the modified boundaries is bounded by $4k(h-1)$. It follows that
\begin{equation}\label{proof_lemma_preservation_of_profile_4}
|\mathcal R_{h,s}| \leq \sum_{i=1}^{2k}\sum_{r=1}^{h-1}|\xi^{(i),r}| +4k(h-1) \leq N-2nk + 4km.
\end{equation}
Combining \eqref{proof_lemma_preservation_of_profile_2} with \eqref{proof_lemma_preservation_of_profile_3} and \eqref{proof_lemma_preservation_of_profile_4} yields, for every $s \in \mathscr{S}_d$ and $h \in \dbl 2,m+1 \dbr$,
\begin{align}
&\left| \big[H_N(\overline{w},\tau_h)\big]_s - \big[ \overline{H}(\tau_h) \big]_s \right| \label{proof_lemma_preservation_of_profile_5}\\
& \qquad \leq \left| \sum_{i=1}^k \frac{C^{(s)}_{n_h}(w^{(2i-1)})}{N-1} - \frac{1}{2} \big[H^{(1)}(\tau_h)\big]_s \right| + \left| \sum_{i=1}^k \frac{C^{(s)}_{n_h}(w^{(2i)})}{N-1} - \frac{1}{2} \big[H^{(2)}(\tau_h)\big]_s \right| + \frac{|\mathcal{R}_{h,s}|}{N-1}. \notag
\end{align}
Using \eqref{proof_lemma_preservation_of_profile_2} and then the fact that $2nk \leq N$ yields 
\begin{align*}
\left| \sum_{i=1}^k \frac{C_{n_h}^{(s)}(w^{(2i-1)})}{N-1}-\frac{1}{2k} \sum_{i=1}^k \big[H_n \big(w^{(2i-1)},\sigma_h\big)\big]_s \right| & \leq \frac{1}{2} \left| \frac{2k(n-1)}{N-1}-1 \right| \\
& \leq \left| \frac{N-2kn}{N} \right| + \frac{1}{n} \\
& < (2C_{d,m}+2)\varepsilon',
\end{align*}
where the last inequality comes from \eqref{eq_condition_n} and \eqref{eq_bound_on_distance_n_N_k}. Thus,  
\begin{align*}
& \left| \sum_{i=1}^k \frac{C^{(s)}_{n_h}(w^{(2i-1)})}{N-1} - \frac{1}{2}\big[H^{(1)}(\tau_h)\big]_s \right| \\ 
& \qquad  < (2C_{d,m}+2)\varepsilon' + \frac{1}{2k} \sum_{i=1}^k d_{\mathscr{H}}\big(H_n(w^{(2i-1)},\cdot),H^{(1)}\big) + \frac{1}{2k} \sum_{i=1}^k \big\| H^{(1)}(\sigma_h) - H^{(1)}(\tau_h) \big\|_{\infty}.
\end{align*}
By definition of the set of trajectories $\mathcal{W}_{n,k}$, we have $H_n(w^{(2i-1)},\cdot) \in B_{\mathscr{H}}(H^{(1)},\varepsilon')$. Moreover, by the second item of Lemma~\ref{lemma_properties_profile_space}, $H^{(1)}$ is 1-Lipschitz. Thus, using these two properties together with \eqref{eq_lemma_sigma_h_and_tau_h_2} of Lemma~\ref{lemma_sigma_h_and_tau_h}, we obtain
\begin{align}
\left| \sum_{i=1}^k \frac{C^{(s)}_{n_h}(w^{(2i-1)})}{N-1} - \frac{1}{2} \big[H^{(1)}(\tau_h)\big]_s \right| & < \left( 4 C_{d,m} + \frac{9}{2} \right) \varepsilon' + \frac{\varepsilon'}{2} \notag \\
& = \left( 4 C_{d,m} + 5 \right) \varepsilon' \label{proof_lemma_preservation_of_profile_6}
\end{align}
Replacing $2i-1$ by $2i$ and $H^{(1)}$ by $H^{(2)}$, the same argument yields
\begin{equation}\label{proof_lemma_preservation_of_profile_7}
\left| \sum_{i=1}^k \frac{C^{(s)}_{n_h}(w^{(2i)})}{N-1} - \frac{1}{2} \big[H^{(2)}(\tau_h)\big]_s \right| < \left( 4 C_{d,m} + 5 \right) \varepsilon'.
\end{equation}
Finally, using Equation~\eqref{proof_lemma_preservation_of_profile_4}, the fact that $2 \leq N$, that $2nk \leq N$ and then \eqref{eq_bound_on_distance_n_N_k} and \eqref{eq_condition_n}, we obtain
\begin{align}
\frac{|\mathcal{R}_{h,s}|}{N-1} & \leq 2 \frac{N-2nk}{N} + m \frac{4}{n} \notag \\
& < \big(4C_{d,m}+m+2\big)\, \varepsilon'. \label{proof_lemma_preservation_of_profile_8}
\end{align}
Substituting \eqref{proof_lemma_preservation_of_profile_6}, \eqref{proof_lemma_preservation_of_profile_7} and \eqref{proof_lemma_preservation_of_profile_8} into \eqref{proof_lemma_preservation_of_profile_5} and then taking the maximum over $s \in \mathscr{S}_d$ gives
\begin{equation}\label{proof_lemma_preservation_of_profile_9}
\big\| H_N(\overline{w},\tau_h)-\overline{H}(\tau_h) \big\|_{\infty} < \big(12 C_{d,m} + m +12 \big) \varepsilon'.
\end{equation}
Finally, fix $t \in [0,1]$ and choose $h \in \dbl 1,m+1 \dbr$ such that $|t-x_h| \leq 1/(2m)$. By Lemma~\ref{lemma_properties_profile_space}, both $H_N(\overline{w},\cdot)$ and $\overline{H}$ are 1-Lipschitz. Hence, using \eqref{eq_lemma_sigma_h_and_tau_h_1} and \eqref{proof_lemma_preservation_of_profile_9},
\begin{align*}
\big\| H_N(\overline{w},t) - \overline{H}(t) \big\|_{\infty} & \leq \big\| H_N(\overline{w},\tau_h) - \overline{H}(\tau_h) \big\|_{\infty} + 2 |t-\tau_h| \\
& < \big(12 C_{d,m} + m +12 \big) \varepsilon' + \frac{1}{m} + 2\left(4C_{d,m}+4\right) \varepsilon' \\
& < \varepsilon,
\end{align*}
where the last inequality comes from the fact that $\varepsilon'< \varepsilon_0'$ together with Equation~\eqref{eq_condition_varepsilon_eta_3}. Since this holds for all $t \in [0,1]$, the desired result is proved.
\end{proof}

Combining Lemmas~\ref{lemma_property_stitching_map_measure_close} and \ref{lemma_preservation_of_profile} with the definition of $\Omega_N\left(\frac{\mu_1+\mu_2}{2}, \frac{H^{(1)}+H^{(2)}}{2}, \varepsilon \right)$ immediately yields the following corollary.

\begin{corollary}[Target-set inclusion]\label{corollary_stitching_well_defined}
Suppose $0<\varepsilon'<\varepsilon_0'$ and that $(n,N) \in \NN^2$ satisfies Condition~\ref{condition_n_N}. Then, 
\begin{equation}\label{eq_corollary_stitching_well_defined}
\psi(\mathcal{W}_{n,k}) \subseteq \Omega_N \left(\frac{\mu_1+\mu_2}{2}, \frac{H^{(1)}+H^{(2)}}{2},\varepsilon \right).
\end{equation}
\end{corollary}
\noindent This establishes that the stitching procedure satisfies the inclusion requirement~\eqref{eq_objectif_inclusion}.

\subsubsection{Combinatorial bound}\label{subsubsection_combinatorial_bound}

In this section, we show that the stitching procedure $\psi$ satisfies the combinatorial bound \eqref{eq_objectif_cardinality_bound}. We set
\begin{equation}\label{eq_def_D_n_N}
D_{n,N} \eqdef (n+1)^{2km}.
\end{equation}
Note that the right-hand side of \eqref{eq_def_D_n_N} depends on $N$ through $k$.

\begin{lemma}[Fibre bound for the
stitching map]\label{lemma_cardinality_preimage_bound}
Suppose that $0<\varepsilon'<\varepsilon_0'$, that $n \in \NN$ and $N \in \NN$ satisfy Condition~\ref{condition_n_N} and that $\psi$ is constructed as in Definition~\ref{def_psi_stitching_map}. For all $\overline{w} \in \Omega_N$, 
\begin{equation}\label{eq_lemma_cardinality_preimage_bound}
\# \psi^{-1}(\{ \overline{w}\}) \leq D_{n,N}.
\end{equation} 
\end{lemma}

\begin{proof}
If $\overline{w} \in \Omega_N$ is not in the image of $\psi$ then there is nothing to prove. We now suppose that $\psi^{-1}(\{ \overline{w}\}) \neq \emptyset$. The idea is to show that an element $(w^{(1)}, \dots , w^{(2k)}) \in \psi^{-1}(\{\overline{w}\})$ is entirely determined by the length of the subwords $u^{(i),h}$ and $\xi^{(i),h}$. More precisely, consider the map
\begin{align*}
\begin{array}{cccc}
\Phi : & \psi^{-1}(\{\overline{w}\}) & \longrightarrow & \dbl 0,n \dbr^{2km}  \\
& (w^{(1)}, \dots , w^{(2k)}) & \longmapsto & \big(|\xi^{(i),h}| \big)_{\substack{1 \leq i \leq 2k \, , \\ 1 \leq h \leq m}} 
\end{array}.
\end{align*}
Note that for all $(i,h) \in \dbl 1,2k \dbr \times \dbl 1,m \dbr$, we have $|\xi^{(i),h}|<6d\varepsilon'n$. Combining this with \eqref{eq_condition_varepsilon_eta_2} yields
\begin{align*}
|\xi^{(i),h}| < 6d\varepsilon'n < n.
\end{align*}
Thus, $\Phi$ is well defined. 

Let us show that $\Phi$ is injective. Consider $(w^{(1)},\dots,w^{(2k)}),(\widetilde{w}^{(1)},\dots,\widetilde{w}^{(2k)}) \in \psi^{-1}(\{ \overline{w}\})$ such that $\Phi(w^{(1)},\dots,w^{(2k)})=\Phi(\widetilde{w}^{(1)},\dots,\widetilde{w}^{(2k)})$. There exist $b,\widetilde{b} \in \Omega_{\text{fin}}$ such that 
\begin{align}
u^{(1),1}\xi^{(1),1}u^{(2),1} \cdots u^{(2k),m} b & = \psi(w^{(1)},\dots,w^{(2k)}) \label{eq_equality_for_injectivity} \\
& = \overline{w} = \psi(\widetilde{w}^{(1)},\dots,\widetilde{w}^{(2k)}) = \widetilde{u}^{(1),1}\widetilde{\xi}^{(1),1}\widetilde{u}^{(2),1} \cdots \widetilde{u}^{(2k),m} \widetilde{b}. \notag
\end{align}
By construction of the subwords $u^{(i),h}$ and $\widetilde{u}^{(i),h}$ in \eqref{eq_def_subwords_for_profile_partition}, we have for all $i$ and $h$,
\begin{align*}
|u^{(i),h}|=\lfloor x_{h+1} n \rfloor - \lfloor x_h n \rfloor =|\widetilde{u}^{(i),h}|.
\end{align*}
Moreover, since $\Phi(w^{(1)},\dots,w^{(2k)})=\Phi(\widetilde{w}^{(1)},\dots,\widetilde{w}^{(2k)})$, we deduce that for all $i$ and $h$,
\begin{align*}
|\xi^{(i),h}|=|\widetilde{\xi}^{(i),h}|.
\end{align*}
Thus, using \eqref{eq_equality_for_injectivity}, we deduce that we must have $u^{(1),1}=\widetilde{u}^{(1),1}$ and by induction we further deduce that for all $i$ and $h$,
\begin{align*}
u^{(i),h}=\widetilde{u}^{(i),h},
\end{align*}
and therefore $w^{(i)}=\widetilde{w}^{(i)}$ for all $i \in \dbl 1,2k \dbr$. In other words, $\Phi$ is injective. We deduce that
\begin{align*}
\# \psi^{-1}(\{ \overline{w}\}) & \leq \# \dbl 0,n \dbr^{2km}  \\
& \leq (n+1)^{2km}.
\end{align*}
Using the definition of $D_{n,N}$, given in \eqref{eq_def_D_n_N}, we have obtained \eqref{eq_lemma_cardinality_preimage_bound}.
\end{proof}

\subsubsection{Lower decoupling}\label{subsubsection_probability_comparaison}

The aim of this section is to show that the stitching procedure $\psi$ satisfies the lower decoupling property stated in \eqref{eq_objectif_comp_proba}. For all $n,N \in \NN$, we define
\begin{equation}\label{eq_def_C_n_N}
C_{n,N} \eqdef \left(\frac{1}{2d} \right)^{N-2kn+2k}.
\end{equation} 
\begin{lemma}[Lower decoupling]\label{lemma_comparaison_probabilities}
Suppose that $0<\varepsilon'<\varepsilon_0'$ and that $(n,N) \in \NN^2$ satisfies Condition~\ref{condition_n_N}. Fix $(w^{(1)},\dots,w^{(2k)}) \in \mathcal{W}_{n,k}$. For all $\overline{w} \in \Omega_N$ such that $\psi(w^{(1)}, \dots, w^{(2k)})=\overline{w}$, we have
\begin{align}
\PP_{d,\textup{sym}}\left( S_{[1:N]}= \overline{w} \right) \geq  C_{n,N} \prod_{i=1}^{2k} \PP_{d,\textup{sym}} \left( S_{[1:n]}=w^{(i)} \right). \label{eq_lemma_comparaison_probabilities}
\end{align}
\end{lemma}

\begin{proof}
Recall that by Lemma~\ref{lemma_MERW_p_half_is_SRW}, under $\PP_{d,\text{sym}}$, the MERW $(S_n)_{n \in \NN}$ has the same law as a simple random walk. Hence, for all $w \in \Omega_n^{(0)}$, we have
\begin{align*}
\PP_{d,\text{sym}}\big(S_{[1:n]}=w\big)=\left( \frac{1}{2d} \right)^{n-1}.
\end{align*}
In particular, we deduce that
\begin{align*}
\PP_{d,\text{sym}}\big(S_{[1:N]}=\overline{w}\big)&=\left( \frac{1}{2d} \right)^{N-1} \\
& = \left( \frac{1}{2d} \right)^{N-1-2k(n-1)} \times  \prod_{i=1}^{2k} \left( \frac{1}{2d} \right)^{n-1}  \\
& = \left( \frac{1}{2d} \right)^{N-2kn+2k-1} \times \prod_{i=1}^{2k} \PP_{d,\text{sym}}\big(S_{[1:n]}=w^{(i)}\big),
\end{align*}
and Equation~\eqref{eq_lemma_comparaison_probabilities} follows from the definition of $C_{n,N}$.
\end{proof}

\subsection{Proof of Proposition~\ref{prop_good_stitching_procedure}}\label{subsection_existence_RL_function_for_couple}

\begin{proof}[Proof of Proposition~\ref{prop_good_stitching_procedure}]
Fix $\varepsilon>0$ and $\kappa>0$. Choose $m \in \NN$ according to \eqref{eq_cdt_m_sufficiently_large}, then choose $\varepsilon'$ sufficiently small such that $0<\varepsilon'<\varepsilon_0'$. Choose $n_0 \in \NN$ such that \eqref{eq_condition_n} holds for every $n \geq n_0$. For each such $n$, choose $N_0(n) \in \NN$ so that \eqref{eq_condition_N} holds for every $N \geq N_0(n)$.

For fixed $n \geq n_0$ and $N \geq N_0(n)$, let $k$, $\psi$, $C_{n,N}$ and $D_{n,N}$ be defined by \eqref{eq_def_of_k}, Definition~\ref{def_psi_stitching_map}, \eqref{eq_def_C_n_N} and \eqref{eq_def_D_n_N} respectively. Corollary~\ref{corollary_stitching_well_defined}, Lemma~\ref{lemma_cardinality_preimage_bound} and Lemma~\ref{lemma_comparaison_probabilities} establish the inclusion, combinatorial and lower-decoupling properties, respectively \eqref{eq_objectif_inclusion}--\eqref{eq_objectif_comp_proba}. It remains to verify the asymptotic bound \eqref{eq_objectif_asympto}. Using the definitions of $D_{n,N}$ and $C_{n,N}$ given in \eqref{eq_def_D_n_N} and \eqref{eq_def_C_n_N}, we have
\begin{align*}
\frac{1}{N} \log \frac{C_{n,N}}{D_{n,N}}& =\frac{N-2kn+2k}{N} \log \left( \frac{1}{2d} \right) - \frac{2k m}{N} \log(n+1). 
\end{align*}
Using \eqref{eq_bound_on_distance_n_N_k} and the fact that $2nk \leq N$, we have
\begin{align*}
\frac{1}{N} \log \frac{C_{n,N}}{D_{n,N}} \geq  \left((2C_{d,m}+1)\varepsilon' + \frac{1}{n} \right)\log \left( \frac{1}{2d} \right) - m \frac{\log(n+1)}{n}.
\end{align*}
Using \eqref{eq_condition_n} and then \eqref{eq_condition_varepsilon_eta_4}, we deduce that
\begin{align*}
\frac{1}{N} \log \frac{C_{n,N}}{D_{n,N}} \geq -\kappa- m \frac{\log(n+1)}{n}.
\end{align*}
Since $\lim_{n \to \infty} m \log(n+1)/n = 0$, 
\begin{align*}
\liminf_{n \to \infty} \liminf_{N \to \infty} \frac{1}{N} \log \frac{C_{n,N}}{D_{n,N}} \geq -\kappa,
\end{align*}
which is precisely the asymptotic bound \eqref{eq_objectif_asympto}. Therefore, the function $\psi$ satisfies the four properties required by Proposition~\ref{prop_good_stitching_procedure}.
\end{proof}


\appendix
\makeatletter
\renewcommand{\@seccntformat}[1]{\textbf{\csname the#1\endcsname}\quad}
\makeatother
\section{List of notation} \label{section_appendix_list_of_notations}

\renewcommand{\arraystretch}{1.12}

\begin{tabular}{llll}
\toprule
\textbf{Symbol} & \textbf{Description} & \textbf{} & \textbf{$\text{1}^{\text{st}}$ use} \\
\midrule

\notation{$\mathscr{S}_d$}{Steps in $\ZZ^d$}{$\{ \pm e_i \}_{i=1}^d$}{p.~\pageref{eq_def_mathscr_S_d}}

\notation{$U_n$, $A_n$}{Construction of the MERW}{}{p.~\pageref{eq_def_U_n}}

\notation{$S_n$, $\PP_{d,p}$}{MERW}{steps : $X_n=A_n X_{U_{n-1}}$}{p.~\pageref{eq_def_ERW}}

\notation{$p_d$}{Diffusive critical parameter}{$p_d=\frac{2d+1}{4d}$}{\eqref{eq_def_diffusive_to_superdiff_crit_param}}

\notation{$\ell_n$}{Empirical measure}{$\ell_n=\frac{1}{n}\sum \limits_{j=1}^n \delta_{S_j}$}{\eqref{eq_def_empirical_measure}}

\notation{$\mathcal{P}(\ZZ^d)$}{Probability measures on $\ZZ^d$}{}{p.~\pageref{eq_def_empirical_measure}}

\notation{$I_{d,p}$}{Weak LDP rate function}{}{\eqref{expression_rate_function}}

\notation{$I_{d,\textup{DV}}$}{DV rate function}{}{\eqref{Donsker_Varadhan_our_case}}

\notation{$p_d^*$}{LDP critical parameter}{$p>p_d^* \Leftrightarrow$ $I_{d,p}$ non-convex}{\eqref{eq_def_new_critical_point}}

\notation{$\| \mu \|$}{KR norm}{}{p.~\pageref{eq_def_BL_norm}}

\notation{$\Omega_n,\Omega_n^{(0)}$}{Trajectories of length $n$}{}{\eqref{def_Omega_n}}

\notation{$\Omega_n(\mu,\varepsilon)$}{Traj.~of length $n$ close to $\mu$}{$w_1=0$, $\ell_n(w) \simeq \mu$}{\eqref{eq_def_Omega_mu_varepsilon}}

\notation{$N_A(w)$}{Occupation time in $A \subseteq \ZZ^d$}{$N_A(w)=\sum_{j=1}^{|w|} \mathbbm{1}_A(w_j)$}{\eqref{eq_def_occupation_time_N_A}}

\notation{$\Delta^{(d)}$}{$(2d-1)$ dimensional simplex}{}{p.~\pageref{eq_def_dimplex}}

\notation{$\mathscr{H}^{(d)}$}{Space of cumulative profiles}{}{\eqref{eq_def_H_d_space_of_profiles}}

\notation{$d_{\mathscr{H}}$}{Uniform distance on profiles}{$\sup_{t \in [0,1]} \| H(t)-\widehat{H}(t) \|_{\infty}$}{\eqref{eq_def_uniform_distance_for_profiles}}

\notation{$\Delta_j(w)$}{$j$-th step of $w$}{$\Delta_j(w) = w_{j+1}-w_j$}{\eqref{eq_Delta_word_j}}

\notation{$C^{(s)}_k(w)$}{Step counter}{$\#\{ j < k:w_{j+1}-w_j=s\}$}{\eqref{eq_def_step_counter}}

\notation{$H_n(w,\cdot)$}{Cumulative profile}{$\left[H_n\left(w,\frac{j-1}{n-1}\right)\right]_s=\frac{C^{(s)}_j(w)}{n-1}$}{\eqref{eq_def_profile_map}}

\notation{$\rho_j(w)$}{Distrib. of increments at $j$}{$[\rho_j(w)]_s=C_j^{(s)}(w)/(j-1)$}{\eqref{eq_def_empirical_distribution_rho}}

\notation{$Q_p$}{Transition probability}{$Q_p(s \,|\, r) = p r_s + \frac{1-p}{2d-1}(1-r_s)$}{\eqref{eq_def_Q_p_r_s}}

\notation{$p_*$}{Transition lower bound}{$p_* = \min \left\{ \frac{1-p}{2d-1} \,  , \, p \right\}$}{\eqref{eq_lower_bound_transition_probability}}

\notation{$\Omega_n(\mu, H,\varepsilon)$}{Traj.~close to $\mu$ and $H$}{$w_1=0$, $\ell_n(w) \simeq \mu$, $H_n(w) \simeq H$}{\eqref{eq_def_Omega_n_mu_f}}

\notation{$\overline{s}_p(\mu,H)$, $\underline{s}_p(\mu,H)$}{RL functions}{}{Def.~\ref{definition_RL}}

\notation{$I_{\textup{A}}^{(p)}$}{Rate function for $(\ell_n,H_n)$}{$I_{\textup{A}}^{(p)}(\mu,H)=-\overline{s}_p(\mu,H)$}{\eqref{eq_def_I_RL_couple}}

\notation{$\mathcal{E}_p$}{Energy-profile map}{}{\eqref{eq_def_energy_profile}}

\notation{$L_p$, $F_p$ and $G_p$}{Scalar energy functions}{}{\eqref{eq_def_L_p}}

\notation{$J_{d,p}$}{Profile rate function}{}{\eqref{eq_corollary_LDP_profile_position_1}}

\notation{$\mathscr{H}^{(d)}_{\textup{bal}}$}{Balanced profiles}{$[H(t)]_s=[H(t)]_{-s}, \; \forall t, \forall s$}{\eqref{eq_def_H_bal}}

\notation{$\mu_1,\mu_2,\mu^*$}{}{}{p.~\pageref{eq_def_mu_1}, \pageref{eq_def_mu_star}}

\notation{$\mathscr{H}^*$}{}{$H \in \mathscr{H}^{(d)}_{\text{bal}}$ and $H(1)=x_{\text{end}}$}{\eqref{eq_def_H_star}}

\notation{$\overline{x}(m)$}{Mesh of $[0,1]$}{$x_h=(h-1)/m$}{\eqref{eq_def_overline_x_m}}

\notation{$u^h$, $\xi^h$, $b$}{Subwords for the stitching}{}{p.~\pageref{eq_def_subwords_for_profile_partition}-\pageref{def_psi_stitching_map}}

\notation{$C_{d,m}$}{}{$C_{d,m}=6d(m+1)$}{\eqref{eq_def_C_d_m}}

\notation{$\psi$}{Stitching procedure}{$u^{(1),1}\xi^{(1),1}u^{(2),1} \cdots u^{(2k),m} b$}{Def.~\ref{def_psi_stitching_map}}

\bottomrule
\end{tabular}

\end{document}